\documentclass[11pt]{elsart-maktor}
\renewcommand{\and}{et}
\usepackage{amssymb,
amsfonts,amssymb,euscript,mathrsfs,pifont
}
\usepackage[english,francais]{babel}
\usepackage{amsxtra}
\usepackage{amscd}
\usepackage[T1]{fontenc}                
\usepackage[all]{xy}

 \DeclareMathAlphabet{\got}{U}{euf}{m}{n}     
\DeclareMathAlphabet{\mat}{U}{msb}{m}{n}     
\DeclareMathAlphabet{\mathbold}{OML}{cmm}{bx}{it} 

\makeatletter

\renewcommand{\section}{\@startsection{section}{1}%
\z@{.5ex \@plus.7ex}{-.5em}%
{\normalfont\bfseries}}

\renewcommand{\subsection}{\@startsection{subsection}{2}%
\z@{.5ex \@plus.7ex}{-.5em}%
{\normalfont\bfseries}}

\renewcommand{\subsubsection}{\@startsection{subsubsection}{3}%
\z@{.5ex\@plus.7ex}{-.5em}%
  {\normalfont\itshape}}
\makeatother
\numberwithin{equation}{section}

\newtheorem{tho}{Th\'eor\`eme}[subsection]
\newtheorem{pro}{Proposition}[subsection]
\newtheorem{lem}{Lemme}[subsection]
\theoremstyle{definition}

\newtheorem{cor}{Corollaire}[subsection]

\newcommand{\goth}{\mathfrak}

\newenvironment{theo*}{\vskip 1em\bf Th\'eor\`eme.\it}{\par\rm}
\newenvironment{dem}{\vskip 1em{\it D\'emonstration} :}%
{\unskip\hfill\null\nobreak\hfill\carre\vskip1em\par}
\newcommand{\carre}{\rule{1ex}{1ex}} 

\renewcommand{\thefootnote}{\fnsymbol{footnote}}

\makeatletter
\providecommand*{\diff}%
{\@ifnextchar^{\DIfF}{\DIfF^{}}}
\def\DIfF^#1{%
\mathop{\mathrm{\mathstrut d}}%
\nolimits^{#1}\gobblespace}
\def\gobblespace{%
\futurelet\diffarg\opspace}
\def\opspace{%
\let\DiffSpace\!%
\ifx\diffarg(%
\let\DiffSpace\relax
\else
\ifx\diffarg[%
\let\DiffSpace\relax
\else
\ifx\diffarg\{%
\let\DiffSpace\relax
\fi\fi\fi\DiffSpace}

\providecommand*{\eu}%
{\ensuremath{\mathrm{e}}}

\providecommand*{\iu}%
{\ensuremath{\mathrm{i}}}

\def \cit { \mathbb C}
  
\def \reel  {\mathbb R}
\def \N {\mathbb N}
\def \Z {\mathbb Z }
\def\G{ \bf G }
          \def\var{\varphi}
\def\U {\bf U}          \def \b{\goth b}
\def\g{ \goth g}                  
\def\u{ \goth u}         \def \j{\goth j}
\def\l{ \ell}           \def \h{\goth h}
    \def \‡{\varsigma}      
            \def \O { \mathcal O}
 \def \W { \mathcal W}

\def \O{\mathcal O}      
\def\var{\varphi}      \def \q{\goth q}
    \def \a{\goth a}
         
                           \def \r{\goth r}
       
           \def \t{\goth t}
     \def \e { \varepsilon }
    \def\m{\goth m} \def \dis{\displaystyle}
\begin{document}

\begin{frontmatter}
\title{La formule du  caract\`ere et la mesure de Plancherel
 pour les groupes de Lie r\'esolubles unimodulaires  sur un corps $p$-adique.}

\author{Khemais Maktouf}
\address{Université de Monastir, 
Faculté des Sciences de Monastir, Département de Mathématiques,
 5019 Monastir, Tunisie}
\ead{khemais.maktouf@fsm.rnu.tn}

\selectlanguage{english}
\begin{abstract}
We establish a  character formula for admissible unitary representations of $p$-adic almost algebraic solvable groups and we deduce the Plancherel
 measure in the  unimodular case.



\medskip
\noindent{\bf R\'esum\'e}
 \medskip
 
\noindent Nous \'etablissons une formule du  caract\`ere pour
les repr\'esen\-tations unitaires admissibles des groupes r\'esolubles
presque alg\'ebriques sur un corps  $p$-adique et nous en 
 d\'eduisons la mesure de Plancherel dans le cas unimodulaire. 

\noindent{\it Mots-cl\'es:} Groupe m\'etaplectique, Repr\'esentation unitaire, Formule du caract\`ere, M\'ethode de  descente, Mesure  Plancherel, Groupe r\'esoluble


\end{abstract}
\end{frontmatter}
\selectlanguage{french}

\begin{section}{Introduction.}\label{introsc}
Le but de ce travail est de donner une d\'emonstration de la formule de Plancherel pour 
une large classe de groupes de Lie r\'esolubles unimodulaires qui s'inscrit dans le cadre de la m\'ethode des orbites de Kirillov-Duflo. Pour ce faire, nous donnons une description
globale des caract\`eres des repr\'esentations unitaires irr\'eductibles g\'en\'eriques~: 
en fait nous \'etablissons au voisinage de chaque \'el\'ement semi-simple une formule du
 caract\`ere. Notre formule s'inspire fortement de celles \'etablies dans le cas r\'eel par Duflo-Heckman-Vergne \cite{[DHV]}  (s\'eries discr\`etes des groupes r\'eductifs), Bouaziz \cite{[Boua]} (repr\'esentations temp\'er\'ees des groupes r\'eductifs), Khalgui-Torasso 
\cite{[KT]} (groupes presque alg\'ebriques r\'eels). Notre d\'emonstration de la formule 
de Plancherel s'inspire \'egalement de celles mises en \oe{}uvre dans le cas r\'eel par
 Duflo-Vergne \cite{[Du-Ve2]} pour les groupes r\'eductifs et par Khalgui-Torasso \cite{[KT2]} pour les groupes presque alg\'ebriques r\'eels.

Soient $k$  un corps local non archim\'edien de
caract\'eristique z\'ero et $\varsigma$  un caract\`ere
  non trivial de  $(k, +)$.  
On consid\`ere $(G, F, {\bf G})$  un groupe r\'esoluble  presque alg\'ebrique sur $k$ (voir le num\'ero \ref{DJ} ci-apr\`es), on note $\g$ son
alg\`ebre de Lie et $\g^* $ le dual de $\g$.
 Dans \cite{[Du]}, M.
Duflo a donn\'e une para\-m\'etri\-sation du dual unitaire de $G$,
l'ensemble des classes des repr\'esentations unitaires
irr\'eductibles de $G$~: pour chaque forme lin\'eaire $g$ sur
$\g$, on note $G(g)$ le stabilisateur de $g$ dans $G$ et $\g(g)$ son
alg\`ebre de Lie~; on d\'esigne par  $Mp(\g/\g(g))$  le groupe
m\'etaplectique de Weil et on consid\`ere $G(g)^\g$ le
produit fibr\'e de $G(g)$ par $Mp(\g/\g(g))$, c'est un rev\^etement d'ordre deux de $G(g)$ dont l'\'el\'ement 
non trivial du noyau est not\'e $ (1,-1)$. On d\'efinit un
caract\`ere $\chi_g$ de ${}^uG(g)^\goth g$, radical unipotent de $G(g)^\g$,
en posant, $$\chi_g(\exp(X))= \varsigma(<g, X>), \, X \in
{}^u\g(g).$$
 On note $Y_G(g)$ l'ensemble des classes  des
repr\'esentations unitaires irr\'eductibles $\tau$ de $G(g)^\g$
dont la restriction \`a ${}^uG(g)^\goth g$ est multiple de
$\chi_g$ et telle que $\tau (1,-1)=-\mbox{Id}$. 
Si $g$ est de type
unipotent (i.e. $g$  est nulle sur un facteur
r\'eductif de $\g(g)$) et $\tau$ un \'el\'ement de $Y_G(g)$, M. Duflo \cite{[Du]} a
associ\'e au couple $(g, \tau)$ une classe de repr\'esentations
unitaires irr\'eductibles de $G$, not\'ee $\pi_{g , \tau}$.
Si on d\'esigne par $Y_G$ l'ensemble des couples $(g, \tau)$, o\`u
$g$ est de type unipotent et $\tau \in Y_G(g)$, le groupe $G$
op\`ere naturellement dans $Y_G$ et la corres\-pondance $(g, \tau) \longmapsto 
\pi_{g, \tau}$
induit une bijection de $G\backslash Y_G$ sur le dual unitaire de
$G$.

Supposons d\'esormais que  $G$ est r\'esoluble  presque connexe, c'est-\`a-dire
$G/(F  .\, {}^uG)$ est ab\'elien,  ${}^uG$ \'etant le radical unipotent de $G$.  Nous montrons qu'il existe un voisinage $\mathcal U$
de $0$ dans $\g$ tel que, pour tout $g \in \g^*$ de type unipotent et $\tau \in Y_G(g)$,
   il existe  une forme
lin\'eaire $\lambda_\tau$ sur $\g(g)$ dont la restriction \`a ${}^u\g(g)$ est \'egale \`a celle de $g$ et tel que, si $X \in \mathcal U \cap \g(g)$, on ait 
\begin{eqnarray}\label{duflo}
\tau(\exp X)&=& \varsigma(<\lambda_\tau, X>)\mbox{ Id }.
\end{eqnarray}
Nous  consid\'erons ainsi
l'orbite co-adjointe $\O_{g,\lambda_\tau }$ construite \`a partir
d'un \'el\'ement $f$ qui a m\^eme restriction que $g$ au radical
unipotent de $\g$ et dont la restriction \`a $\g(g)$ est
$\lambda_\tau$.

Soit $s$ un \'el\'ement semi-simple de $G$, $G(s)$ son centralisateur dans $G$, 
et $\g(s)$ son alg\`ebre de Lie. On note  $\O_{g,\lambda_\tau, s}$ l'ensemble des points fixes de $s$ dans  $\O_{g,\lambda_\tau}$.  Lorsqu'il n'est pas vide, 
 $\O_{g,\lambda_\tau, s}$ est une sous-vari\'et\'e localement ferm\'ee 
de $\g(s)^*$, 
r\'eunion d'un nombre 
fini de $G(s)$-orbites. On munit  $\O_{g,\lambda_\tau, s}$ 
de la  mesure canonique (de Liouville) $d\mu_{\O_{g,\lambda_\tau,
s}}$, celle  dont la restriction \`a chaque $G(s)$-orbite qu'il contient
est la mesure de Liouville.
Nous
d\'efinissons une fonction $\phi_{g, \tau,s}$ sur $\O_{g,\lambda_\tau, s}$ :
si $l= x.f \in \O_{g,\lambda_\tau, s}$ tel que $(1-s)$ soit un endomorphisme   inversible de $\g/\g(l)$, on a,   
\begin{eqnarray}\label{thall}
\phi_{g, \tau,s}(l) &=& \Phi(\rho_x(\tilde s))\mbox{ Tr} (\, {}^x \tau(\tilde s)),
 \end{eqnarray}
o\`u $\tilde s $ est un relev\'e de $s$ dans $G(l)^\g$, 
$\rho_x(\tilde s)$ est l'image de $\tilde s $ dans $Mp(\g/\g(l))$, et 
$ \Phi$ est le caract\`ere de la repr\'esentation m\'etaplectique de
  $Mp(\g/\g(l))$ associ\'ee \`a $\varsigma$ (le membre de droite de (\ref{thall}) ne d\'epend ni du choix de $x \in G$ tel que $l = x.f$, ni du choix de $\tilde s $ situ\'e au-dessus de $s$). La fonction $\phi_{g, \tau,s}$ est  $G(s)$-invariante. Voir le num\'ero \ref{hcgg}.

Depuis \cite{[M3]}, 
on sait que si l'orbite co-adjointe $\O_g= G.g$ est ferm\'ee dans $\g^*$ alors
 $\pi_{g , \tau}$ est admissible
(\`a trace).  On note $\Theta_{g, \tau}$ son caract\`ere,  c'est
la fonction g\'en\'eralis\'ee, $G$-invariante,
sur $G$, d\'efinie par
$$
\Theta_{g, \tau}(\var d_Gx) = \mbox{Tr}(\pi_{g, \tau}(\var
d_Gx)), \,\,\,  \var \in C_c^\infty(G),
$$
o\`u $d_Gx$ est une mesure de Haar sur $G$.

 Si  $\mathcal V$ est un voisinage ouvert $G(s)$-invariant assez petit de $1$ dans $G(s)$, l'image   $\mathcal W$ de 
$G \times  \mathcal V$ par l'application $\Psi : (x, y) \longmapsto
 xsy x^{-1}$ est un voisinage ouvert $G$-invariant de $s$ dans $G$ et 
$\Psi$
 induit 
un diff\'eomorphisme 
 de l'espace fibr\'e $G \times_{G(s)}  \mathcal V$ sur 
$\mathcal W$. Par  suite, suivant la m\'ethode de descente, 
$\Theta_{g, \tau}$ poss\`ede une restriction $\theta_{g, \tau,s}$ \`a~$\mathcal V$, symboliquement d\'efinie par :
$$
\quad \theta_{g, \tau,s}(y) =\Theta_{g, \tau}( s y), \,\, \,  y \in 
 \mathcal V. 
$$
Le r\' esultat principal de la premi\`ere partie du pr\' esent
travail est le suivant.
\begin{tho} Soit $s$ un \'el\'ement semi-simple de $G$. Il existe 
    un voisinage  ouvert $\mathcal V_s$, 
 $G(s)$-invariant,  de $1$ dans $G(s)$ tel que 
pour tout 
$(g, \tau) \in Y_G$,  pourvu que   l'orbite co-adjointe $\O_{g}$ soit ferm\'ee dans $\g^*$,  
l'on ait, pour tout $   \beta \in C_c^\infty(\mathcal V_s)$, 
\begin{eqnarray}\label{x}  \theta_{g,\tau,s}(\beta d_{G(s)}x)&=& \int_{
\O_{g,\lambda_\tau, s}} \widehat{(\beta \circ \exp
d_{\g(s)}X)}_{\g(s)} \, (l)\phi_{g,\lambda_\tau,s}(l)
d\mu_{\O_{g,\lambda_\tau, s}}(l),  
\end{eqnarray}
 o\`u $\lambda_\tau$ est une forme lin\'eaire sur $\g(g)$ v\'erifiant la condition (\ref{duflo}),  $d_{G(s)}x$ et $d_{\g(s)}X$ sont  des mesures de Haar sur $G(s)$ et $\g(s)$ respectivement
qui se correspondent.
\end{tho}

\medskip

Il convient de rappeler ici que dans \cite{[M3]}, nous avons d\'emontr\'e pour les
 groupes de Lie \`a radical co-compact  une
 formule du caract\`ere au voisinage d'\'el\'ements semi-simples  
 en termes 
de transform\'ees de Fourier de l'ensemble de points fixes de $s$ dans $\O_g$,
 d'une fonction bien d\'efinie sur ce dernier, et du caract\`ere de la
 repr\'esentation $\tau$.  Cela dit, l'ouvert  de  validit\'e  de la formule du caract\`ere   d\'epend de la repr\'esentation $\tau$. De ce fait, cette formule est insuffisante pour d\'emontrer la formule de Plancherel de $G$.

Dans la deuxi\`eme partie, comme  application   de la formule (\ref{x}), nous donnons, dans le cas o\`u $G$ est unimodulaire,  la
mesure de  Plancherel de $G$.
Pour cela,   nous  \'etudions  les formes lin\' eaires
fortement r\' eguli\`eres sur l'alg\`ebre de Lie $\g$. Un \' el\'
ement $g \in \g^*$ est dit   r\' egulier si $\dim \g(g)$ est
minimale. Alors $\g(g)$ est commutative, on note $\j_g$ l'unique
facteur r\' eductif de $\g(g)$. On dit que $g$ est fortement
r\'egulier si de plus $\j_g$ est de dimension maximale.
 L'ensemble des \' el\' ements
 fortement r\'eguliers de $\g^*$, not\'e $\g^*_{t.r}$, 
 est un ouvert de Zariski, $G$-invariant, non vide de $\g^*$.
On montre qu'il existe un  ouvert de Zariski $\Omega_G$, $G$-invariant, non vide de
$\g^*$, contenu dans $\g^*_{t.r}$ et invariant par translation par les \'el\'ements de $({}^u\g)^\perp$ tel que, pour tout $g \in \Omega_G$, l'orbite co-adjointe $G.g$ est ferm\'ee dans $\g^*$. 
Notons $\Omega_U$, image de $\Omega_G$ par l'application restriction de $\g^*$ \`a $({}^u\g)^*$.   Alors,  $\Omega_U$ est  un ouvert de Zariski, non vide, $G$-invariant de $({}^u\g)^*$. 

\noindent Les sous-alg\`ebres de Lie $\j_g$, $g \in \g^*_{t.r}$,  sont
appel\'ees  sous-alg\`ebres de Cartan-Duflo.  Toutes les
sous-alg\`ebres de Cartan-Duflo dans $\g$ sont $G$-conjugu\' ees 
(Proposition \ref{yreg3}). Fixons un repr\'esentant $\j$ de cette
classe que l'on munit d'une  mesure de Haar $d_\j X$. On se donne une  mesure de Haar $d_\g X$ (resp. $d_{{}^u\g}X$) sur $\g$ (resp. $ {}^u\g$).  
Par la suite, nous montrons que, pour chaque $u \in  \Omega_U$, l'orbite $G.u$ admet une mesure positive $G$-invariant $d\beta_{G.u}$ (cette mesure d\'epend de  $d_\g X$, $d_{{}^u\g}X$, et $d_\j X$). Ainsi, il existe une unique mesure bor\'elienne positive $\mu_u$  sur
 $G\backslash ({}^u\g)^*$
 telle que l'on ait,
 \begin{eqnarray}\label{thal}
 \int_{({}^u\g)^*}\psi(l) d_{({}^u\g)^*}l &=& \int_{G\backslash ({}^u\g)^*}d\mu_u(\omega)\int_{\omega}\psi(l)d\beta_\omega(l)
  \end{eqnarray}
  pour toute fonction $\psi$,  bor\'elienne positive, ou int\'egrable  sur $({}^u\g)^*$.
 
Pour  $u \in \Omega_U$, soit $ g \in \Omega_{G}$ de type unipotent dont la restriction \`a ${}^u\g$ est $u$. Notons  $\j_g$  l'unique facteur r\'eductif de $\g(g)$. Il existe alors
$x \in G$ tel que $\j_g= x.\j$. On munit ainsi  $\j_g$ de la mesure de Haar $d_{\j_g}X$, image de la mesure de Haar $d_\j X$ par  l'application $ \j \longrightarrow \j_g, X \longmapsto x.X$.  
Soit  $R_g$ un facteur r\'eductif de
 $G(g)$.
On d\'esigne par $d_{ R_g^\g}x$ la mesure de Haar sur $ R_g^\g$
tangente \`a  la mesure de Haar $d_{\j_g} X$ sur $\j_g$ et par  $d_{ \widehat{(R_g^\g)}}\tau$
la mesure de Plancherel de ${R_g^\g}$ correspondante.  On note $\widehat{(R_g^\g)}_{-}$
l'ensemble des classes  des repr\'esentations unitaires
irr\'eductibles $\tau$ de $ R_g^\g$  telles que $\tau (1,-1)= -\mbox{Id}$.
  On note  $d_g\tau$ la mesure image de  $d_{ \widehat{(R_g^\g)}}\tau$ sur $Y_G(g)$ par l'application   $ \widehat{(R_g^\g)}_{-}\longrightarrow  Y_G(g), \, \,  \tau \longmapsto \tau \otimes \chi_g$ (voir le num\'ero \ref{plancar}).
On d\'efinit une fonction g\'en\'eralis\'ee, $G$-invariante,   $\Theta_{G.u}$ sur $G$ par~: 
\begin{eqnarray*} 
\Theta_{G.u}(\var d_Gx) = 2\int_{Y_G(g)}  \Theta_{g, \tau}(\var d_Gx)\, d_g\tau,\quad \var
\in C_c^\infty(G),
 \end{eqnarray*}
  o\`u  $d_Gx$ est  la mesure de Haar sur $G$
tangente \`a la mesure de Haar $d_\g X$ sur $\g$.
Le r\'esultat principal de la deuxi\`eme partie du pr\'esent travail est le suivant.
\begin{tho} Pour tout $\var
\in C_c^\infty(G)$, on a~:
\begin{eqnarray}\label{thale}
\var(1) &=& \int_{ G \backslash ({}^u\g)^*}\Theta_\omega(\var d_Gx)
d\mu_u(\omega ).
 \end{eqnarray}
\end{tho}
\noindent Pour \'etablir la formule (\ref{thale}), nous devons  d\'emontrer le r\'esultat suivant : 
 \begin{eqnarray} \label{tha}
 \Theta_{G.u}(\var  d_Gx)&=& \int_{G.u}\widehat{(\var_{|{}^uG}\circ \exp d_{{}^u\g}X)}_{{}^u\g}(l) d\beta_{G.u}(l), \, \, \forall \, \var \in C_c^\infty(G). 
 \end{eqnarray}
Pour ce faire,  la formule (\ref{tha}) s'interpr\`ete  comme une \'egalit\'e entre fonctions g\'en\'eralis\'ees $G$-invariantes.  En utilisant la formule (\ref{x}), nous d\'emontrons 
 qu'elle est satisfaite sur les
voisinages semi-simples de chaque \'el\'ement semi-simple de $G$.
  Les formules (\ref{thal}) et (\ref {tha}) donnent alors  la formule (\ref{thale}).
\end{section}
\begin{section}{Notations et g\'en\'eralit\'es}\label{aick1}
\begin{subsection}{}
Dans la suite, $k$  d\'esigne   un corps local non archim\'edien de
caract\'eristique z\'ero, $\O$  l'anneau des entiers, et
$\varpi$ une uniformisante de $\O$.  Le  corps r\'esiduel $\O/\varpi \O$ est fini  de cardinal $ q $ et de 
   caract\'eristique $p$. 
On note  $v$ la valuation  normalis\'ee de $k$~: $v(x) = \sup\{n \in \Z \mbox{ tel que } x \in
\varpi^n\O\}$. On
d\'efinit une valeur absolue sur $k$~:
$|x|_p = q^{-v(x)}, x \in k$.
On fixe  une cl\^oture alg\'ebrique
$\bar k$ de $k$ et on prolonge la norme $p$-adique $| \, . \, |_p$ de $k$ \`a $\bar
k$ de la mani\`ere suivante : si $x \in {\bar k}^\times$, on
consid\`ere un corps $L \subset \bar k$ contenant $k$ et $x$ et
tel que l'on ait $[L:k] < \infty$. On pose $$|x|_{\bar k}
= |N_{L/k}(x)|_p^{\frac{1}{[L:k]}},$$ o\`u $N_{L/k}$ d\'esigne l'application norme
de $L$ sur $k$.

\noindent On fixe un caract\`ere
$\varsigma$   de $k$ dans $\cit^\times$ tel que $$
\varsigma_{|\O}= 1  \quad  \mbox{ et } \quad  \varsigma_{|\varpi^{-1}\O}\neq 1 . $$
Pour $a \in k$, on note $\varsigma_a$ le caract\`ere de $k$ dans $\cit^\times$ d\'efini par : $\varsigma_a(t)= \varsigma(at)$, $t \in k$.
 D'apr\`es \cite{[We2]}, l'application $a \in k \longmapsto \varsigma_a $ est un isomorphisme du groupe $(k, +)$ sur le groupe de ses caract\`eres.
Enfin, on d\'esigne par $d\mu$ la mesure de Haar sur $k$ telle que $d\mu(\O)= 1$.

\end{subsection}

\begin{subsection}{}
Si $X$ est un espace  topologique totalement discontinu
(i.e. chaque point de $X$ admet une base
de voisinages form\'ee par des parties compactes ouvertes) on note
 $C_c^\infty(X) $ le $\cit$-espace des fonctions
localement constantes \`a support compacts. Si $
A $ est une partie de $X$, on d\'esigne par $1_A$  sa fonction caract\'eristique.
\end{subsection}
\begin{subsection}{}\label{ddf} Soit $V$ un espace vectoriel de
dimension finie sur $k$. On d\'esigne par $V^*$ l'espace  des formes lin\'eaires sur
$V$. 
 Un r\'eseau de $V$ est un  $\mathcal O $-module, compact et ouvert.
On sait qu'un    $\mathcal O $-module   $L$  est un r\'eseau de $ V $ si et seulement s'il existe une base $(e_1 ,  \ldots,
 e_n)$  de $V$ telle que  $L = \O e_1 +  \cdots + \O e_n$.
 On  note
 $L^\perp = \{l \in V^* \, /\, \varsigma(<l, v>) =1, \mbox{ pour tout } v \in L\}$
 le r\'eseau de $V^*$  que l'on  appelle r\'eseau dual de $L$ relativement au
  caract\`ere $\varsigma$.

 Si $\varphi \in C_c^\infty(V)
$ et $ d_{V}v$  est une mesure de Haar sur $V$, on note $\widehat{(\varphi d_{V}v)}_V$ ou
 $(\varphi d_{V})\, \widehat{}_V$ la transform\'ee de Fourier de la densit\'e
$\var d_{V}v$  relativement \`a $\varsigma$, c'est l'\'el\'ement de
 $C_c^\infty(V^*)$
  d\'efini par~: $$\widehat{(\varphi d_{V}v)}_V(l)
=\int_V\varphi(v)\varsigma(< l , v>)  d_{V}v \mbox {, pour tout } l
\in V^* .$$
La base $(e_1 ,  \ldots,
 e_n)$  de $V$ est dite univolumique pour la mesure de Haar $d_Vv
 $  si l'on a $d_Vv(L) =~1 $.  Si tel est le cas,  on a : 
$$\widehat{(1_L d_{V}v)}_V =1_{L^\perp}. $$
  On appelle mesure duale (relativement \`a $\varsigma$) de la
mesure de Haar  $d_{V}v$, la mesure de Haar $d_{V^*}l$ sur $V^*$ telle que
l'on ait~:
 $$ \var(0) =
\int_{V^*}\widehat{(\varphi d_Vv)}_V (l)  d_{V^*}l \mbox{, pour
tout }
   \var \in C_c^\infty(V). $$
 Soit $W$ un
sous-espace vectoriel de $V$ et $d_{W}w$ une mesure de Haar sur
$W$. On appelle mesure quotient de $d_{V}v$ par $d_{W}w$, la
mesure de Haar $d_{V/W}\dot v$ sur $V/W$ donn\'ee par :
 $$ \int_V \var (v) d_{V}v = \int_{V/W} \int_W \var(v+w) d_{W}wd_{V/W}\dot
 v \mbox{ , pour tout } \var \in C_c^\infty(V).
 $$
D'autre part $W^\perp =
\{l \in V^*, \, l_{|W}= 0 \}$  s'identifie naturellement \`a $(V/W)^*$. On munit alors $W^\perp$ de la mesure
de Haar $d_{W^\perp}\lambda$ duale de la mesure de Haar $d_{V/W}\dot v$. On a,
pour tout $l \in W^*$,
 $$ \widehat{(\var_{|W} d_{W}w)}_W(l) = \int_{W^\perp}\widehat{(\var d_{V}v)}_V(\tilde{l} +
 \lambda)d_{W^\perp}\lambda,\quad  \var \in C_c^\infty(V),
 $$
o\`u $\tilde{l}$ est un \'el\'ement de $V^*$ dont la restriction
\`a $W$ est $l$.

\noindent Soit $Q$  une forme quadratique  non d\'eg\'en\'er\'ee  sur $V$. Il existe un nombre
complexe $\gamma(Q)$ de module $1$ satisfaisant l'\'equation
suivante~:
\begin{eqnarray}\label{quadratique}
\quad \int_V\int_V\varphi(v - w)\varsigma(\frac{1}
{2}Q(w))d_Vwd_Vv &=& c_Q \gamma(Q) \int_V\varphi(v) d_Vv, \quad
\var \in C_c^\infty(V),
\end{eqnarray}
o\`u $c_Q$ est une constante strictement positive (voir \cite{[We1]}).
\end{subsection}
\begin{subsection}{} On appelle $k$-espace symplectique  tout
couple  $(V , B)$  o\`u $ V$ est un $k$-espace vectoriel de
dimension finie et $B$ une forme bilin\'eaire altern\'ee non
d\'eg\'en\'er\'ee sur $V$. Une base
 $(e_1,
\ldots, e_n, f_1, \ldots, f_n)$ de $V$ est dite  symplectique
si $$B(e_i, e_j) = B(f_i, f_j) = 0 \mbox{ et } B(e_i,f_j) =
\delta_{ij}, \, \, \mbox{ pour tout } 1\leq i, j \leq n, $$
$\delta_{ij}$ \'etant le symbole de  Kronecker. Si on  identifie $V$ et  son dual $V^*$ par
 $B$ alors le r\'eseau $r~=  \oplus_{i =1}^{n}\O e_i
\oplus \oplus_{i =1}^{n}\O f_j $ est auto\-duale (i.e. $r^\perp= r$) et 
 la mesure de
Haar sur $V$ telle que $r$ soit de mesure $1$ est la
mesure de Haar  auto\-duale, relativement \`a $\varsigma$.

\noindent On note $Sp(V,B)$ ou plus simplement $Sp(V)$ le groupe symplectique associ\'e \`a 
$(V ,B)$ : 
$$Sp(V) = \{x \in GL(V)\, / \, B(xv, xw)= B(v, w),\, \mbox{ pour tout } v , w
\in V\}.$$

\noindent Un lagrangien de $V$ est un sous-espace vectoriel totalement isotrope par
rapport \`a $B$ et de dimension maximale. Pour qu'un sous-espace
vectoriel  de $V$ soit un lagrangien il faut et il suffit qu'il
soit \'egal \`a son orthogonal par $B$. Le groupe symplectique
$Sp(V)$ op\`ere transitivement sur les lagrangiens de $V$.
 
\end{subsection}
\begin{subsection}{} Soit $G$ un groupe localement compact. Par une mesure de Haar sur $G$, on entend une mesure de Haar invariante \`a gauche.  
 On se donne    une mesure de Haar $d_Gx$ sur $G$.
 On note $\Delta_G$ la fonction module de $G$~:  $$\int_G\var(xy^{-1})d_Gx =
\Delta_G(y)\int_G\var(x)d_Gx, \, \,  \var \in
C_c^\infty(G), \, y \in G.$$
 Soit $H$  un sous-groupe  ferm\'e   de $G$.
 On d\'esigne par
$\Delta_{H,G}$ le morphisme de groupes de $H$ dans $\reel_+^*$
d\'efini par : $$\Delta_{H,G}(y) =
\frac{\Delta_{H}(y)}{\Delta_{G}(y)} \mbox{, pour tout } y \in H,$$
et par $C_c(G ; H)$ le $\cit$-espace des fonctions $\var$ sur
$G$, continues \`a support compact modulo $H$, v\'erifiant  $$\var(xy) = \Delta_{H,G}(y)\var(x) \mbox{, pour tout } x \in
G, y\in H.$$ Le groupe $G$ op\`ere par translations \`a gauche
dans $C_c(G ; H)$ et il existe, \`a une constante positive pr\`es,
une unique  forme lin\'eaire positive, $G$-invariante, sur cet
espace. Une telle forme lin\'eaire est appel\'ee dans la suite une
mesure invariante sur $G/H$. En particulier,
si $d_Hx$ est  une mesure de
Haar  sur $H$,  il existe une unique mesure invariante
$d_{G/H}\dot x$ sur  $G/H$ telle que $$\int_G \var(x)d_Gx =
\int_{G/H}\{\int_H \var(xy)\Delta_{H,G}(y)^{-1} d_Hy\} d_{G/H}\dot x, \, \, 
\mbox{ pour tout } \var \in C_c^\infty(G). $$ La mesure invariante
 $d_{G/H}\dot x$ s'appelle le quotient des
mesures de Haar   $d_Gx$ et $d_Hx$.

\end{subsection}

\begin{subsection}{}\label{am} Soit $G$ un groupe localement compact.
 Par repr\'esentation unitaire de $G$, on  entend une
repr\'esentation unitaire continue de $G$ dans un espace de
Hilbert s\'eparable. Soit $\pi$ une repr\'esentation unitaire de
$G$  dans  $\mathcal H$.
Si $\var \in C_c^\infty(G)$ et $d_Gx$ est  une mesure de Haar sur
$G$, on d\'efinit l'op\'erateur $\pi(\var d_Gx)=
\int_G\var(x)\pi(x)d_Gx$ agissant dans $\mathcal H$, en posant,
pour tout $v, w \in \mathcal H$, $$<\pi(\var d_Gx)v, w>=
\int_G\var(x)<\pi(x)v,w>d_Gx,$$ o\`u $<\,,\,
>$ d\'esigne le produit scalaire dans $\mathcal H$. 

\noindent On suppose  de plus que la topologie sur $G$ est totalement discontinue. Pour chaque
sous-groupe $K$ de $G$, on d\'esigne par 
$${\mathcal H}^K = \{v \in \mathcal H,\, \,  \pi(x)v = v, \mbox{ pour tout } x \in K \} . $$ 
 La repr\'esentation $\pi$
est dite admissible si, pour tout sous-groupe compact ouvert $K$
de $G$, on a  $\dim {\mathcal H}^K<+\infty$. On voit que $\pi$
est admissible si et seulement si les op\'erateurs $\pi(\var d_Gx)$,  $\var \in C_c^\infty(G)$, sont de rang
finis. Si tel est le cas, on
d\'efinit  le caract\`ere de $\pi$ comme \'etant la fonction
g\'en\'eralis\'ee  sur $G$~: $$ \Theta_\pi(\var d_Gx) = \mbox{ Tr}
(\pi(\var d_Gx)), \quad  \var \in C_c^\infty(G).$$
\end{subsection}
\begin{subsection}{}
Soit ${\G}$  un groupe alg\'ebrique d\'efini sur $k$ dont
l'ensemble des points rationnels est  not\'e ${\G}_k$. On
d\'esigne par ${}^u{\G}$ son radical unipotent, c'est le plus
grand sous-groupe normal ferm\'e unipotent de  ${\G}$. Il est d\'efini sur $k$.  Un
facteur r\'eductif d\'efini sur $k$ de ${\G}$  est un sous-groupe
r\'eductif $\bf R$ (i.e. ${}^u{\bf R}= \{1\}$) d\'efini sur $k$
tel que ${\G}$ soit  produit semi-direct de $\bf R$ et de
${}^u{\G}$. Le sous-groupe $\bf R$  est appel\'e aussi un
$k$-facteur r\'eductif de ${\G}$. D'apr\`es (\cite{[Bo-Se]}, Proposition
5.1) ${\G}$ poss\`ede un $k$-facteur
 r\'eductif $\bf R$ et tout
sous-groupe r\'eductif d\'efini sur $k$ de ${\G}$ est conjugu\'e
\`a un sous-groupe de $\bf R$ par un \'el\'ement de~${}^u{\G}_k$.
\end{subsection}
\begin{subsection}{}\label{DJ}
Un groupe presque alg\'ebrique est un triplet $(G, F, {\G})$, o\`u
${\G}$ est un groupe  al\-g\'e\-bri\-que d\'efini sur $k$, $G$ un
groupe localement compact et  $F$  un sous-groupe  fini du centre de
$G$ tels que 
 $G/F$  soit un
sous-groupe  d'indice fini de ${\G}_k$, dense pour la topologie de
Zariski,  et la
projection canonique $p_G$ de $G$ dans ${\G}_k$ soit continue. 

\noindent Comme  $G/F$ est d'indice fini dans ${\G}_k$, il est  ouvert dans ${\G}_k$ et   $p_G$
 est un hom\'eomorphisme local. Si bien que   l'on peut munir
$G$ d'une structure de groupe de Lie sur $k$ de telle sorte
que $p_G$ soit un diff\'eomorphisme local.  L'alg\`ebre de Lie de $G$ est canoniquement isomorphe \`a  l'alg\`ebre de Lie $\g$ de  ${\G}_k$. Nous dirons que $\g$ est l'alg\`ebre de Lie du groupe presque alg\'ebrique  $(G, F, {\G})$.
Etant donn\'e qu'un sous-groupe unipotent ${\bf U}_k$  de ${\G}_k$  n'a pas de
sous-groupe propre d'indice fini,  il est contenu
dans  $G/F$. D'apr\`es  (\cite{[Du]},  Lemme II.11) $G$ contient un
unique sous-groupe ferm\'e $U$ isomorphe \`a ${\bf U}_k$ par $p_G$.
 Les \'el\'ements de $U$  sont appel\'es \'el\'ements unipotents de $G$.
Il est clair que  $p_G$ r\'ealise une bijection
de l'ensemble des \'el\'ements unipotents de $G$ sur ceux de
${\G}_k$.
 Lorsque  ${\bf U}_k$  est le radical unipotent ${}^u{\G}_k$ de ${\G}_k$, le  sous-groupe correspondant   de $G$, not\'e ${}^uG$,  est  appel\'e le radical unipotent de $G$.
Si $\bf R$ est un $k$-facteur r\'eductif de ${\G}$  alors $G$ est
produit semi-direct de $ p_G^{-1}({\bf R}_k)$ et de ${}^uG$. Le
sous-groupe $ p_G^{-1}({\bf R}_k)$  est appel\'e facteur r\'eductif de $G$. Si $ p_G^{-1}({\bf R}_k)/F$ est ab\'elien, $G$ est dit 
un groupe r\'esoluble presque connexe .

\noindent Un \'el\'ement    de $G$ est dit semi-simple si son image
par $p_G$ est semi-simple  de ${\G}_k$.  On obtient une d\'ecomposition de
Jordan dans $G$~ : chaque \'el\'ement $x$ de $G$ s'\'ecrit de
mani\`ere unique sous la forme
 \begin{eqnarray}
\label{jordan} x& =& x_s.x_u= x_u.x_s ,
\end{eqnarray}
  o\`u $x_s$ est semi-simple appel\'e
  partie semi-simple de $x$ et $x_u$ est unipotent appel\'e partie unipotente de $x$.
  
\noindent Sauf mention explicite du contraire le groupe $G$ op\`ere sur
lui-m\^eme par automorphismes int\'erieurs, sur $\g$ par l'action adjointe, et sur $\g^*$  par l'action co-adjointe.
Si $\Omega$ est une orbite co-adjointe  dans $\g^*$, nous d\'esignons par $d\mu_\Omega$ la mesure de Liouville sur $\Omega$ (voir \cite{[M3]}, Paragraphe 11).
\end{subsection}
\end{section}

\begin{section}{Voisinages semi-simples, Application exponentielle.}\label{amn}
 \begin{subsection}{}\label{alaa} Soit $(G, F, {\G})$  un groupe presque alg\'ebrique  d'alg\`ebre de Lie $\g$.
On rappelle  qu'un ouvert $\mathcal W$ de $G$ est dit semi-simple
(ou $G$-semi-simple s'il est n\'ecessaire de pr\'eciser le groupe
$G$) s'il est invariant par conjugaison, et si, pour
tout $x \in G$, on a $x \in \mathcal W$ si et seulement si $x_s
\in \mathcal W$, o\`u $x_s$ est la partie semi-simple de $x$.

De m\^eme un ouvert $\mathcal V$ de $\g$ est dit semi-simple (ou
$G$-semi-simple s'il  est n\'ecessaire de pr\'eciser $G$) s'il est
invariant par l'action adjointe de $G$, et si, pour tout $X \in
\g$, on a $X \in \mathcal V$ si et seulement si $X_s \in \mathcal
V$, o\`u $X_s$ est la partie semi-simple de $X$.

Les ouverts semi-simples de $G$ (resp. $\g$) sont les ouverts
d'une topologie, appel\'ee la toplologie semi-simple, sur $G$
(resp. $\g$).

\begin{subsubsection}{} \label{vois}
On choisit d\'esormais une r\'ealisation de  ${\G}$ comme un
sous-groupe alg\'ebrique  de $GL_m(\bar k)$, d\'efini sur $k$,
 et donc $\g$ comme une sous-alg\`ebre de Lie de
$\goth{gl}_m(k)$.  On pose, pour $X \in \goth{gl}_m(k)$,   $$\rho(X)
= \sup_{\lambda \in spec(X)}|\lambda|_{\bar k},$$
o\`u $ spec(X)$ d\'esigne
l'ensem\-ble des   valeurs propres de l'action de  $X$ dans ${\bar k}^m$.
Pour  $\e >0$, on pose
 $$\g_\e= \{X \in \g,  \rho(X) \leq \e\},
{\G}^{ \e}_{k} = \{x \in {\G}_{k},  \rho(x-1) \leq \e\}.$$ Le
r\'esultat suivant est d\'emontr\'e  dans \cite{[M3]}.
\begin{tho}\label{yvois}
 i) Les ouverts
$ {\G}^{ \e}_{k}$, $\e>0$, forment une base  de voisinages
semi-simples de $1$  dans $ {\G}_{k}$.

 ii) Les ouverts $\g_\e$, $\e>0$ sont invariants par translations par ${}^u\g$ et  forment une base  de voisinages
semi-simples de  $0$ dans $\g$.
\end{tho}

\medskip

\end{subsubsection}
\begin{subsubsection}{}\label{Adem Ch}
 On fixe $0<a<1$ tel que $\dis \lim_{n \rightarrow \infty}\frac{a^n}{|n! |_p} = 0$ et que
  $a>\frac{a^n}{|n! |_p}$, pour tout $n \geq 2$ (voir  \cite{[Se4]}, Lemme V.4).  Le lemme suivant
  est  d\^u \`a Harish-Chandra (non publi\'e).
\begin{lem} \label{yAdem Ch} L'application exponentielle
$$\exp :\goth{gl}_m(k)_a \longrightarrow GL_m( k)^a, X \longmapsto
\sum_{n \geq 0} \frac{X^n}{n! } , $$ est  un diff\'eomorphisme,
dont l'inverse est l'application logarithme~: $$\log(x) = \sum_{n
\geq 1 }(-1)^{n-1}\frac{(x-1)^n}{n }, \quad x \in GL_m( k)^a .$$
Si $X \in \goth{gl}_m(k)_a$ et $ x \in GL_m( k)^a$, on a :
$$\exp(yXy^{-1}) =y \exp(X)y^{-1} \, , \, \, \,  \log(yxy^{-1})= y
\log(x)y^{-1}\, , \, \, \mbox{ pour tout } y \in GL_m( k). $$
\end{lem}

\medskip

Pour $0< \e\leq a$,
on pose ${\bf G}_{k, \e}= \exp(\g_\e)$.
 On a la proposition  suivante~: 
\begin{pro} \label{yyAdem Ch} Pour tout $0< \e\leq a$, ${\bf G}_{k, \e}$ est un
voisinage ${\bf G}_{k}$-semi-simple de $1$ dans  ${\bf G}_{k}$ et
l'application exponentielle $\exp $  r\'ealise  un diff\'eomorphisme de $\g_\e$ sur ${\bf
G}_{k, \e}$.
\end{pro}
\begin{dem} Soit $0< \e\leq a$. D'apr\`es (\cite{[Ch]}, Paragraphe  II.12), si $X= X_s+ X_u \in
\g_\e$  alors $x= \exp(X) \in {\bf G}_{k}$ et on a $x_s =
\exp(X_s)$, $x_u = \exp(X_u)$. De plus, 
l'application exponentielle induit une bijection de l'ensemble des
\'el\'ements   nilpotents de $\g$ sur l'ensemble des \'el\'ements
unipotents de $ {\bf G}_{k}$. Ceci prouve que ${\bf G}_{k, \e}$
est un voisinage ${\bf G}_{k}$-semi-simple de $1$ dans  ${\bf
G}_{k}$.  D'apr\`es le lemme \ref{yAdem Ch},  $\exp : \g_{\e} \longrightarrow {\bf G}_{k,
\e}$  est un diff\'eomorphisme de $\g_\e$ sur ${\bf G}_{k, \e}$.
\end{dem}
\end{subsubsection}
\begin{subsubsection}{} \label{red} Maintenant, nous allons remonter
l'application exponentielle \`a $G$.
   On note    $$ (G/F)_{inv}= \bigcap_{x \in {\G}_{k}} x  (G/F) x^{-1}.$$
 Alors $(G/F)_{inv}$ est un sous-groupe normal, ferm\'e, et d'indice fini de ${\G}_{k}$~;
  il contient tous les \'el\'ements unipotents de ${\G}_{k}$. Ainsi $(G/F)_{inv}$ est
   un ouvert ${\G}_{k}$-semi-simple dans ${\G}_{k}$. On pose
 $$a_G = \sup\{0< \e \leq a, \, \, \mbox{ tel que } {\G}_{k, \e} \subset (G/F)_{inv} \}, \, \,  n_F =  \mbox{  cardinal de } F ,$$
 et pour $0<\e <a_G$,
 $$G_{\e} = \{x^{2n_F} \, , \, \, \, x \in p_G^{-1}({\G}_{k, \e})\}.$$
On a le r\'esultat suivant.
\begin{pro} \label{yred} Pour tout $0< \e <  a_{G} $, on a~:

  $i)$ $G_{\e}$ est un  ouvert $G$-semi-simple de $G$.
   
   $ii)$ la restriction de $p_G$ \`a  $G_{\e}$ est un diff\'eomorphisme  de $G_{\e}$
     sur ${\G}_{k, \e |2n_F|_p }$.
 \end{pro}
\begin{dem}  On note $\alpha_1: p_G^{-1}( {\G}_{k, \e}) {\longrightarrow} p_G^{-1}( {\G}_{k, \e |2n_F|_p}),\, x 
\longmapsto x^{2n_F},$
$$
 \alpha_2 :   {\G}_{k, \e}   {\longrightarrow}  {\G}_{k, \e|2n_F|_p}, \, x \longmapsto x^{2n_F}, \mbox{ et } \alpha_3 : \g_\e \longrightarrow  \g_{\e|2n_F|_p},\,  X\longmapsto 2n_F X. $$
   On consid\`ere le diagramme commutatif suivant~:
\begin{displaymath}
  \xymatrix{ \g_\e\ar[r]^{\exp}\ar[d]_{\alpha _3}&
   {\G}_{k, \e}\ar[d]_{\alpha _2} &  p_G^{-1}( {\G}_{k, \e}) \ar[d]_{\alpha _1}\ar[l]_{p_G} \\
   \g_{\e|2n_F|_p} \ar[r]_{{\exp}}&  {\G}_{k, \e|2n_F|_p}
  & p_G^{-1}( {\G}_{k, \e|2n_F|_p} )\ar[l]^{p_G} }
\end{displaymath}
Puisque $\exp$ et $\alpha_3$ sont des diff\'eomorphismes donc   $\alpha_2$ l'est aussi. Comme  $p_G$ est un diff\'eomorphisme local et  $\alpha_1$ est continue, on d\'eduit que   $\alpha_1$ est un diff\'eomorphisme local.
Il s'en suit que   $G_\e$ est un  ouvert de $G$,  invariant par conjugaison par les \'el\'ements
de $G$.
 Soit $x= x_sx_u \in G$, o\`u
$x_s$ est la partie semi-simple et $x_u$ la partie unipotente de
$x$.  Si $x_s \in G_\e$  alors  il existe  $y \in
p_G^{-1}({\G}_{k, \e})$ tel que $x_s = y^{2n_F}$. \'Ecrivons $p_G(y) =
\exp(Y)$, avec $Y \in \g_\e$ semi-simple, et $p_G(x_u) = \exp(2n_FZ)
= z^{2n_F}$, avec $Z \in \g$ unipotent. On a~:  $[Y, Z] = 0$. En effet,
puisque $p_G(x_s)$ commute avec $\exp(2n_FZ)$ donc
$p_G(x_s)$ commute avec $\exp(tZ), t \in k$.  Si bien que l'on a 
$\exp(t \, adZ)Y = Y, \mbox{ pour tout }   t \in k.$ Par suite $[Y, Z] = 0$. On en d\'eduit  que
 $p_G(y)z = zp_G(y) \in {\G}_{k, \e} $.  Soit $\tilde z$
l'\'el\'ement unipotent de $G$ correspondant \`a $z$. Alors, il
existe $f \in F$ tel que $y \tilde z = f\tilde z y$. Mais $f^{n_F(2n_F-1)}=1$, on a  $x=
y^{2n_F}\tilde z^{2n_F} = (y\tilde z)^{2n_F} \in G_\e$.  R\'eciproquement, il est clair que
si $x \in G_\e$ alors $x_s \in G_\e$.
Concernant le deuxi\`eme point, on a~:  $p_G(G_\e) = {\G}_{k, \e
|2n_F|_p }$. L'injectivit\'e de la restriction de $p_G$ \`a $G_\e$
se d\'emontre sans difficult\'e. Le r\'esultat se d\'eduit du fait
que $p_G$ est un diff\'eomorphisme local et de la proposition \ref{yyAdem Ch}.
\end{dem}
 Pour  $0< \e <  a_{G} $, on vient de d\'emontrer que  $p_G : G_{\e}
    \longrightarrow{\G}_{k, \e |2n_F|_p }$ est un diff\'eomor\-phis\-me ;
 on note $q_G$ l'application inverse.
On d\'efinit ainsi une application, appel\'ee  application exponentielle, $\exp_G : \g_{\e |2n_F|_p}
 \longrightarrow G_{\e},  \quad
 \exp_G = q_G \circ \exp .$
\end{subsubsection}
\begin{subsubsection}{}\label{khem1}
 Soit  $ x$  un \'el\'ement  de $G$.   On note  ${\G}_k(p_G(x))$
le centralisateur de $p_G(x)$ dans   ${\G}_k$. C'est un sous-groupe alg\'ebrique
de ${\G}_k$ d'alg\`ebre de Lie $\g(x)$. Soit $0<\e \leq a$. On a~: ${\G}_k(p_G(x))_{ \e} = {\G}_{k ,\e} \cap {\G}_k(p_G(x))$ et  l'application exponentielle $\exp$ induit un diff\'eomorphisme de $\g(x)_\e$ sur ${\G}_k(p_G(x))_{ \e}$.  D'autre part, 
l'application~:  $$p_G^{-1}({\G}_k(p_G(x))) \longrightarrow F, \, y \longmapsto xyx^{-1}y^{-1}, $$
est un morphisme de groupes. Son noyau est $G(x)$, le centralisateur de $x$ dans $G$, qui est un sous-groupe ferm\'e et d'indice fini de $p_G^{-1}({\G}_k(p_G(x)))$.  De plus, si  $y \in G(x)$ et $y=
y_sy_u$ sa d\'ecomposition de Jordan,  on a~:  $y_s, y_u \in G(x)$.
Pour
$0<\e < a_{G}$, 
 on pose
  $$G(x)_\e = \{y^{2n_F} \, , \, y \in p_G^{-1}(({\G}_k(p_G(x)))_{ \e})\}.$$
\begin{lem} \label{ykhem1} Soit  $ x \in G$ et $0<\e < a_{G}$. On a
$G(x)_\e = G_\e \cap G(x)$ et la restriction de $\exp_G$ \`a  $ \g(x)_{\e |2n_F|_p}$ est un diff\'eomorphisme
    de  $ \g(x)_{\e |2n_F|_p}$ sur  $G(x)_\e$.
\end{lem}
\begin{dem} Si $y \in p_G^{-1}({\G}_k(p_G(x))_{ \e})$, il existe $f \in F$ tel que $y.x = f x.y$~; par suite $y^{2n_F } \in G(x)$. Ainsi  $G(x)_\e  \subset G_\e \cap G(x)$. Montrons maintenant l'inclusion inverse. Soit $z \in  G_\e \cap G(x)$. Par d\'efinition, il existe $y \in   p_G^{-1}({\G}_{ k, \e})$ tel que $z = y^{2n_F}$.  Ecrivons $p_G(y) = \exp(Y)$
avec $Y \in \g_\e$.
 On voit que $Ad(p_G(x)).Y = Y$, et par suite  $\exp(Y) \in {\G}_k(p_G(x))_{ \e}$. Si bien que $z \in G(x)_\e$.

La deuxi\`eme assertion du lemme r\'esulte du fait que $p_G$ induit un diff\'eomorphisme de $G(x)_\e$ sur 
${\G}_k(p_G(x))_{ \e|2n_F|_p}$,  d'inverse  la restriction de $q_G$ \`a ${\G}_k(p_G(x))_{ \e|2n_F|_p}$,  et du fait que $\exp$ r\'ealise un diff\'eomorphisme de $\g(x)_{ \e|2n_F|_p}$ sur ${\G}_k(p_G(x))_{ \e|2n_F|_p}$.
\end{dem}
\end{subsubsection}
\begin{subsubsection}{}\label{khem2} Soit $\u$ un id\'eal alg\'ebrique
$G$-invariant  de $\g$ et $u$ une forme lin\'eaire sur $\u$. On
note $H = G(u)$ et $\h = \g(u)$. Si 
$0<\e < a_{G}$, 
 on pose
  $$H_\e = \{y^{2n_F} \, , \, y \in p_G^{-1}(({\G}_k(u))_{ \e})\}.$$
\begin{lem} \label{ykhem2}
Soit $0<\e < a_{G}$.  On a $H_\e = G_\e \cap H$
    et la restriction de $\exp_G$ \`a  $ \h_{\e |2n_F|_p}$ est un diff\'eomorphisme
    de  $ \h_{\e |2n_F|_p}$ sur  $H_\e$.
\end{lem}

\begin{dem} 
 Soit $x \in  G_\e \cap
 H$.  Il existe  $y \in p_G^{-1}({\bf G}_{k, \e})$ tel que $ x = y^{2n_F}$.
 Ecrivons $p_G(y)= \exp(Y)$, \, $Y \in \g_\e$. On a~: $p_G(y^{2n_F})
 = \exp(2n_FY) \in  {\bf G}_k(u)$. Puisque ${\bf G}_k(u)$ est un sous-groupe alg\'ebrique de ${\bf G}_k$, on  a $\exp(2n_FY_s), \exp(2n_FY_u) \in {\bf G}_k(u)$.  Il en r\'esulte  que $Y_s, Y_u \in \h$ et que $Y\in \h$.
 Ainsi $\exp(Y) \in ({\G}_k(u))_{ \e}$. Si bien que $x \in H_\e$. D'o\`u $H_\e = G_\e \cap H$. 
\end{dem}
\end{subsubsection}
\begin{subsubsection}{}\label{resoluble} Dans ce paragraphe, on
suppose que $\g$ est r\'esoluble.  Soit $\t$ un facteur r\'eductif
de $\g$ et  ${\bf T}_k$ un facteur r\'eductif de ${\bf G}_{k}$
d'alg\`ebre de Lie $\t$.  Pour  $0<\e \leq a$, on note ${\bf T}_{k, \e} = \exp(\t_\e)$, \,  $T = p_G^{-1}({\bf T}_k)$, et $T_\e = \{y^{2n_F} \, , \, y \in p_G^{-1}({\bf T}_{k, \e})\}$.

\begin{pro} \label{yresoluble} Si $a$ est suffisamment petit alors, pour tout $0<\e \leq a$, on a  : 
$${\bf G}_{k, \e} = {\bf T}_{k, \e}{}^u{\bf G}_{k}, $$ 
qui est  un
sous-groupe normale de  ${\bf G}_{k}$. La loi dans ${\bf
G}_{k, \e}$ est donn\'ee par la formule de Campbell Hausdorff. Si
de plus $0<\e<a_G$,  on
 a $G_\e = T_\e {}^uG$, qui est un sous-groupe normal de $G$.
\end{pro}
\begin{dem} Il est clair que, pour tout $\e>0$, $\t_\e$ est un r\'eseau de $\t$ et que la famille
$(\t_\e)_{\e>0}$ est une base de voisinages de $0$ dans $\t$. Il
r\'esulte  de (\cite{[Ho3]}, Lemme 1.1) et du th\'eor\`eme \ref{yvois} que si $a$ est
suffisamment petit, pour tout $0<\e \leq a$, ${\bf G}_{k, \e}$ est un
sous-groupe normal de  ${\bf G}_{k}$ et  la loi dans ${\bf
G}_{k, \e}$ est donn\'ee par la formule de Campbell Hausdorff.
Comme ${\bf T}_{k, \e}$ et ${}^u{\bf G}_{k}$ sont des sous-groupes
de ${\bf G}_{k, \e}$, on a ${\bf T}_{k, \e}{}^u{\bf G}_{k} \subset
{\bf G}_{k, \e}$. D\'emontrons  l'inclusion inverse.
D'apr\`es le th\'eor\`eme \ref{yvois}, on a $\g_\e = \t_\e +{}^u\g$.
Soit  $X \in \t_\e$ et $Y \in {}^u\g$. En utilisant la formule de
Campbell Hausdorff, on voit que $\exp(-X)\exp(X + Y) \in {}^u{\bf
G}_{k}$. Si bien que $\exp(X + Y) \in {\bf T}_{k, \e}{}^u{\bf
G}_{k}$.

\noindent Supposons maintenant $0<\e<a_G$. Soit $x \in G_\e$ et $y \in
p_G^{-1}({\bf G}_{k, \e})$ tel que $x = y^{2n_F}$. \'Ecrivons $y=
y_sy_u$ sa d\'ecomposition de Jordan. Quitte \`a conjuguer $y$ par
un \'el\'ement de ${}^uG$, on peut supposer que $y_s \in
p_G^{-1}({\bf T}_{k, \e})$. Ainsi $y_s^{2n_F} \in T_\e$ et par suite 
 $x =y_s^{2n_F}y_u^{2n_F} \in T_\e{}^uG$. Si bien que $G_\e
\subset T_\e{}^uG$. D'autre part la restriction de $p_G$ \`a
$T_\e{}^uG$ r\'ealise une bijection  de $T_\e{}^uG$ sur
${\bf T}_{k, \e|2n_F|_p}{}^u{\bf G}_{k} ={\bf G}_{k, \e|2n_F|_p}$.
D'apr\`es la proposition \ref{yred}, on a $G_\e =T_\e{}^uG$. La derni\`ere assertion de la proposition d\'ecoule du fait que $T_\e$ est un sous-groupe normal  de $T$ et ${}^uG$ est un sous-groupe normal de $G$.
\end{dem}
\end{subsubsection}
\end{subsection}
\begin{subsection}{} \label{YY}
Dans ce paragraphe, on se donne    un \'el\'ement semi-simple  $s$ de $G$. On 
d\'esigne par $\lambda_1, \ldots, \lambda_l$ (resp. $ \nu_1, \ldots,  \nu_{l'}$) les valeurs
propres distinctes de $p_G(s)$ dans ${\bar k}^{\, m} $ (resp. de $Ads$ dans $\g_{\bar k} =\bar k\otimes_k \g $).
On d\'efinit les nombres $a'_\g(s)$ et $a_G(s)$ que l'on note simplement $a'(s)$ et $a(s)$  lorsque aucune confusion n'est possible, 
\begin{eqnarray} \label{conditionss}
&& a(s) = \inf_{ i \neq j}\{|\lambda_i \lambda_j^{-1} -1|_{\bar k},  a_G  \}, \, \,
 a'(s)=  \inf_{ \nu_i \neq 1 }\{| \nu_i -1|_{\bar k},
   | \nu_i^{-1} -1|_{\bar k},\,  a_G  \}.
\end{eqnarray}
Il est clair que $a(s) \leq a'(s)$, $a(1) = a_G$, et que si $Ads =
\mbox{Id}_\g$, $a'(s)= a_G$.
 On a le r\'esultat suivant~:
\begin{lem} \label{yYY}
 Soit $0<\e'(s)< a'(s)$ et $0<\e(s)< a(s)$.

$ i)$  Soit $y \in G(s)_{\e'(s)}$. Si un vecteur d'un sous-quotient
$s$-invariant de $\g$ ou de $\g^*$ est fix\'e par $sy$, alors il
est fix\'e par $s$.

$ ii)$ Soient  $g \in G, \, y, y' \in G(s)_{\e(s)}$ tels
que l'on ait $sy'= gsyg^{-1}$. Alors $g\in G(s)$.

\end{lem}
\begin{dem}
On commence par d\'emontrer la premi\`ere assertion du lemme. Soit  $W$
un sous-quotient
$Ads$-invariant de $\g$  et $v \in W$ un vecteur fix\'e par $Ad(sy)$. On peut supposer que $v  \neq 0$.  Ecrivons $v = \sum_ \nu v_ \nu$, o\`u $ \nu$ parcourt l'ensemble des valeurs propres de $Ads$ dans $W_{\bar k}$ et $v_ \nu$ appartient au sous-espace propre correspondant. Soit $ \nu$ tel que $v_ \nu \neq 0$. Comme $Ads$ et $Ady$ commutent, on a $Ad(sy). v_ \nu = v_ \nu$ et par suite $v_ \nu$ est vecteur propre de $Ady$ pour la valeur propre $ \nu^{-1}$.  Ainsi,  $ \nu$ est le quotient de deux valeurs propres $\alpha_1$ et $\alpha_2$ de $p_G(y)$ dans ${\bar  k}^m$. On a~:
 $$ | \nu^{-1} -1|_{\bar k} = |\frac{\alpha_2}{\alpha_1} -1|_{\bar k} \leq 
\max \{|{\alpha_1} -1|_{\bar k}, |{\alpha_2} -1|_{\bar k}\} \leq \e'(s) .  $$
Compte tenu du choix de $\e'(s)$, on a n\'ecessairement $  \nu = 1$. Si bien que l'on a $Ads.v = v$.
On d\'emontre  de  la m\^eme  fa\c{c}on que si un sous-quotient $Ad^*s$-invariant de $\g^*$  est fix\'e par $Ad^*(sy)$ alors il
est fix\'e par $Ad^*s$.

D\'emontrons maintenant  la deuxi\`eme assertion~:  si $\lambda$ est une   valeur propre  de $p_G(s)$,
  on note $V_\lambda$ le sous-espace propre correspondant. Chacun  de ces  sous-espaces propres
  est stable par $p_G(y_s)$ et $p_G(y'_s)$.  Soit $ 1 \leq i \leq l$ et  $ v \in V_{\lambda_i}$
   un vecteur propre de $p_G(sy'_s)$
  pour une valeur propre $ \lambda  $.
   Comme $sy'_s$ est conjugu\'e \`a  $sy_s$,   $ \lambda$  est aussi une valeur
    propre de $p_G(sy_s)$.
   Il existe donc $1 \leq j \leq l$  et $\alpha$
     (resp. $\alpha'$)  une valeur propre de $p_G(y_s)$ (resp. $p_G(y'_s)$) tels
     que $\lambda =\lambda_j \alpha = \lambda_i\alpha'$.
     On a :  $$|\frac{\lambda_{i}}{\lambda_j}-1 |_{\bar k} =
      |\frac{\alpha}{\alpha'}-1|_{\bar k} \leq \max\{|{\alpha}-1|_{\bar k}, \,
       |{\alpha'}-1|_{\bar k}\} \leq \e(s)|2n_F|_p \leq \e(s) .$$
     Vu le choix de $\e(s)$, $ j= i$ et donc $\alpha = \alpha'$. On en d\'eduit que $p_G(g).v \in V_{\lambda_i}$.
     Si bien que, $p_G(g)$ commute \`a $p_G(s)$.
      Ainsi, il existe  $f \in F$ tel que $gsg^{-1} = fs$.
         Mais, ceci implique que $y'= fgyg^{-1}$ et par suite $f=
         1$. D'o\`u $g\in G(s)$. \end{dem}
\end{subsection}
\begin{subsection}{}\label{reductive}
Soit $(M, F, {\bf M})$ un groupe presque alg\'ebrique r\'eductif d'alg\`ebre de Lie $\goth m $.

\begin{pro} \label{yreductive}
On a :

$ i)$ Si $X \in \goth m$ est nilpotent, $0$ est dans
l'adh\'erence de  l'orbite de $X$  sous l'action adjointe de $M$ (pour  la
topologie $p$-adique).

 $ii)$ Si $x \in M$  alors $x_s$, la partie semi-simple de
$x$, appartient \`a l'adh\'erence de  l'orbite de $x$ sous
l'action conjugaison de $M$ (pour la topologie $p$-adique).

$ iii)$ Si $X \in \goth m$ alors  $X_{s}$,  la partie
semi-simple de $X$, appartient \`a  l'adh\'erence de l'orbite de $X$
sous l'action adjointe de $M$ (pour la topologie $p$-adique).
\end{pro}
\begin{dem}
D'apr\`es (\cite{[Bo-Se]}, Corollaire 6.6), les \'el\'ements
nilpotents de $ \goth m$ forment un nombre fini
 de classes pour
l'action adjointe  de ${\bf M}_k$. Etant
donn\'e que $M/F$ est d'indice fini de ${\bf M}_k$, on a le m\^eme
r\'esultat pour l'action adjointe  de $M$.
Soit $X \in \goth m$  nilpotent et $t \in k^\times$ tel que $|t|_p
<1$. On consid\`ere la suite $(t^nX)_{n \in \N}$. D'apr\`es ce qui
pr\'ec\`ede, il existe $n_1, n_2 \in \N$ avec $n_1 < n_2$ tels que
$t^{n_1}X$ et $t^{n_2}X$  soient dans une m\^eme $M$-orbite~; soit
$y \in M$ tel que $y.X = t^{n_2-n_1}X$. Alors la suite $(y^m.X)_{m \in
\N}$ converge vers $0$.

\noindent Montrons maintenant l'assertion ii).  Soit $x= x_sx_u \in M$,  o\`u $x_s$ est la partie semi-simple et
$x_u$ la partie unipotente de $x$. \'Ecrivons $x_u = \exp(X)$, $X \in \goth m(x_s)$
nilpotent. Comme  $\goth m(x_s)$  est r\'eductive, en utilisant le raisonnement  pr\'ec\'edent, on en
d\'eduit qu'il existe $y \in (M/F)(p_M(x_s))$ tel
que la suite $(y^m.X)_{m \in \N}$ converge vers $0$,  $p_M$ \'etant la
projection canonique de $M$ dans ${\bf M}_k$.  
  Ainsi la suite $(y^mx_uy^{-m})_{m \in \N}$ converge vers $1$. Soit
$\tilde y \in M$ tel que $p_M(\tilde y) = y$. Alors $\tilde y^{n_F}$ commute avec
$x_s$ et la suite
  $({\tilde y}^{n_Fm}x {\tilde y}^{-n_Fm})_{m \in \N}$ converge vers $x_s$. La preuve de iii) est de m\^eme
style que celle de ii).
\end{dem}
\begin{lem} \label{yyreductive}
 Soit
$\Omega$ un voisinage ouvert $M$-invariant de $0$ (resp. $1$) dans $\goth m$ (resp. $M$).
 Alors, il existe
$\e>0$ tel que $\goth m_\e \subset \Omega$ (resp. $M_\e \subset \Omega$).
\end{lem}
\begin{dem}
 En effet, si ce n'\'etait
pas le cas, on pourrait construire une suite $(X_m)_{m \in {\N^\times}}$
d'\'el\'ements de $\goth m$ telle que, pour tout $m \in \N^\times$,
\begin{enumerate}
\item $X_m \in \goth m_{\frac{1}{m}}$ ;
\item $X_m \notin \Omega $.
\end{enumerate}
On note $X_{m, s}$ la partie semi-simple de $X_m$, $m \in {\N^\times}$.  On a $X_{m, s} \in \goth m_{\frac{1}{m}},  m \in {\N^\times}$.
Comme $X_{m, s}$ est dans l'adh\'erence de l'orbite de $X_m$
(Proposition \ref{yreductive}),
 on voit que $X_{m, s} \notin {\Omega} $, $m \in {\N^\times}$.
  On peut 
 supposer que, pour tout $m \in {\N^\times}$,  $X_m$ est semi-simple.
 Puisque $\goth m$ ne poss\`ede qu'un nombre fini de classes de
 conjugaison de sous-alg\`ebres de Cartan (voir \cite{[Bo-Se]}, Corollaire 6.5), on peut supposer
 \'egalement qu'il existe une sous-alg\`ebre de Cartan $\h$ telle
 que $X_m \in \h_{\frac{1}{m}}$, pour tout $ m \in {\N^\times}$. Comme
 $(\h_{\frac{1}{m}})_{  m \in {\N^\times}}$ est une base de voisinages de $0$
 dans $\h$ donc $(X_m)_{m
\in {\N^\times}}$ converge vers $0$. Si bien que $X_m \in \Omega$ pour $m$
suffisamment grand.

\noindent  Soit $0<\e'<a_M$. Alors $\exp_M(\m_{\e'|2n_F|_p}) \cap\Omega$ est un voisinage ouvert $M$-invariant de $1$ dans $M$. D'apr\`es le raisonnement pr\'ec\'edent, il existe $0<\e<a_M$ tel que 
$\exp_M(\m_{\e|2n_F|_p}) \subset \exp_M(\m_{\e'|2n_F|_p}) \cap\Omega.$
\end{dem}

\end{subsection}
\begin{subsection}{} Soit $(G, F, {\G})$ un groupe  presque alg\'ebrique  d'alg\`ebre de Lie $\g $
  et  $0< \e <  a_{G} $.  Soit  $\m$ est un r\'eseau de $\g$ contenu dans $\g_{\e |2n_F|_p}$ et tel que $[\m, \m] \subset  \m$. Si  $m$ est un entier suffisamment grand alors $\exp_G(\varpi^m\m)$ est un sous-groupe compact ouvert  de $G$ contenu dans $G_\e$. 
 On se donne une mesure de Haar  $d_Gx$ (resp. $d_\g X$)  sur $G$ (resp. $\g$).
 On dit que $d_Gx$ et $d_\g X$ sont tangentes (ou se correspondent)  si 
 $$ \int_G1_{\exp_G(\varpi^m\m)}(x)d_Gx =\int_\g 1_{\varpi^m\m}(X)d_\g X. $$
Soit $H$ un sous-groupe ferm\'e de $G$ d'alg\`ebre de Lie $\h$, et  $d_Hx$ (resp.  $d_\h X$) une mesure de Haar sur $H$ (resp. $\h$) telles que ces deux mesures soient tangentes.
Alors la mesure invariant  $d_{G/H}\dot x$ ne d\'epend de $d_\g X$
et $d_\h X$ que par la mesure quotient $d_{\g/\h}\dot X$ sur $\g/\h$. On dit
aussi que
 $d_{G/H} \dot x$ et $d_{\g /\h} \dot  X$ sont tangentes.
 \end{subsection}
\end{section}
\begin{section}{M\'ethode de descente.}\label{aick2}
Dans ce paragraphe, nous allons  exposer la m\'ethode de descente due \`a
 Harish-Chandra dans le cas
r\'eductif $p$-adique (voir \cite{[Ha]}), \`a  M. Duflo et M.
Vergne dans le cas des groupes presque alg\'ebriques r\'eels (voir
\cite{[Du-Ve]}), et \'etendu par l'auteur  au cas des groupes
 presque alg\'ebriques  $p$-adiques (voir  \cite{[M3]}).
 \begin{subsection}{}\label{Harish}
On a le r\'esultat suivant qui est   d\^u \`a
Harish-Chandra (voir \cite{[Ha]}).
\begin{pro}\label{yHarish}
Soient $X$ et $Y$ deux vari\'et\'es $k$-analytiques et $f : X
\longrightarrow Y$ une submersion surjective. Soient $\mu_X$ et $\mu_Y$ deux formes volumes  sur $X$ et $Y$ respectivement. Alors, pour tout
$\alpha \in C_c^\infty(X)$, il existe une unique fonction
$f_\alpha \in C_c^\infty(Y)$ telle que, pour tout $\beta \in
C_c^\infty(Y)$, on ait
\begin{eqnarray} \label{descentes}
  \int_X \alpha(x)
(\beta \circ f)(x) d\mu_X(x) &=& \int_Y f_\alpha(y)
\beta(y)d\mu_Y(y). \end{eqnarray}
 De plus, $\mbox{supp} (f_\alpha)\subset
f(\mbox{supp }\alpha)$ et si $\beta$ est une fonction mesurable
sur $Y$ alors $\beta$ est localement int\'egrable sur $Y$ si et
seulement si $\beta \circ f$ l'est sur $X$ et, dans ce cas,
l'\'egalit\'e (\ref{descentes}) reste vraie. L'application $\alpha
\mapsto f_\alpha$ est surjective  de $C_c^\infty(X)$ sur
$C_c^\infty(Y)$.
\end{pro}
\medskip
\end{subsection}
\begin{subsection}{}\label{semis}
Consid\'erons  la situation suivante :
 soit $(G, F, {\bf G})$ un groupe presque alg\'ebrique  d'alg\`ebre de Lie $\g$ et  $s$ un
automorphisme de $G$ dont la diff\'erentielle, not\'ee encore $s$,
est un automorphisme semi-simple de $\g$. On note $G(s)$ le
sous-groupe des points fixes de $s$ dans $G$ et $\g(s)$ son
alg\`ebre de Lie.
On d\'esigne par  $ \nu_1, \ldots,  \nu_l$
les valeurs propres distinctes de $s$ dans $\g_{\bar k} = \bar k \otimes_k  \g
 $. On pose
 $$a'(s) =\inf_{  \nu_i \neq 1}\{| \nu_i  -1|_{\bar k}, 
| \nu_i^{-1} -1|_{\bar k},    a_{G}\}.$$
 Pour $0<\e< a'(s)$, on consid\`ere
l'application $\psi : G \times G(s)_{ \e} \longrightarrow G$, $(y,
z) \mapsto  yzs(y^{-1})$. C'est une submersion en tout point.
Ainsi l'image de $\psi$, not\'ee $\W(s , \e)$, est un ouvert  de
$G$, $*$-invariant, o\`u  $*$ est l'op\'eration de $G$   sur
lui-m\^eme  d\'efinie par~: $$y*z= yzs(y^{-1}), \,  y, z \in G.$$
 On munit $G$ (resp. $G(s)$) d'une mesure de Haar $d_Gx$ (resp.
$d_{G(x)}y$). Si $ \alpha \in C_c^\infty(G \times G(s)_{ \e})$, on
d\'esigne par $\beta_\alpha $ l'\'el\'ement de $C_c^\infty(G(s)_{
\e})$, d\'efini par $$\beta_\alpha(y) = \int_G \alpha(x,y)d_Gx,\,
y\in G(s)_{ \e}. $$
\begin{tho} \label{ysemis}
On suppose que $G$ est unimodulaire. Soit  $0<\e< a'(s)$. 
 Soit $\Theta$ une fonction g\'en\'eralis\'ee $*$-invariante sur
$\W(s , \e)$.
 Il existe une unique fonction g\'en\'eralis\'ee $\theta_s$ sur $G(s)_{ \e}$,
 $G(s)$-invariante  telle que, pour toute $ \alpha \in C_c^\infty(G \times  G(s)_{\e})$, on ait
 $$\Theta(\psi_\alpha d_Gx) = \theta_s(\beta_\alpha d_{G(s)}y).$$
De plus, si $\theta_s = 0$ alors $\Theta = 0$.
\end{tho}

\medskip

{\it Exemple :} Soit $\pi$ une repr\'esentation admissible de $G$
et  $S$ un op\'e\-ra\-teur born\'e et inversible  dans l'espace de $\pi$ tels
que l'on ait $$ S\pi(x) S^{-1}= \pi(s(x)), \mbox{ pour tout } x
\in G.$$ La fonction g\'e\-n\'e\-ra\-li\-s\'ee $\Lambda$  sur $G$
d\'efinie par : $$\Lambda(\var d_Gx)= \mbox{ Tr}(S\circ \pi(\var
d_Gx)), \,  \var \in C_c^\infty(G)$$ est $*$-invariante. On obtient ainsi une fonction
g\'e\-n\'e\-ra\-li\-s\'ee $\lambda_s$ sur $G(s)_{ \e}$,
 $G(s)$-invariante, qui d\'etermine enti\`erement $\Lambda$ sur $\W(s , \e)$.
\end{subsection}

\begin{subsection}{}\label{PTo} Dans ce paragraphe, on se donne  un groupe presque alg\'ebrique
 $(G, F, {\G})$ d'alg\`ebre de Lie $\g$, $s$ un \'el\'ement
 semi-simple de $G$. On reprend les notations de la section \ref{YY}.
Pour $0<\e(s)< {a'(s)}$,  on consid\`ere l'application
analytique~: 
  $$ \Psi : G \times G(s)_{\e(s)} \longrightarrow  G, \, \, 
                         (x , y)   \longmapsto 
                   xs y x^{-1}.$$
Alors, $\Psi$ est une submersion et son image,  not\'e $\W_G(s , \e(s))$ ou simple\-ment $\W(s ,
\e(s) )$ lorsque aucune confusion n'est possible,  est un  ouvert  $G$-semi-simple. De plus,
d'apr\`es le lemme \ref{yYY}, si $\e(s) < a(s)$, alors    
$\Psi$
induit un diff\'eomorphisme ${\bar \Psi}$ de l'ouvert $G
\times_{G(s)} G(s)_{\e(s)} $ du fibr\'e vectoriel $G \times_{G(s)}
G(s)$ sur $\W(s , \e(s) )$.
  On choisit  une mesure de Haar $d_Gx$ (resp.
$d_{G(s)}y$) sur $G$ (resp. $G(s)$). On munit $G /G(s)$ de la
mesure invariante $d_{G/G(s)} \dot x$.
 On a le r\'esultat suivant.
\begin{tho}\label{yPTo}
 Soit $0<\e(s)< {a'(s)}$
  tel que  $\bar \Psi$ soit  un diff\'eomorphisme.
 Soit $\Theta$ une fonction g\'en\'eralis\'ee $G$-invariante sur
$\W(s , \e(s) )$.
 Il existe une unique fonction g\'en\'eralis\'ee $\theta$,  $G(s)$-invariante, sur
$G(s)_{\e(s) }$,
telle que, pour toute fonction  $ \var \in C_c^\infty(\W(s , \e(s)))$, on ait
\begin{eqnarray}\label{opp}
  \int_{\W(s , \e(s))}\Theta (x) \var (x) d_Gx &=&
\int_{G/G(s)}|\det Adx{_\g}|_p  \{\int_{G(s)_{\e(s) }} \theta
(y)\var(xsyx^{-1})\times\nonumber \\&&
 \times|\det(Ad(sy)^{-1}
-1){_{(Ads-1)\g}}|_pd_{G(s)}y\}d_{G/G(s)}\dot{x}. 
\end{eqnarray}
 R\'eciproquement: si $\theta$ est une fonction
g\'en\'eralis\'ee, $G(s)$-invariante, sur $ G(s)_{\e(s) }$, la formule
(\ref{opp}) d\'efinit une fonction  g\'en\'eralis\'ee  $\Theta$,
$G$-invariante, sur l'ouvert $\W(s ,  \e(s))$.
\end{tho}
\end{subsection}
\end{section}
\begin{center}{\bf Entensions de repr\'esentations du groupe de Heisenberg}
\end{center}
\begin{section}{}\label{islam}
Dans cette section, nous  rappelons  certains r\'esultats  sur la
repr\'esentation de Weil et le groupe  m\'e\-tap\-lec\-ti\-que.
On peut consulter \cite{[M3]}, \cite{[We1]},   \cite{[Pe]} et
\cite{[L-V]}.
\begin{subsection}{}\label{ddk}
Soit $G$ un groupe localement compact, $H$ un sous-groupe
 ferm\'e de $G$, et $\pi$ une repr\'esentation unitaire de $H$ dans
 un espace de Hilbert ${\mathcal H}$. On note $\mbox{Ind}_H^G\pi$
 la repr\'esenta\-tion induite que  l'on   r\'ealise   dans l'espace des fonctions $\var $ sur $G$ \`a valeurs
 dans ${\mathcal H}$ mesurables telles que
 $$\var(xy) = \Delta_{H, G}(y)^{\frac{1}{2}}\pi(y^{-1})\var(x)\, ,
 \mbox{ pour } x \in G \, , y \in H
 \mbox{ et } \int_{G/H}||\var(x)||_{\mathcal H}^2 d_{G/H} \dot x < \infty.  $$
Le groupe $G$ agit sur cet espace par translations \`a gauche.
\end{subsection}
\begin{subsection}{}\label{represz}
 Soit  $\U$  un groupe alg\'ebrique unipotent d\'efini
 sur $k$  dont $ {\U}_k$ est  l'ensemble  des points rationnels. 
 On d\'esigne par  $\u$ l'alg\`ebre de Lie de  $ {\U}_k$.
 \`A toute  forme
 lin\'eaire  $u$
 sur  $ \u$, on associe une classe de
  repr\'esentations unitaires irr\'eductibles de $ {\U}_k$ par la m\'ethode des orbites de  Kirillov \cite{[Ki1]}.    Soit  $ \goth l$ une polarisation
    en  ${ u}$ et  on note ${\it L}
= \exp( \goth l)$ le sous-groupe unipotent de  ${{\U}_k}$ d'alg\`ebre de
Lie  $ \goth l$.
 On d\'efinit un caract\`ere  $ \chi_{u,  \goth l} $  de ${\it L}$, en posant~:
  $$  \chi_{u,  \goth l}(\exp X) =
 \varsigma( < u , X >) \mbox{,  pour tout } X \in  \goth l .$$
On consid\`ere    la
 repr\'esentation
 induite $ {\rm Ind}^ {{\U}_k}_{\it L}\chi_{u,  \goth l}$, not\'ee $ \pi_{u, \goth l}$, r\'ealis\'ee dans l'espace de Hilbert $\mathcal H_{ \goth l}$  d\'efini en num\'ero \ref{ddk}. Alors $ \pi_{u, \goth l}$ est
  une repr\'esentation  unitaire  irr\'eductible de $ {\U}_k$.  Sa classe  ne d\'epend pas de $\goth l$, on la note  $\pi_u$. Rappelons  que $\pi_u$ ne d\'epend que
 de l'orbite co-adjointe de $u$. 
 
\end{subsection}
\begin{subsection}{}\label{Perrin} On reprend les notations du num\'ero \ref{represz}. Soit $u\in \u^*$,  $ \goth l$ et $ \goth l'$ deux polarisations  en $u$. On note $F_{ \goth l\, ',  \goth l}$ l'op\'erateur d'entrelacement canonique de $\mathcal H_{ \goth l}$  dans $\mathcal H_{ \goth l\, '}$. Il est d\'efini par la propri\'et\'e suivante~:
il existe une mesure $d\dot y$, $L'$-invariante,  sur
$L'/L \cap L' $  telle que, pour tout $\alpha \in  \mathcal H_{ \goth l}$, continu \`a support compact modulo $L$, on ait : 
$$F_{ \goth l\, ',  \goth l} \alpha(x) = \int_{L'/L \cap L\, '} \alpha(x.y) \chi_{u,  \goth l\, '}(y) d\dot y,\,\,  x \in {\U}_k.$$
La mesure $ d\dot y$ sur  $L'/L \cap L'$ est d\'etermin\'ee par le fait que 
$F_{ \goth l\, ',  \goth l}$ est une isom\'etrie de $\mathcal H_{ \goth l}$ dans  $\mathcal H_{ \goth l\, '}$.
\end{subsection}

\begin{subsection}{}\label{ddn} Soit  $(V, B)$ un espace symplectique. On note   $H= V\times
k$ le groupe de Heisenberg associ\'e. La loi dans $H$ est donn\'ee
par
$$
(v,t).(v',t') = (v+v', t+t' + \frac{1}{2}B(v,v')),  \mbox{  pour
tout } v, v' \in V \mbox { et } t,t' \in k.
$$
Alors, $H$ est un groupe  unipotent
d'alg\`ebre de Lie  $\h = V \times k$, o\`u  le
crochet de Lie est donn\'e par
$$
[(v, t), (v', t')] =  (0,B(v,v')) \mbox{, pour tous } v, v' \in V \mbox
{ et } t,t' \in k.
$$
On note $E= (0, 1) \in \h$. C'est un \'el\'ement central de $\h$
 et on a~:  $\h= V \oplus kE$, o\`u on a identifi\'e $V$ avec le sous-espace vectoriel de $\h$
  constitu\'e des vecteurs $(v, 0)$,  $v \in V$.
L'application exponentielle $\exp : \h \longrightarrow H$ est donn\'ee par
$$
\exp(v+tE) = (v,t),
\mbox{ pour tout } v \in V\, , \,  t \in k.
$$
Le groupe symplectique $Sp(V)$ op\`ere  dans $\h$ par la
formule $$g.(v+ tE) = g.v + tE , \mbox{ pour tout } v \in V\, , \,
t \in k. $$
En composant avec l'application exponentielle, $Sp(V)$
op\`ere aussi dans $H$~:
$$g.\exp(v+ tE) =\exp( g.v + tE) , \mbox{ pour tout } v \in V\, ,
\,  t \in k. $$
Soit   $a_0 \in k^\times$ et  $E_{a_0}^*$ la
forme lin\'eaire sur $\h$ d\'efinie par : $$E_{a_0}^* (E)= a_0 \, ,
\,{ E_{a_0}}_{| V} = 0. $$
Soit \'egalement  $\l$ un lagrangien de $(V,
B)$. Alors $ \goth l = \l \oplus kE$ est une
polarisation en $E_{a_0}^*$ et 
 la
repr\'esentation  $(\pi_{E_{a_0}^*, \goth l},
 \mathcal H_{\goth l})$ de $H$ est unitaire irr\'eductible. Elle   v\'erifie la relation
 \begin{eqnarray} \label{Von}
\pi_{E_{a_0}^*, \goth l}(0, t)& =& \varsigma_{a_0}(t) \mbox{Id}_{\mathcal H_{\goth l}}
\, , \, \mbox{ pour tout } t \in k.
\end{eqnarray}
Pour $x \in Sp(V)$,    on pose  $${}^x\pi_{E_{a_0}^*, \goth l}(h)=
\pi_{E_{a_0}^*, \goth l}(x^{-1}. h), \quad h \in H.$$
Alors ${}^x\pi_{E_{a_0}^*, \goth l}$ est aussi une repr\'esentation  unitaire irr\'eductible de $H$  et 
v\'erifie   la formule (\ref{Von}).  D'apr\`es le th\'eor\`eme de Stone-Von Neumann, ${}^x\pi_{E_{a_0}^*, \goth l}$ est  \'equivalente \`a
   $\pi_{E_{a_0}^*,\goth l}$.  Soit $A(x)$  l'op\'erateur de 
$\mathcal H_{\goth  l}$ dans $\mathcal H_{x. \goth l}$ d\'efini par : $A(x) \alpha(y) = \alpha(x^{-1}.y)$, \, $\alpha \in \mathcal H_{\goth l}$, \, $y \in H$. On pose
\begin{eqnarray}\label{Lion-Perin}
S'_{E_{a_0}^*, {\l}}(x) &=& \frac{1}{||A(x)||}F_{\goth l,\,  x.\goth l} \circ A(x).
\end{eqnarray}
 Alors, on a :
$$ \qquad \qquad  S'_{E_{a_0}^*,  {\l}}(x)^{-1}
\pi_{E_{a_0}^*, \goth l}(y) S'_{E_{a_0}^*,  {\l}}(x)={}^x\pi_{E_{a_0}^*, \goth l}(y)
 \, , \, \mbox{ pour tout } \, y\in H.$$
De plus, $ S'_{E_{a_0}^*,  {\l}} : Sp(V)
\longrightarrow \mathcal U( \mathcal H_{\goth l})$ est une
repr\'esentation projective de   $Sp(V)$ dans  $\mathcal U(
\mathcal H_{\goth l})$, le groupe des op\'erateurs unitaires de $\mathcal H_{\goth l}$.  On d\'esigne par   $Mp(V)$ 
 l'ensemble des couples
$(x,\phi)$, o\`u $ x \in Sp(V) \mbox{ et } \phi$ est une fonction sur les
lagrangiens de $V$ v\'erifiant les propri\'et\'es (15) et (16) de
(\cite{[M3]}, Paragraphe~24.4). On munit $Mp(V)$ de  la loi de
composition interne introduite dans (\cite{[M3]}, Paragraphe~24.4). Il
r\'esulte de \cite{[We1]} que $Mp(V)$ est un groupe topologique,
localement compact, dont la premi\`ere projection fait  un
rev\^etement \`a deux feuillets de $Sp(V)$. Le groupe $Mp(V)$
s'appelle le groupe m\'etaplectique et 
           $$ 
          S_{E_{a_0}^* ,\l}   :  Mp(V)     \longrightarrow
          \mathcal U(\mathcal H_{\goth l}), \, \, 
                           (x , \phi)  \longmapsto   \phi(\l)
                           S'_{E_{a_0}^*,\l}(x) $$
 est une repr\'esentation fid\`ele.  Elle ne d\'epend pas du
lagrangien $\l$, on la note simplement $ S_{E_{a_0}^*}$. C'est la
repr\'esentation m\'etaplectique de $Mp(V)$. Le r\'esultat suivant
est bien connu.
\begin{tho}
\

 $i)$ Les diff\'erents groupes $Mp(V)$, quand $a_0$
d\'ecrit $k^\times$, sont canoniquement  isomorphes.

 $ii)$ La repr\'esentation m\'etaplectique $S_{E_{a_0}^*}$ ne
d\'epend de  $a_0$ que par sa classe  dans~$k^\times  /
{k^\times} ^2$.
\end{tho}
\end{subsection}
\begin{subsection}{Une fonction sur le groupe m\'etaplectique.
}\label{Kllov}
 Soit $ x \in Mp(V)$ et $x =
\tilde {s}.x_u$ sa d\'ecomposition de Jordan, $\tilde {s}= (s, \psi_s)$ \'etant
la partie semi-simple de $x$ et $x_u$ \'etant la partie unipotente
de $x$.  Soit \'egalement
  $ (1-s).V  = W_1 \oplus W_2$ une d\'ecomposition   en somme directe orthogonale par
rapport \`a   $B$, o\`u   $W_1$  est un  sous-espace symplectique $s$-invariant poss\'edant un
lagrangien $\l_1$ stable sous l'action de $s$ et $W_2$  est un  sous-espace symplectique $s$-invariant poss\'edant
un lagrangien $\l_2$ tel que $(s.\l_2)\cap \l_2 = 0$ (voir \cite{[M3]}, Lemme 27). On consid\`ere  la  forme quadratique $Q_{s , \l_2}$
sur $ \l_2$~:  $$Q_{s , \l_2}(v) = B(v,(s^{-1} - 1 )^{-1}. v), \quad \mbox{ pour tout }
 v \in \l_2.$$ On note
$\gamma_{a_0}(Q_{s , \l_2}) = \gamma(a_0Q_{s , \l_2})$ le nombre complexe de module $1$
associ\'e \`a la  forme quadratique  $a_0Q_{s , \l_2}$ par la formule (\ref{quadratique}). On pose
\begin{eqnarray}
\label{ammar} \Phi_{a_0}( x) = \Phi_{a_0}(\tilde s) = \psi_s(\l_0
+\l_1 +\l_2) \overline{\gamma_{a_0}(Q_{s , \l_2})},
\end{eqnarray}
$\l_0$ \'etant un lagrangien de $V(s)$. D'apr\`es (\cite{[M3]},
Corollaire 27) le nombre $ \Phi_{a_0}( x)$, est bien d\'efini,
ne d\'epend ni  du choix de la d\'ecomposition  $ (1-s).V= W_1 \oplus W_2$ ni du choix
 des lagrangiens $\l_0, \l_1, \l_2$
de  $V(s)$,  $W_1$ et $ W_2$ res\-pectivement. La formule (\ref{ammar}) d\'efinit une fonction $ \Phi_{a_0}$ sur  $Mp(V)$,  invariante par conjugaison par les \'el\'ements de $Mp(V)$.
\end{subsection}
\begin{subsection}{Une formule du caract\`ere.}\label{Air}
   On se donne un \'el\'ement
semi-simple $s$ de $Sp(V)$. On choisit  $\tilde s = (s, \psi_s)
$ un rel\`evement de $s$ dans $ Mp(V)$. On consid\`ere l'action de  $H$ sur lui-m\^eme d\'efinie par 
 $$
 x*y = xys(x^{-1}), \mbox{ pour tout } x, y \in
H.
$$
L'application $\Psi : H \times H(s) \longrightarrow H$, $(x, y) \longmapsto x*y$ est   une submersion
surjective.
 On pose
  $$\Lambda_{\tilde s}(\var d_Hx) =
\mbox{Tr}(S_{E_{a_0}^*}(\tilde{s})\pi_{E_{a_0}^*}(\var d_Hx)), \mbox{ pour tout } \var \in C_c^\infty(H),   $$
 $d_Hx$ \'etant une mesure de Haar sur $H$. Ceci a un sens puisque $\pi_{E_{a_0}^*}$
est admissible et  $S_{E_{a_0}^*}(\tilde{s})$ est un op\'erateur
born\'e.  Alors $\Lambda_{\tilde s}$ est  fonction g\'en\'eralis\'ee sur $H$,   $*$-invariante. 
On d\'esigne par $\lambda_{\tilde s}$ la fonction g\'en\'eralis\'ee sur $H(s)$ associ\'ee \`a $\Lambda_{\tilde s}$ par le th\'eor\`eme
\ref{ysemis}. Alors $\lambda_{\tilde s}$ d\'etermine enti\`erement $\Lambda_{\tilde s}$, et elle  est $H(s)$-invariante. 
On note  $\O_{E_{a_0}^*}$  l'orbite co-adjointe de $E_{a_0}^*$ sous l'action de  $H$ et  $\O_{E_{a_0}^*, s} = \O_{E_{a_0}^*} \cap \h^*(s)$ l'ensemble des points fixes de $s$ dans $\O_{E_{a_0}^*}$. Alors $\O_{E_{a_0}^*, s}$ est une $H(s)$-orbite.
\begin{tho}\label{yAir}
 Pour tout $\beta \in C_c^\infty(H(s))$, on a :
\begin{eqnarray*}
 \lambda_{\tilde s} (\beta d_{H(s)}y)
 &= & \Phi_{a_0}(\tilde s)|\det(1 - s)_{(1-s)V}|_{p}^{-\frac{1}{2}}
 \int_{\O_{E_{a_0}^*,s}}(\beta \circ
\exp d_{\h(s)}Y)\,  \widehat{}_{\h(s)} \, (l)
d\mu_{\O_{E_{a_0}^*,s}}(l),
\end{eqnarray*}
o\`u $d_{H(s)}y$ est une  mesure de Haar  sur $H(s)$ et   $d_{\h(s)}Y$ est la mesure de Haar sur $\h(s)$ tangente \`a $d_{H(s)}y$.
\end{tho}
\noindent Lorsque $a_0 = 1$, le th\'eor\`eme ci-dessus est le  th\'eor\`eme 28 de \cite{[M3]}.
\begin{subsubsection}{}\label{aca}
 Dans la suite, on  utilise les notations suivantes :
$$ V_s = (1 - s)V, \,\,  \h_{2,s} =  V_s + kE, \, \,  H_{2,s} =
\exp(\h_{2,s}), \, \,  E^*_{2,s}= E^*_{a_0|\h_{2,s}}, $$ $s_2 =s_{|V_s}$, la
restriction de $s$ \`a $V_s$,  et $\tilde{s_2}$
un rel\`evement de $s_2$ dans $Mp(V_s)$.
 On se donne un lagrangien $\l_2$ de
 $V_s$. On note $\goth l_2 =\l_2 \oplus kE $. C'est une polarisation en $E^*_{2,s}$.
On d\'esigne par  $(\pi_{E^*_{2,s}}, \mathcal H_2)$ la repr\'esentation
unitaire irr\'eductible de $H_{2,s}$
 associ\'ee \`a $E^*_{2,s}$ par la
m\'ethode des orbites de Kirillov et r\'ealis\'ee dans l'espace de Hilbert $
\mathcal H_2=  \mathcal H_{\l_2}$, et par $S_{E^*_{2,s}}$ la
repr\'esentation m\'etaplectique  de $Mp(V_s)$ associ\'ee au
caract\`ere $\varsigma_{a_0}$ et r\'ealis\'ee dans le m\^eme
espace de
 Hilbert $\mathcal H_2$.

\noindent D\'esignons par $ \nu_1, \ldots,  \nu_l$  les valeurs propres
distinctes de $s$ dans $V_{\bar k}={\bar k} \otimes_k  V $. On
pose~: $$a'(s)_{Sp(V)} = \inf_{ \nu_i\neq 1}\{| \nu_i -1|_{\bar k},
| \nu_i^{-1} -1|_{\bar k}, 
\,  a\}.$$ Soit   $0<\e(s)<a'(s)_{Sp(V)}$. 
On note $$\W(\tilde{s_2} , \e(s))= \{\tilde x\tilde{s_2}\tilde y
{\tilde x}^{-1}, \, \, \tilde y \in  Mp(V_s)(\tilde{s_2})_{\e(s)},
\, \, \tilde x \in Mp(V_s) \}.$$
 Pour $\tilde z
\in \W(\tilde{s_2} , \e(s))$ et   $\alpha \in C_c^\infty(V_s)$, on consid\`ere  l'op\'erateur, de
rang fini, agissant dans   $\mathcal H_2$~:
\begin{eqnarray} \label{op}
 J_{\alpha}(\tilde z)= S_{E^*_{2,s}}(\tilde z)
 \int_{V_s}\alpha(v)\pi_{E^*_{2,s}}(\exp({\tilde z}^{-1}.v)\exp(-v))d_{V_s}v,
\end{eqnarray}
 $d_{V_s}v$ \'etant   une mesure de
Haar sur $V_s$.
On a le r\'esultat suivant. 
\begin{pro} \label{yaca}
On suppose que  $ \tilde z $ est  semi-simple.
 Pour tout $\alpha \in C_c^\infty(V_s)$  telle que
$\int_{V_s}\alpha(v) d_{V_s}v =1$, on a :
\begin{eqnarray}
\mbox{Tr}(J_{\alpha}(\tilde z)) &=&\Phi_{a_0}(\tilde z)
|\det(1-\tilde z)_{V_s}|_p^{-\frac{1}{2}}.
\end{eqnarray}
\end{pro}
\begin{dem} C'est une cons\'equence de  (\cite{[M3]}, Proposition 29.3.3).
\end{dem}
\end{subsubsection}
\begin{subsubsection}{}\label{exten} Dans la suite, on a besoin du r\'esultat suivant.
\begin{lem}\label{yexten}
 Soit $\alpha \in C_c^\infty(V_s)$. Il existe un sous-groupe ouvert $K$ de $GL(V_s)$ tel que
$$\alpha(x.v) = \alpha(v), \, \,  \mbox{ pour tous }  v \in V_s, \, x \in
K. $$
\end{lem}
\begin{dem}
Puisque $\alpha \in C_c^\infty(V_s)$, il existe $r$ un r\'eseau de
$V_s$, $v_1, \ldots, v_m \in V_s$, et $c_1, \ldots, c_m \in \cit$
tels que $$\alpha = \sum_{1 \leq j \leq m} c_j1_{v_j+ r}. $$
L'ensemble $$K= \bigcap_{1 \leq j \leq m} \{x \in GL(V_s) , \, \,
x(v_j)-v_j \in r, \, x (r) = r\} $$ est  un sous-groupe ouvert  de
$GL(V_s)$ et   remplit les conditions voulues.
\end{dem}
\begin{pro} \label{yyexten}
Soit $\alpha \in C_c^\infty(V_s)$ et  $\tilde {y_0} \in Mp(V_s)(\tilde{s_2})_{\e(s)}$.  Il
existe un voisinage $\mathcal V_{ \tilde y_0}$ de $ \tilde y_0 $
dans $Mp(V_s)(\tilde{s_2})$ contenu dans
$Mp(V_s)(\tilde{s_2})_{\e(s)}$
 tel que
$$ J_{\alpha}(\tilde s_2 \tilde y) =J_{\alpha}(\tilde s_2 \tilde
y_0), \, \,   \mbox{ pour tout } \,  \tilde y \in \mathcal V_{\tilde y_0}.
$$
\end{pro}
\begin{dem} D'apr\`es  le 
lemme \ref{ykhem1}, la projection canonique  de $Mp(V_s)$ dans
$Sp(V_s)$ induit un diff\'eomorphisme $p_s : Mp(V_s)(\tilde{s_2})_{\e(s)} \longrightarrow Sp(V_s)(s_2)_{\e(s)|4|_p}$.
  Pour $y \in Sp(V_s)(s_2)_{\e(s)|4|_p}$, on note $\tilde y = p_s^{-1}(y)$. On a : 
$$J_{\alpha}(\tilde s_2 \tilde y) =S_{E^*_{2,s}}(\tilde s_2 \tilde y)J'_{\alpha}(s_2y), $$
o\`u
$$J'_{\alpha}(s_2y)=\int_{V_s}\alpha(v)\pi_{E^*_{2,s}}(\exp((s_2y)^{-1}.v)\exp(-v))d_{V_s}v.
$$  Un calcul imm\'ediat donne 
\begin{eqnarray*}J'_{\alpha}(s_2 y)
&=&\int_{V_s}\alpha(v)\pi_{E^*_{2,s}}(\exp(((s_2y)^{-1}-1).v))\varsigma_{a_0}(-\frac{1}{2}
B(((s_2y)^{-1}-1).v,v))d_{V_s}v\\
 &=& |\det
((s_2y)^{-1}-1)_{V_s}|_p^{-1}\\&&
 \int_{V_s}\alpha(((s_2y)^{-1}-1)^{-1}.v)
\pi_{E^*_{2,s}}(\exp(v))\varsigma_{a_0}(-\frac{1}{2}B(v,((s_2y)^{-1}-1)^{-1}
v))d_{V_s}v.
\end{eqnarray*}
D'apr\`es le lemme ci-dessus, il existe un sous-goupe ouvert
$K$ de $GL(V_s)$ tel que $$\alpha(((s_2y_0)^{-1}-1)^{-1}x.v) =
\alpha(((s_2y_0)^{-1}-1)^{-1}v), \, \, \forall \,  v \in V_s, \,\forall \, x
\in K. $$
Il en r\'esulte qu'il existe un sous-groupe
compact ouvert $K'$ de  $Sp(V_s)(s_2)$ contenu dans
$Sp(V_s)(s_2)_{\e(s)|4|_p}$ tel que l'on ait  $y_0K'
\subset Sp(V_s)(s_2)_{\e(s)|4|_p}$ et 
$$\alpha(((s_2y)^{-1}-1)^{-1}.v) = \alpha(((s_2y_0)^{-1}-1)^{-1}.v), \forall \,  y \in y_0K', \forall \,  v \in V.$$
 Maintenant, soit $r_s$ un r\'eseau de $V_s$ contenant
$((s_2y_0)^{-1}-1).\mbox{supp}(\alpha)$. L'application $\eta
: ( y, v) \longmapsto \varsigma_{a_0}(-\frac{1}{2}B(v,((s_2y_0y)^{-1}-1)^{-1}
v))$ 
de
$K' \times r_s$ dans $\cit^\times$, est localement
constante.  On en
d\'eduit qu'il existe un sous-groupe compact ouvert $K_0$
contenu dans $K'$ tel que l'on ait  $$\eta(y,v)=  \eta(1,v), \mbox{ pour tout } (y,v) \in K_0 \times r_s .$$ 
Si bien que l'on a : 
$$J'_{\alpha}(s_2  y_0y)=J'_{\alpha}( s_2y_0), \,\mbox{ pour tout } y \in
K_0. $$
 On note alors  $\tilde {K_0} = p_s^{-1}(K_0)$.

\noindent  D'autre part,
pour $m \in \N$ et   $w \in \omega^m r_s$,  on a :
\begin{eqnarray*}
\pi_{E^*_{2,s}}(\exp(w))J'_{\alpha}( s_2 y_0) &=& |\det
((s_2y_0)^{-1}-1)_{V_s}|_p^{-1}
 \int_{V_s}\alpha(((s_2y_0)^{-1}-1)^{-1}.(v-w))\\&& \times
\pi_{E^*_{2,s}}(\exp(v))
\varsigma_{a_0}(-\frac{1}{2}B(v-w,((s_2y_0)^{-1}-1)^{-1}
(v-w)))\\&& \times \varsigma_{a_0}(\frac{1}{2}B(w, v))d_{V_s}v.
\end{eqnarray*}
Ainsi,  si $m$ est  suffisamment grand, on a : 
$$\pi_{E^*_{2,s}}(\exp(w))J'_{\alpha}( s_2 y_0) =J'_{\alpha}(s_2 y_0)\, ,
\, \mbox{ pour tout } w \in \omega^m r_s. $$
 Posons $K_m= \exp(\omega^mr_s + \O a_0E)$. Si $m$ est  suffisamment grand, $K_m$ est un
 sous-groupe compact ouvert de $H_{2,s}$.
 Comme $\pi_{E^*_{2,s}}$ est admissible, $\mathcal H_2^{K_m}$ est de dimension
finie. Notons $M=Sp(V_s)(r_s)$ et $M^{V_s} $ son image
r\'eciproque  dans $Mp(V_s)$.
 Alors,   $ (S_{E^*_{2,s}})_{|M^{V_s}}$ laisse
 invariant
  $\mathcal H_2^{K_m}$  et elle
   induit une repr\'esentation continue du groupe
  $M^{V_s}$ dans
 $\mathcal H_2^{K_m}$.  Puisque
$GL(\mathcal H_2^{K_m})$ n'a pas de sous-groupes non triviaux
arbitrairement petits, il existe un sous-groupe ouvert
$\tilde{ M}$ de $M^{V_s}$ tel que
 $ (S_{E^*_{2,s}})_{|\tilde{ M}} = \mbox{Id}_{\mathcal
 H_2^{K_m}}$.
 L'ensemble $\tilde{y_0}.(\tilde{ M} \cap \tilde {K_0})$ est un voisinage de
 $\tilde{y_0}$ ayant les propri\'et\'es voulues.
\end{dem}
\end{subsubsection}
\begin{subsubsection}{}\label{extent} On reprend  les notations des sections \ref{aca} et \ref{exten}.
Soit $\alpha \in C_c^\infty(V_s)$  tel que $\int_{V_s}\alpha(v)
d_{V_s}v =1$. On d\'efinit une  application : $$T_\alpha :
   \W(\tilde{s_2} , \e(s))\longrightarrow \cit, \, \,  \tilde z \longmapsto |\det(1- \tilde
z )_{V_s}|_p^{\frac{1}{2}}\mbox{Tr}(J_{\alpha}(\tilde z)).$$
\begin{pro}\label{yextent}
On a~:
\begin{eqnarray}\label{ramd}
T_\alpha = {\Phi_{a_0}}_{|\W(\tilde{s_2} , \e(s))}.
\end{eqnarray}
Si  $\e(s)$ est suffisamment petit, alors
\begin{eqnarray}\label{ramdda}
\quad T_\alpha( \tilde  z) = T_\alpha(\tilde{s_2}), \forall \,  \tilde z \in \W(\tilde{s_2} , \e(s)) , \forall
a_0 \in k^\times.
\end{eqnarray}
\end{pro}
\begin{dem} On a~:
\begin{eqnarray}\label{ramdan}
T_\alpha(  \tilde x \tilde z  {\tilde x}^{-1})&= &
T_{\alpha^{\tilde x}}(\tilde z), \, \tilde z \in
\W(\tilde{s_2} , \e(s)), \, \tilde x \in Mp(V_s),
\end{eqnarray}
o\`u $\alpha^{\tilde x}$ est l'\'el\'ement de  $C_c^\infty(V_s)$ d\'efini par :  $\alpha^{\tilde x}(v)= \alpha({\tilde x}.v)$, pour tout $v
\in V_s$.
Il r\'esulte du lemme  \ref{yexten}, de la proposition \ref{yyexten},  et de  l'\'egalit\'e (\ref{ramdan})
que  l'application $T_\alpha $ est localement constante.
D'apr\`es la proposition \ref{yaca}, l'\'egalit\'e (\ref{ramd}) est vraie sur les  \'el\'ements semi-simples de
$\W(\tilde{s_2} , \e(s))$.
 Mais l'ensemble des \'el\'ements semi-simples de $Mp(V_s)$ qui sont contenus dans
 $\W(\tilde{s_2} , \e(s))$ est dense dans $\W(\tilde{s_2} , \e(s))$.  On en d\'eduit que l'application $T_\alpha $ ne d\'epend pas de $\alpha$ pourvu que $\int_{V_s}\alpha(v) d_{V_s}v =1$. Il en r\'esulte que
 $T_\alpha$ est invariante par conjugaison par les \'el\'ements de $Mp(V_s)$.  Soit $\tilde  z \in \W(\tilde{s_2} , \e(s))$.  Ecrivons $ \tilde  z = \tilde  z_s \tilde  z_u$, o\`u $\tilde  z_s$ est la partie semi-simple de $\tilde  z$ et $\tilde  z_u$ est la partie unipotente de $\tilde  z$.  En utilisant la proposition
 \ref{yreductive}, on~a~:
 $$T_\alpha( \tilde  z) = T_\alpha( \tilde  z_s)= \Phi_{a_0}( \tilde  z_s)=\Phi_{a_0}( \tilde  z) .$$
Pour d\'emontrer la deuxi\`eme assertion, on remarque    que
l'ensemble  $$ A: = \{\tilde z \in \W(\tilde{s_2} , \e(s)) \mbox{ tel
que } T_\alpha( \tilde z) = T_\alpha(\tilde{s_2})
\}$$ est un voisinage ouvert, invariant, de $\tilde{s_2}$ dans $Mp(V_s)$.  On en d\'eduit qu'il existe un voisi\-nage ouvert, $Mp(V_s)(\tilde{s_2})$-invariant,   $\tilde U$   de $1$ dans   $Mp(V_s)(\tilde{s_2})$ contenu dans $Mp(V_s)(\tilde{s_2})_{\e(s)}$ tel que $\tilde{s_2}.\tilde U \subset A$. 
D'apr\`es le lemme \ref{yyreductive}, il existe $0 < \e' \leq \e(s)$ tel que 
$Mp(V_s)(\tilde{s_2})_{\e'}\subset \tilde U $.  Si bien que  $\W(\tilde{s_2} , \e') \subset A$.
\end{dem}
\end{subsubsection}
\begin{subsubsection}{} \label{exc}
Dans cette section, on se donne un $k$-espace vectoriel $W$ de
dimension finie et $x$ un automorphisme semi-simple  de $W$. Si $V$ est  un  sous-espace vectoriel stable par $x$ 
et $B$ une forme symplectique  sur $V$ fix\'ee par $x$, 
 on note $V_x =(1-x).V $ et
$\tilde{ x_{|V_x}}$ un rel\`evement de $x_{|V_x}$ dans $Mp(V_x)$.
On a le r\'esultat suivant.
\begin{pro}  \label{yexc} Il existe $\e(x) >0$ tel que, pour tout sous-espace
vectoriel $V$ stable par $x$ et pour toute forme symplectique sur $V$ fix\'e par $x$,   on ait
\begin{eqnarray}\label{ramadann}
\quad  {\Phi_{a_0}}_{| \W( \tilde{ x_{|V_x}}, \e(x))}= \Phi_{a_0}( \tilde{ x_{|V_x}}) ,  \mbox{ pour tout } \, 
a_0 \in k^\times.
\end{eqnarray}
\end{pro}
\begin{dem}
On d\'esigne par $\mathcal B_i$ l'ensemble des sous-espaces
symplectiques $(V, B)$,  de dimension $i$, de $W$. Le groupe
$GL(W)$  op\`ere \`a gauche  dans $\mathcal  B_i$ par : si $g \in
GL(W)$ et $(V, B) \in \mathcal B_i $ alors  $g.(V, B)= (g.V,
B^{g})$, o\`u $B^{g}(v, w)= B( g^{-1}v,  g^{-1}w)$.
 Choisissons un
\'el\'ement $(V_0, B_0) \in \mathcal B_i$, fix\'e par $x$ (s'il en
existe) et posons $K=  GL(W)((V_0, B_0))$ le stabilisateur de $(V_0, B_0)$ dans $ GL(W)$.
Alors $K$ est un sous-groupe alg\'ebrique  (ferm\'e) de  $GL(W)$. Comme l'action de $GL(W)$ dans $\mathcal B_i$ est transitive alors $\mathcal B_i$
s'identifie canoniquement   avec l'espace homog\`ene
$GL(W)/K$. Consid\'erons l'action  par translations \`a gauche  de
$GL(W)$ sur $GL(W)/K$. Alors, l'application $\goth
T : GL(W)/K \longrightarrow \mathcal B_i$,
$g.K\longmapsto g.(V_0, B_0)$ est une bijection,
$GL(W)$-\'equivariante. Soit $H = GL(W)(x)$, le centralisateur de
$x$ dans  $GL(W)$. D'apr\`es (\cite{[Ri]}, Th\'eor\`eme A), l'ensemble
de points fixes $(GL(W)/K)^x$  pour l'action de $x$ dans
$GL(W)/K$   est une r\'eunion
finie de $H$-orbites. L'image de $(GL(W)/K)^x$
par $\goth T $, not\'ee $\mathcal A_i$, est l'ensemble de
sous-espaces symplectiques $(V, B)$ pour lequel $V$ est stable par $x$ et   $x_{|V} \in Sp(V, B)$,
de dimension $i$ de $W$. Il en r\'esulte que  $\mathcal A_i$ est
une r\'eunion finie de $H$-orbites. Prenons, maintenant, une
$H$-orbite $\O_i$ dans $\mathcal A_i$. Par transport de structure et la proposition 
\ref{yextent},
on peut choisir $\e_{\O_i}(x)
>0$ tel que l'on ait la formule (\ref{ramadann}) pour tout $(V, B)
\in \O_i$. Il s'ensuit qu'il existe $\e_i(x)
>0$ tel que l'on ait la formule (\ref{ramadann}) pour tout $(V, B)
\in \mathcal A_i$. Il suffit alors de prendre  $\e(x) = \inf\{\e_i(x) \mbox{ tel que } 0 \leq i \leq \dim W \mbox{ et } \mathcal A_i \neq \emptyset\}$. 
\end{dem}
\end{subsubsection}
\end{subsection}
\end{section}

\medskip

\begin{center}{ \bf Formule du  caract\`ere.}
\end{center}
\begin{section}{}\label{intro} Soit $(G, F, {\G})$ un groupe presque alg\'ebrique d'alg\`ebre de Lie $\g$.
Pour
chaque forme lin\'eaire $g$ sur $\g$, on
consid\`ere
l'extension m\'etaplectique
$G(g)^\goth g $  correspondant \`a  l'action  de $G(g)$ dans
$(\goth g, \beta_g) $ dont l'\'el\'ement central non trivial du noyau est
not\'e $ (1,-1)$ (voir le num\'ero \ref{ddp}). On note $\chi_g$  le caract\`ere  de ${}^uG(g)^\goth g$,
radical unipotent de $G(g)^\g$, d\'efini par la formule :
$\chi_g(\exp(X))= \varsigma(<g, X>), \, X \in {}^u\g(g)$. On
note $Y_G(g)$ l'ensemble des classes  des repr\'esentations
unitaires irr\'eductibles $\tau$ de $G(g)^\g$ dont la
restriction \`a ${}^uG(g)^\goth g$ est multiple de $\chi_g$ et telle que
$\tau (1,-1)= -\mbox{Id}$. Si $g$ est de type unipotent (voir num\'ero \ref{kh24}) et $\tau$ est un
\'el\'ement de $Y_G(g)$, M. Duflo  a associ\'e au couple $(g,
\tau)$ une classe de repr\'esentations unitaires irr\'eductibles
$\pi_{g , \tau}$ de $G$ telle que :
\begin{itemize}
\item pour tout automorphisme $a$ de $(G , F , \G)$, on ait~:
$\pi_{ag , {}^a\tau} =  {}^a\pi_{g , \tau}$ (o\`u
 $ {}^a\tau = \tau
\circ a^{-1}$ et ${}^a\pi_{g , \tau} =\pi_{g , \tau}\circ a^{-1}$
)~;

\item  soit $\tau' \in Y_G(g)$ ;  $\pi_{g , \tau}$ et  $\pi_{g ,
\tau'}$ sont \'equivalentes si et seulement si $\tau$ et $\tau'$
sont \'equi\-va\-lentes~;

\item soit $g'$ une forme de type unipotent
et $\tau' \in  Y_G(g')$. On suppose que $g$ et $g'$ ne sont pas
dans une m\^eme $G$-orbite, alors $\pi_{g , \tau}$ et $\pi_{g' ,
\tau'}$ sont in\'equivalentes.
\end{itemize}
D\'esignons par $Y_G$ l'ensemble des couples $(g, \tau)$, o\`u
$g$ est de type unipotent et $\tau \in Y_G(g)$. Le groupe $G$
op\`ere naturellement dans $Y_G$ et la corres\-pondance ci-dessus
induit une bijection de $G\backslash Y_G$ sur $\widehat{G}$, l'ensemble des classes
des repr\'esentations unitaires irr\'eductibles de $G$.
 Cette description s'appelle la m\'ethode des orbites de
Kirillov-Duflo pour la construction des repr\'esentations
unitaires irr\'eductibles de $G$. Elle  permet d'obtenir  tous les \'el\'ements  de $\widehat{G}$.

Si $\pi_{g , \tau}$ est admissible, on note $\Theta_{g, \tau}$ son caract\`ere~;   c'est la
fonction g\'en\'eralis\'ee sur $G$ d\'efinie par
\begin{eqnarray}\label{anniver}
\Theta_{g, \tau}(\var d_Gx) &=& \mbox{Tr}(\pi_{g, \tau}(\var
d_Gx)), \,  \var \in C_c^\infty(G),
\end{eqnarray}
o\`u $d_Gx$ est une mesure de Haar sur $G$.

\noindent Dans \cite{[M3]}, nous avons  d\'emontr\'e que si l'orbite
co-adjointe $\O_g$ de $g$  est ferm\'ee dans $\g^*$ et si $\tau$ est de dimension finie  alors
$\pi_{g , \tau}$ est admissible  et que $\Theta_{g, \tau}$ est
donn\'e, au voisinage de  chaque \'el\'ement semi-simple $s$ de
$G$, par une formule \`a la Duflo-Heckman-Vergne, en termes de
transform\'ees de Fourier de l'ensemble des points fixes $\O_{g,
s}$ de $s$ dans $\O_g$, du caract\`ere de  $\tau$, et d'une
fonction  bien d\'efinie, constante sur chaque $G(s)$-orbite sur
$\O_{g, s}$. N\'eanmoins le voisinage semi-simple  sur lequel nous avons   la
validit\'e de la formule \'etabli dans \cite{[M3]} d\'epend de  la
repr\'esentation  $\pi_{g , \tau}$.
Dans cette partie du pr\'esent travail nous allons   rem\'edier  \`a cette
 d\'ependance pour les groupes r\'esolubles  presque connexes et dans le cas o\`u $p$,
  la caract\'eristique r\'esiduelle de $k$,  est diff\'erente de
  $2$.

\begin{subsection}{Une partition de l'ensemble des orbites co-adjointes.} \label{kh24} 
On suppose d\'esormais  que $\g$ est   r\'esoluble. 
 Suivant Duflo (voir
\cite{[Du]}), une forme  lin\'eaire $g$ sur $\g$ est dite de type
unipotent si $\r_g$,  facteur r\'eductif de $\g(g)$,  est contenu
dans $\ker g$. D'apr\`es \cite{[Du]}, l'application qui \`a $g \in \g^*$
associe $ g_{|{}^u \g}$ induit une bijection entre l'ensemble
des $G$-orbites de type unipotent dans $\g^*$ et l'ensemble des
$G$-orbites dans $({}^u \g)^*$.

\noindent   Soit $g \in \g^*$ de
type unipotent. On note $L(g)$ l'ensemble des formes lin\'eaires
$\lambda$ sur $\g(g)$ dont la restriction \`a ${}^u\g(g)$ est
\'egale \`a celle de $g$. Si   $\r_g$  est un facteur r\'eductif
de $\g(g)$ alors
 l'application $\lambda  \in (\g(g))^* \longmapsto \lambda_{|\r_g} $ r\'ealise une bijection
entre $L(g)$ et  $\r_g^*$. On note $\goth D$ l'ensemble des
couples $(g, \lambda)$ o\`u $g$ est de type unipotent et o\`u
$\lambda \in L(g)$. Le groupe $G$ op\`ere naturellement sur $\goth
D$. Si $(g, \lambda) \in \goth D$ et $\b$ est une sous-alg\`ebre
de type fortement unipotent relativement \`a $g$ (i.e.  $\b$ est
co-isotrope par rapport \`a $g$, alg\'ebrique, et $\b = \g(g)
+{}^u\b$), on consid\`ere une forme lin\'eaire $f$ sur $\g$ telle
que
\begin{eqnarray}
f_{|{{}^u\b}} = g_{|{{}^u\b}} \mbox { et } f_{|{\g(g)}} = \lambda.
\end{eqnarray}
Le r\'esultat suivant est d\^u \`a M. Duflo \cite{[Du]}.
\begin{pro}  On a~ :

$ i)$
L'orbite $G.f$ de $f$ sous l'action de $G$ ne d\'epend pas des
choix de $\b$ et $f$. On la note~$\O_{g,\lambda}$.

$ ii)$
L'application  $(g, \lambda) \longmapsto \O_{g,\lambda}$ induit
une bijection de $ G\backslash{\goth D}$ sur $ G\backslash \g^*$.
\end{pro}
\end{subsection}
\begin{subsection}{}\label{meet}
On se donne  $g \in \g^*$ de type unipotent. On pose $\b_g = \g(g) + {}^u\g$.  C'est une
sous-alg\`ebre de type fortement unipotent relativement \`a $g$.
C'est  la sous-alg\`ebre acceptable canonique associ\'ee \`a $g$
(voir \cite{[Du]}, Chapitre 1). Soit $\lambda \in
L(g)$ et $f \in
\g^*$ tel que $$f_{|{}^u\g} = g_{|{}^u\g} \mbox { et }
f_{|{\g(g)}} = \lambda.$$ Alors, 
on a~: $\O_{g,\lambda} = G.f$.  

Soit  $\u$   un id\'eal ab\'elien de $\g$,
$G$-invariant,  et contenu dans ${}^u\g$.
 On note $U$ le
sous-groupe unipotent de $G$ d'alg\`ebre de Lie $\u$.
 On pose $$ u =g_{|\u}, \,  \h = \g(u),\,  h= g_{|\h}, \, H= G(u).$$ Alors, $\h$ est une
sous-alg\`ebre de Lie, alg\'ebrique,  de $\g$. D'apr\`es
(\cite{[Pu2]}, p. 500-501),  $\h(h)= \g(g) +\u$, si bien que  le
radical unipotent de $\h(h)$ est ${}^u\g(g) +\u$.

 En consid\'erant $\lambda$ comme un \'el\'ement de $L(h)$, en la prolongeant par $g_{|\u}$ sur $\u$,
  alors $\O_{h,{\lambda}} $ est l'orbite
co-adjointe de $f_{|\h}$ sous l'action de $H$.

Notons $r : \g^* \longrightarrow \h^*$  (resp. $r' : \g^*
\longrightarrow \u^*$ )
 l'application restriction de $\g^*$ dans $\h^*$ (resp. dans $\u^*$).
 Alors, l'image r\'eciproque de $\{u\}$ dans  $\O_{g,{\lambda}}$ par $r'$  est $H.f$.
Il en r\'esulte que si $\O_{g,{\lambda}}$ est ferm\'ee dans $\g^*$
alors $H.f$ l'est aussi. Comme $U.x.f = x.f + \h^\perp$, pour tout
$x \in H$, $H.f$ est satur\'ee par rapport \`a $r$. Si bien que si
$\O_{g,{\lambda}}$ est ferm\'ee, il en est de m\^eme pour
  $\O_{h,{\lambda}} =
r(H.f) $   dans $\h^*$. 

 On munit   $\O_{g,{\lambda}}$ (resp.
$\O_{h,{\lambda}} $) de la mesure de Liouville
$d\mu_{\O_{g,{\lambda}}}$ (resp. $d\mu_{\O_{h,{\lambda}}}$). On
choisit  $d_Gx$ (resp.  $d_Hy$) une mesure de Haar sur $G$ (resp.
$H$).  On note $d_\g X$ (resp. $d_\h Y$) la mesure de Haar sur
$\g$ (resp. $\h$)  tangente \`a $d_Gx$ (resp.  $d_Hy$).
 On consid\`ere la mesure invariante $d_{G/H}\dot x = d_Gx/d_Hy$ sur $G/H$.
 On munit  l'espace vectoriel $\h^\bot$
 de la mesure de Haar $d_{\h^\bot}t$ duale de  ${d_\g X}/{d_\h Y}$. On a le
r\'esultat suivant qui   se d\'emontre de la m\^eme fa\c{c}on que le  lemme 7 de \cite{[Ki]}.
\begin{pro}\label{ymeet}
 Avec les notations ci-dessus on a, pour toute
 $  \var$ mesurable positive ou  int\'egrable   sur $\O_{g,{\lambda}}$,
 $$ \int_{\O_{g,{\lambda}}} \var(l)d\mu_{\O_{g,{\lambda}}} (l)  =
\int_{G/H}\int_{\O_{h,{\lambda}}} \int_{\h^\bot}\var(x.(\tilde w +
t))d_{\h^\bot}td\mu_{\O_{h,{\lambda}}}(w)d_{G/H}\dot{x},$$ o\`u
$\tilde w$ est un \'el\'ement de $ \g^*$ dont la restriction \`a $\h$ est $w$.
\end{pro}
\end{subsection}
\end{section}
\begin{section}{Techniques de r\'ecurrence dans  la construction de $\widehat{G}$. }
\label{ddkk}
\begin{subsection}{}\label{ddp}
Soit $V$ un espace vectoriel de dimension finie sur $k$ muni d'une
forme bilin\'eaire altern\'ee $B$. Soit $H$ un groupe op\'erant
dans $V$ par des automorphismes fixant $B$. On note $V^{\bot_B} = \{v \in V\, / \, B(v , w ) = 0, \, \, \mbox{ pour tout }w \in V\}$. Alors $B$ induit sur $V/{V^{\bot_B}}$ une
structure symplectique invariante sous l'action de $H$, not\'ee
encore $B$. Si $ V \neq V^{\bot_B} $, le groupe m\'etap\-lec\-tique $Mp(V/{V^{\bot_B}})$ est bien
d\'efini. Si $ V = V^{\bot_B} $, on pose
 $ Mp(V/{V^{\bot_B}}) = \{\pm 1\} $. On note $H^V$ l'ensemble des
couples $(x , m)$, o\`u $x \in H $ et $ m \in Mp(V/{V^{\bot_B}})$ tels que
$ x$ et $m$ aient m\^eme
image dans le groupe symplectique $Sp(V/{V^{\bot_B}}) $. On a une projection naturelle de  $H^V$  sur
$H$ qui au couple $ (x , m)$ de  $H^V$ fait correspondre l'\'el\'ement $x$ de $H$.
Si $H$ est un groupe localement compact et si son action dans $V$ est continue,
 cette projection  est continue et son noyau est constitu\'e  de deux \'el\'ements qui sont
centraux dans $H^V$. Autrement dit, $H^V$ est une extension
centrale d'ordre deux de $H$, dite aussi extension m\'etaplectique associ\'ee \`a
l'action de $H$ dans $V$.
 On peut \'egalement d\'ecrire  $H^V$ comme l'ensemble
des couples  $ (x , \phi )$, o\`u $ x \in H $ et $\phi $ est une
fonction sur les lagrangiens de $V/{V^{\bot_B}}$, tels que, notant
$\bar{x}$ l'image de $x$ dans $ Sp(V/{V^{\bot_B}})$, on ait $
(\bar{x }, \phi) \in Mp(V/{V^{\bot_B}}) $. Remarquons que  si $W$ est un
sous-espace vectoriel $H$-invariant de $V$ contenu dans $V^{\bot_B}$,
le groupe $H^{V/W}$ est bien d\'efini et il est \'egal \`a $ H^V$.
\end{subsection}
\begin{subsection}{}\label{ddq}
Les notations sont celles du num\'ero pr\'ec\'edent. Soit $W$ un
sous-espace $H$-invariant de $V$ tel que son orthogonal 
par $B$  soit contenu dans $ W + V^{\bot_B} $. Les groupes
$H^V $ et $H^W$ sont d\'efinis. Si $ m$  est un sous-espace totalement 
isotrope de dimension maximale dans $W$ alors $ m + V^{\bot_B} $ est  totalement isotrope  maximal dans $V$.
Rappelons que  l'application $ \l \mapsto \l /{V^{\bot_B}}$ est une bijection, canonique,  de l'ensemble des sous-espaces totalement isotropes maximaux de $V$ sur l'ensemble
des lagrangiens
de $ V/{V^{\bot_B}} $. On identifie donc ces deux ensembles au moyen de
cette bijection.
Le r\'esultat suivant est d\^u \`a M. Duflo \cite{[Du]}.
\begin{lem}\label{yddq}
\

$i)$  Soient $ ( x , \phi)$ dans $H^V$ et $ ( x , \psi)$ dans
$H^W$. Le nombre $ \phi ( m + V^{\bot_B})\psi(m)^{-1} $ ne d\'epend pas du choix
du lagrangien $m$ de $W$. On le note $\phi \psi^{-1}$.

 $ii)$ Soient $(x , \phi),   ( x' , \phi') \in H^V$ et
$ (x , \psi),  ( x' , \psi') \in  H^W $. On pose
$ (x , \phi)( x' , \phi') = ( x x' , \phi'')$ et  $ (x , \psi) ( x' ,
\psi') =( x x' ,  \psi'')$.
On a~:   $ \phi'' {\psi''}^{-1}  =  (\phi  \psi^{-1} )
                                     ( \phi' {\psi'}^{-1})$.
\end{lem}
\end{subsection}
\begin{subsection}{}\label{m22} Dans ce paragraphe, nous allons d\'ecrire les techniques de r\'ecurrence dans  la const\-ruc\-tion de $\widehat{G}$.
 \begin{subsubsection}{}\label{ddv} Soit $g \in {\goth g}^* $ et $\beta_g$   la forme  bilin\'eaire altern\'ee  sur $\g$  d\'efinie par : $$\beta_g(X, Y) = <g, [X, Y]> \mbox{,
pour tout } X, Y \in \g.$$
  Elle est  invariante par $G(g)$, si bien que  l'extension m\'etaplectique
$G(g)^\goth g $  correspondant \`a  l'action  de $G(g)$ dans
$(\goth g, \beta_g) $ est bien d\'efinie.  L'\'el\'ement non
trivial du noyau de la projection canonique de $ G(g)^\goth g $ sur $G(g)$
 est not\'e $ (1,-1)$.
On d\'efinit un
caract\`ere $\chi_g$ de ${}^uG(g)^\goth g$, radical unipotent de $G(g)^\goth g$,  par la formule
 $$\chi_g(\exp(X))= \varsigma(<g, X>), \, X \in
{}^u\g(g).$$
On d\'esigne par $Y_G(g)$ l'ensemble des classes
d'\'equivalence des repr\'esentations unitaires irr\'eductibles de
$G(g)^\goth g $ dont la restriction \`a ${}^uG(g)^\goth g$
est multiple de $\chi_g$  et telle que $\tau (1,-1)= -{\rm Id}$.
 Soit $R_g$ un facteur r\'eductif de $G(g)$  d'alg\`ebre de Lie $\r_g$ et
$R_g^\goth g $  son  image r\'eciproque  dans $G(g)^\goth g $. Alors $R_g^\goth g
$ est un facteur r\'eductif de  $G(g)^\goth g $ et l'application $\tau \longrightarrow \tau_{|R_g^\g}$  permet d'identifier  $Y_G(g)$  \`a  l'ensemble des
classes des repr\'esentations unitaires irr\'eductibles $\tau$ de
$R_g^{\goth g} $  telles que $\tau(1 , -1) = - {\rm Id}$.

Maintenant, on se place dans les conditions  
du num\'ero \ref{meet}.
  Le groupe $H(h)$ (resp. $G(g)$)  op\`ere
dans $\h$ en fixant la forme lin\'eaire $h$. Donc, les rev\^etements  m\'etaplectique $H(h)^{\h}$
et  $G(g)^{\h}$  sont bien d\'efinis et on a~: $$H(h)^{\h}
=G(g)^{\h}U. $$  En cons\'equence  $R_g^\h$,   l'image r\'eciproque de $R_g$
dans $G(g)^{\h}$, est un facteur r\'eductif de~$H(h)^{\h}$.

 \noindent Soient $ ( x , \var)$ et $(x ,
\psi)$ deux \'el\'ements de $G(g)^\g$ et $G(g)^\h$
respectivement repr\'esentant $x$. On d\'efinit le scalaire $\var\psi^{-1}$ comme
dans le lemme  \ref{yddq} . On pose, pour $ \tau \in Y_G(g)$, 
\begin{eqnarray}
\label{vc} \tilde \tau (x , \psi) = \var\psi^{-1}\tau (x , \var) .
\end{eqnarray}
 $\tilde \tau (x , \psi)$ ne d\'epend pas du choix du couple $(x ,
 \var)$, de plus,
 la formule (\ref{vc}) d\'efinit  un \'el\'ement de $Y_H(h)$.

On note  $\goth q = \ker(u)$  et $Q = \exp(\goth q)$.  Alors, $\goth q$ est un id\'eal de $\h$ et $Q$ est un  sous-groupe unipotent normal dans $H$.  On pose  $G_1 = H/Q$,  $\g_1 = \h/ \goth q$,
et $g_1$ l'\'el\'ement de $\g_1^*$ obtenu par passage au quotient
 de $h$. Avec ces notations, on a~ :  $G_1(g_1)= H(h)/Q$ et
$\g_1(g_1)= \h(h)/\goth q$. De plus, comme $H(h)^{\g_1}=
H(h)^{\h}$, on a  $G_1(g_1)^{\g_1} = H(h)^{\h}/Q$, si bien que  $R_g^\h$ (resp.  $\tilde \tau$)  s'identifie canoniquement \`a un facteur r\'eductif
de $G_1(g_1)^{\g_1}$ (resp. \`a un \'el\'ement de $Y_{G_1}(g_1)$).

\end{subsubsection}
\begin{subsubsection}{}\label{repr}
 Soit $g \in \g^* $
  de type unipotent et $ \tau \in Y_G(g) $.
On suppose que $\g$ contient un id\'eal ab\'elien $\u$, $G$-invariant, et
contenu dans ${}^u\g$. En reprenant les notations du paragraphe \ref{ddv},
 on note   $\pi_{h,  \tilde \tau}$ (resp.  $\pi_{g_1,  \tilde \tau}$)
  la classe de rep\-r\'e\-sen\-tations unitaires
  irr\'eductibles de
 $H$ (resp. $G_1$) associ\'ee \`a la donn\'ee $(h,  \tilde \tau)$ (resp. $(g_1,  \tilde \tau)$) par la m\'ethode
  des orbites  de Kirillov-Duflo. On a :
\begin{eqnarray}\label{represe}
  \pi_{h,  \tilde \tau} &=& \pi_{g_1,  \tilde \tau} \circ p_1,
\end{eqnarray}
  o\`u $p_1$ est la projection canonique de $H$ dans $G_1$.
L'\'el\'ement  de $\widehat{G}$
associ\'e \`a la donn\'ee $(g, \tau)$ par la m\'ethode   des orbites  de Kirillov-Duflo  est alors
\begin{eqnarray}\label{repres}
\pi_{g , \tau }
= \mbox{Ind}_H^G \pi_{h,  \tilde \tau}.
\end{eqnarray}
\end{subsubsection}
\end{subsection}
\begin{subsection}{ }\label{kh22} Dor\'enavant, on suppose   que $G$ est   r\'esoluble  presque connexe et que $p$,
  la carac\-t\'eristique r\'esiduelle de $k$,  est diff\'erente de
  $2$.  Soit $0< \e < a_G$.
Soit $g \in \g^*$  de type unipotent et  $ \tau \in Y_G(g) $.
Nous allons montrer
   qu'il existe  une forme
lin\'eaire $\lambda_\tau$ sur $\g(g)$ dont la restriction \`a ${}^u\g(g)$ est \'egale \`a celle de $g$ et telle que, si $X \in \mathcal \g(g)_{ \e |n_F|_p}$, on ait $\tau(\exp X)= \varsigma(<\lambda_\tau, X>)\mbox{ Id }$ (dans cette situation, nous dirons que $\lambda_\tau$ est associ\'ee \`a~$\tau$).
\begin{subsubsection}{}
 Soit  
 $R_g$ un facteur r\'eductif de $G(g)$ d'alg\`ebre de Lie $\r_g$. On d\'esigne par $R_g^\g$
 l'image r\'eciproque de $R_g$  dans $G(g)^\g$.  Comme  $ R^\g_{g, \e}=\exp_{R_g^\g}(\r_{g, \e |n_F|_p})$ est un sous-groupe   central de $R_g^\g$ donc,  d'apr\`es le lemme de Schur,
il existe un caract\`ere  $\Psi_\tau$ de $R^\g_{g, \e}$ tel que
$\tau_{| R^\g_{g, \e}}=\Psi_\tau \mbox{Id}$. Puisque $\exp_{R_g^\g}: X \in \r_{g, \e |n_F|_p} \longmapsto
\exp_{R_g^\g}(X) \in R^\g_{g, \e} $ est un isomorphisme de groupes,
   on d\'efinit 
un caract\`ere $\psi_{\tau}$ de $\r_{g, \e |n_F|_p}$, en posant,
$$
\psi_\tau(X) = \Psi_\tau (\exp_{R_g^\g} (X)),  \mbox{   pour tout }
X \in \r_{g, \e |n_F|_p}.
$$
Il existe alors une forme
lin\'eaire  $\lambda_{\tau}$  sur $\r_g$ v\'erifiant
\begin{eqnarray}
\label{isra}
\psi_\tau(X) &= & \‡(<\lambda_{\tau}, X>) \mbox{, pour tout  } X
\in \r_{g, \e |n_F|_p}.
\end{eqnarray}
En la  prolongeant  par $g_{|{}^u\g(g)}$ sur ${}^u\g(g)$, on a : 
 \begin{eqnarray}\label{new}
\tau(\exp_{G(g)^\g}(X))& = &    \‡(<\lambda_{\tau}, X>)\mbox{ Id } \mbox{, pour tout  } X
\in \g(g)_{ \e |n_F|_p}.
\end{eqnarray}
Remarquons que $\lambda_{\tau}$ est d\'efinie modulo $(\g(g)_{ \e |n_F|_p})^\perp$.
\end{subsubsection}
\begin{subsubsection}{}
Soit $\u$  un id\'eal ab\'elien  de  $\g$, $G$-invariant, et
contenu dans ${}^u\g$. Avec les notations du paragraphe \ref{ddv},  $\r_g$ (resp. $R_g$)  est un facteur r\'eductif de $\h(h)$ (resp. $H(h)$).  De plus, $ R^\h_{g, \e}$ est un sous-groupe   central de $R_g^\h$ et 
 $\exp_{R_g^\h} : X \in \r_{g, \e |n_F|_p} \longmapsto
\exp_{R_g^\h}(X) \in R_{g, \e}^\h$ est un isomorphisme de groupes. 
 Pour $X \in \r_{g, \e |n_F|_p}$, on pose   $\exp_{R_g^\h}(X)= (x_X , \psi_X)$ et
$\exp_{R_g^\g}(X)= (x_X', \var_X)$. Il est \'evident  
que $x_X =x_X'$. Il r\'esulte du lemme \ref{yddq} que
l'application $X \in \r_{g, \e |n_F|_p} \longmapsto \var_X \psi_X^{-1}$ est un
morphisme de groupes  de $\r_{g, \e |n_F|_p}$ dans le groupe des racines huiti\`eme de l'unit\'e.
 Mais $8\r_{g, \e |n_F|_p}= \r_{g, \e |n_F|_p}$ donc $\var_X \psi_X^{-1} = 1$, pour tout
 $X \in \r_{g, \e |n_F|_p}$. Si bien que l'on a~: $$ \tilde
\tau(\exp_{R_g^\h}(X)) =  \‡(<\lambda_{\tau}, X>)\mbox{Id},
\mbox{ pour tout } X \in \r_{g, \e |n_F|_p} .$$
Ainsi, si on prolonge $\lambda_{\tau}$ en une forme lin\'eaire sur $\h(h)$   par $g_{|\u}$ sur $\u$ alors $\lambda_{\tau}$ est associ\'ee \`a  $ \tilde
\tau$.

\end{subsubsection}
\end{subsection}
\end{section}

\begin{section}{Formule du caract\`ere au voisinage des \'el\'ements
semi-simples.}\label{4}
\begin{subsection}{D\'efinition  de la fonction  $\phi_{g,\tau, s}$.}\label{hcgg} Les donn\'ees et les notations
sont celles des paragraphes \ref{kh22}, \ref{kh24}, et \ref{meet}.  Soit $s \in G$
 semi-simple,  $0< \e < a_G$, et $\lambda_\tau$ une forme lin\'eaire sur $\r_g$ v\'erifiant l'\'egalit\'e (\ref{isra}), consid\'er\'ee comme
un \'el\'ement de $L(g)$, en la prolongeant par $g_{|{}^u\g(g)}$ sur ${}^u\g(g)$.
 On note $\mathcal O_{g,\lambda_\tau, s} =\mathcal O_{g, \lambda_\tau} \cap
 \g^* (s) $, l'ensemble des points fixes de $s$ dans $\mathcal O_{g, \lambda_\tau} $.  On va
d\'efinir une fonction $ \phi_{g,\tau, s}$ sur $\mathcal O_{g,
\lambda_\tau, s}$ comme suit~:  soit $l \in \O_{g,\lambda_\tau,
s}$ et
 $x \in G$ tel que $l= x.f= f-g + x.g$. On a $G(l) = G(x.g)$ et
 $\beta_l = \beta_{x.g}$.  L'application $\alpha : \g/\g(g)\longrightarrow \g/\g(x.g), \, X \longmapsto x.X$, est un isomorphisme d'espaces symplectiques.  Soit $\alpha_x$ l'isomorphisme de   $ Sp(\g/\g(g))$ sur $ Sp(\g/\g(x.g))$ d\'efini par $ \alpha_x(y) = \alpha y \alpha^{-1} $.  Alors, $\alpha_x$ se rel\`eve de mani\`ere unique en un isomorphisme de  $ Mp(\g/\g(g))$ sur $ Mp(\g/\g(x.g))$, not\'e encore $\alpha_x$. Maintenant,  
 l'application $\tilde \alpha_x :   G(g)^\g \longrightarrow
 G(x.g)^\g$, \, $(y,m) \longmapsto (xyx^{-1},\alpha_x(m)) $ est un isomorphisme de
 groupes.  On pose ${}^x\tau = \tau \circ ({\tilde \alpha_x})^{-1}$.
 La valeur de la fonction
$\phi_{g,\tau,s}$  en $l$ est donn\'ee par~:
\begin{eqnarray}\label{functmetap}
\phi_{g,\tau, s}(l) &=& {\rm Tr} ({}^x\tau  (\tilde{s}))\Phi_1(\rho_x(\tilde{s}))|\det(1-s^{-1})_{(1-s)\g/\g(l)}|_p^{-\frac{1}{2} },
\end{eqnarray}
o\`u $\tilde{s}$  est un \'el\'ement  du rev\^etement
m\'etaplectique $G(x.g)^\g$ repr\'esentant $s$ et $\rho_x$ est
l'application canonique de $G(x.g)^\g$ dans $Mp(\g/\g(x.g))$. On
v\'erifie que le membre de droite de (\ref{functmetap}) ne d\'epend
ni du choix de $x \in G$ tel que $l=x.f$, ni du choix de
l'\'el\'ement $\tilde s$ de $G(x.g)^\g$ situ\'e au-dessus de
$s$.  La fonction $\phi_{g,\tau, s}$ est $G(s)$-invariante.
 Remarquons que si $s$ est central alors la fonction
$\phi_{g,\tau, s}$ est constante sur l'orbite  $\mathcal
O_{g,\lambda_\tau, s}=  \mathcal O_{g,\lambda_\tau}$ et elle est
\'egale \`a  ${\rm Tr}(\tau (s, \psi)) \psi(\l)$, o\`u $(s, \psi) \in G(g)^\g$ un relev\'e de $s$ et $\l$ est un
lagrangien de $\g/\g(g)$. Remarquons aussi que si $x \in G$ alors  $x.\lambda_{\tau}$ est associ\'ee \`a ${{}^x\tau}$ et on a : 
\begin{eqnarray} \label{kl}
 \phi_{x.g,{}^x\tau, s} &=& \phi_{g,\tau, s}.
\end{eqnarray}
\end{subsection}
\begin{subsection}{} \label{mlhgf} Dans cette section,  on suppose   que $\g$ contient un id\'eal ab\'elien $\u$,
$G$-invariant, et contenu dans ${}^u\g$. On reprend les donn\'ees et les notations
des  num\'eros  \ref{ddv} et \ref{hcgg}. On pose
$$G_{s,H } = \{x \in G / xsx^{-1} \in H\}.$$
 On va d\'emontrer le r\'esultat suivant :
\begin{pro} \label{ymlhgf}
 On suppose  $s \in H$. Soient $x \in G_{s,H }$ et
$l \in \O_{g, \lambda_\tau, xsx^{-1}}$ tels que  $r(l) \in \O_{h,
\lambda_\tau }$. Alors, on a~ :
\begin{eqnarray}\label{funor}
 \phi_{g,\tau, s}(x^{-1}.l) &=& \phi_{h, \tilde \tau,
xsx^{-1}}(r(l))   \Delta_{H,
G}^{-\frac{1}{2}}(xsx^{-1})|\det(1-(xsx^{-1})^{-1})_{(1-xsx^{-1})(\g/\h)}|_p^{-1}, \quad \quad
\end{eqnarray}
$ r$ \'etant l'application restriction de $\g^*$ \`a $\h^*$.
\end{pro}
\begin{dem}
Par transport de structure, on a~: $$\phi_{g,\tau, s}(x^{-1}.l) =
\phi_{x.g, {}^x\tau, xsx^{-1}}(l)=\phi_{g, \tau, xsx^{-1}}(l) ,$$
o\`u la deuxi\`eme \'egalit\'e r\'esulte de la formule (\ref{kl}).
 Comme par hypoth\`ese $r(l) \in
\O_{h, \lambda_\tau }$ et comme $U.l = l + \h^\perp$,
  il existe $y \in H$ tel que
 $l= y.f$.
 Posons $s' = xsx^{-1} \in H$ et soit  $V$ un suppl\'ementaire $s'$-invariant de
 $\h(r(l))/\g(l)$ dans $\h/\g(l)$. Alors $V$ est un sous-espace symplectique de
 $(\g/\g(l), \beta_l)$ et l'application
 $\alpha_V : V \longrightarrow \h/ \h(r(l))$, $X \mod \g(l) \longmapsto X \mod \h(r(l))$,
  est un isomorphisme symplectique $s'$-equivariant.
 D\'ecomposons $
(1-s').V$   en somme directe ortho\-go\-nale par rapport \`a sa
structure symplectique  en  $W_1 \oplus W_2$ o\`u $W_1$ (resp. $W_2$) est un sous-espace symplectique $s'$-invariant 
  poss\'edant 
un lagrangien $\l_1$ (resp. $\l_2$) tel que $ \l_1 = s'. \l_1$ (resp.  $(s'.\l_2)\cap \l_2 = 0$).
 D'autre part, \'ecrivons $\g/\g(l) = V \oplus V^\perp $. Alors $V^\perp$ est un
  sous-espace symplectique
 de $\g/\g(l)$, $s'$-invariant et  on v\'erifie que $\h(r(l))/\g(l)$ est un
 lagrangien, bien entendu, $s'$-invariant
  de $V^\perp$.
La formule (\ref{vc}) donne
$$
\psi(\alpha_V(\l_0+ \l_1 + \l_2)){}^y\tilde \tau (s' , \psi) = \var(\l_0 +(\h(r(l))/\g(l)) +\l_1 +
\l_2 ){}^y\tau (s' , \var),
$$
o\`u $\l_0$ est un lagrangien de $V(s')$. 
Il s'ensuit que l'on a :
$$
 \phi_{g, \tau, s'}(l)|\det(1-s'^{-1})_{(1-s')\g/\g(l)}|_p^{\frac{1}{2} } = \phi_{h, \tilde \tau, s'}(r(l))
|\det(1-s'^{-1})_{(1-s')\h/\h(r(l))}|_p^{\frac{1}{2} } .
$$
D'autre part, les espaces $\g/\h$ et $\h(r(l))/ \g(l)$ sont en
dualit\'e, $s'$-invariante,  par $\beta_l$. Donc, $|\det
s'_{\g/\h}|_p = |\det s'_{\h(r(l))/ \g(l)}|_p^{-1}$ et $|\det
(1-s')_{(1-s').(\g/\h)}|_p = |\det (1-s'^{-1})_{(1-s').(\h(r(l))/
\g(l)})|_p$. Il en r\'esulte que
$$
|\det(1-s'^{-1})_{(1-s')\g/\g(l)}|_p^{\frac{1}{2}}
=|\det(1-s'^{-1})_{(1-s')\h/\h(r(l))}|_p^{\frac{1}{2}}
\Delta_{H, G}^{\frac{1}{2}}(s')|\det(1-s'^{-1})_{(1-s')\g/\h}|_p.
$$
Si bien que l'on a
$
\phi_{g, \tau,
s'}(l)  =  \phi_{h, \tilde \tau, s'}(r(l)) \Delta_{H,
G}^{-\frac{1}{2}}(s')|\det(1-s'^{-1})_{(1-s')\g/\h}|_p^{-1}.
$
\end{dem}
\end{subsection}

\begin{subsection}{Choix de voisinages semi-simples.}\label{vbv}
 Pour chaque $n \in \N$ tel que  $2n \leq \dim \g$,  on fixe :
\begin{itemize}
\item un espace symplectique $(V_n, B_n)$ de dimension \'egale
$2n$. 
\item   des repr\'esentants  $(T_i)_{i \in I_n}$ de l'ensemble des
classes de conjugaison  des tores maximaux de $Sp(V_n)$.
\item$(e_1, \ldots, e_n, f_1, \ldots , f_n )$  une base
symplectique de  $(V_n, B_n)$. On pose~ : $$ r_{ n, +} =\O e_1 +
  \cdots \O e_n + \O f_1 + \cdots \O f_n , \, \, \,   r_{ n,-} =\O e_1+
   \cdots \O e_n + \O \varpi^{-1} f_1 + \cdots \O
 \varpi^{-1} f_n ,$$
 $$ Sp(V_n)(r_{n, +, -}) = \{ x \in Sp(V_n) \, , \, x.r_{ n,+} = r_{
n,+}\, \mbox{ et }\, \, \,  x.r_{n,-} = r_{ n,-} \}. $$
On note
   $$ a_n = \sup\{0<\e<a \, / \,  \forall \,  i \in I_n , \, \, \exists \,  x \in Sp(V_n),\, \,   T_{i, \e} \subset  xSp(V_n)(r_{n,+,-})x^{-1}\}. $$ 
   \end{itemize}
   On
suppose que $a$ v\'erifie les conditions de la proposition
\ref{yresoluble}. On pose
\begin{eqnarray} \label{cond}
{\bf a}_G = \inf_{ n \leq \frac{1}{2}\dim \g}\{a_n, \, a_G
\}.
\end{eqnarray}
 Soit $s \in G$, semi-simple. D'apr\`es la proposition \ref{yexc}, il existe  $0<\e(s)<{\bf a}_G$ 
tel que, pour tout  $(V, B)$
sous-espace symplectique, $s$-invariant,  de ${}^u\g$ et
 pour lequel $s_{|V} \in Sp(V, B)$,
 en notant $V_s =(1-s).V $ et
$\tilde{ s_{|V_s}}$ un rel\`evement de $s_{|V_s}$ dans $Mp(V_s)$,
 on
ait la formule~(\ref{ramadann}).
On note  $a''_G(s) = \sup\{0<\e(s)<{\bf a}_G\}$. 
 Remarquons que si $H $ est un
sous-groupe presque alg\'ebrique de $G$ contenant $s$ alors
$a''_H(s) \geq a''_G(s)$. S'il n'y a pas de confusion \`a
craindre, on note $a''(s)$ le nombre r\'eel $a''_G(s)$.
\end{subsection}

\begin{subsection}{}\label{hd}   Fixons $0 < \e <{\bf a}_G$. Soit  $s \in G$ semi-simple, $G(s)$
le centralisateur de $s$ dans $G$, et  $0<\e(s)< \inf\{a'(s), a''(s), \e\} $  tel que l'application
 $\bar \Psi :  G \times_{G(s)} G(s)_{\e(s)}\longrightarrow \W(s , \e(s))$, \,
  $[x, y]\longmapsto  xs y x^{-1}$ soit un diff\'eo\-mor\-phisme.
Soit $g \in \g^*$ de type unipotent et $\tau \in Y_G(g)$. 
 Si $\pi_{g, \tau}$ est admissible alors la restriction de $\Theta_{g, \tau}$ \`a $\W(s , \e(s) )$    est
 une fonction g\'en\'eralis\'ee  $G$-invariante, 
 on d\'esigne par $\theta_{g, \tau,s}$ la fonction g\'en\'eralis\'ee sur
 $G(s)_{\e(s)}$ qui lui est associ\'ee  par le th\'eor\`eme \ref{yPTo}.

 Choisissons une forme lin\'eaire $\lambda_\tau$  sur $\g(g)$ v\'erifiant l'\'egalit\'e (\ref{new}). 
Alors $\O_{g,\lambda_\tau, s}$
 est une sous-vari\'et\'e localement ferm\'ee, r\'eunion finie de $G(s)$-orbites co-adjointes. Si elle n'est pas vide, elle supporte une mesure positive canonique $d\mu_{\O_{g,\lambda_\tau, s}}$ dont la restriction \`a chaque $G(s)$-orbite qu'elle contient est la mesure de Liouville de celle-ci.  Lorsque  $\O_{g,\lambda_\tau, s}$  est vide, nous conviendrons  que $d\mu_{\O_{g,\lambda_\tau, s}}$ est la mesure nulle.

Soit
 $d_{G(s)}x$ et $d_{\g(s)}X$ deux mesures de Haar sur $G(s)$ et $\g(s)$ respectivement
qui se correspondent.
\begin{tho} \label{yhd} 
On suppose que $\O_{g, \lambda_\tau}$ est ferm\'ee dans $\g^*$.
Alors, on a~:
\begin{eqnarray}
\label{ Torasso} \theta_{g,\tau,s}(\beta d_{G(s)}x)&=& \int_{
\O_{g,\lambda_\tau, s}} (\beta \circ \exp
d_{\g(s)}X)\,\widehat{}_{\g(s)} \, (l)\phi_{g,\tau,s}(l)
d\mu_{\O_{g,\lambda_\tau, s}}(l),
\end{eqnarray}
 pour tout $\beta \in C_c^\infty(G(s)_{\e(s)})$.
\end{tho}
Comme  application  de ce th\'eor\`eme, nous donnons dans la deuxi\`eme partie du pr\'esent travail, sous l'hypoth\`ese suppl\'ementaire que $G$ est unimodulaire,  une
 d\'emonstration de  la mesure de Plancherel pour $G$.
 
 \noindent{\it {\bf  Remarque : } $\O_{g, \lambda_\tau}$ est ferm\'ee  dans $\g^*$ si et seulement si $\O_{g}$ est ferm\'ee dans $\g^*$.}
\end{subsection}

\begin{subsection}{}\label{m}
En traitant \`a part  le cas  ${}^uG =\{ 1\}$,  la d\'emonstration  du
th\'eor\`eme \ref{yhd}  se fait  par r\'ecurrence sur la dimension de
$G$. Compte tenu du lemme suivant, il suffit d'une part de
d\'emontrer le th\'eor\`eme \ref{yhd} dans le cas du num\'ero
\ref{ma}, c'est le cas du produit semi-direct d'un groupe de
Heisenberg $H$ par un groupe r\'eductif agissant trivialement sur
le centre de $H$, et d'autre part de d\'emontrer dans le cas du
num\'ero \ref{mlk}, que si le th\'eor\`eme \ref{yhd} est vrai pour
$H$ alors il est vrai pour $G$. Soit $l$ une forme lin\'eaire sur $\g$.
\begin{lem}\label{ym}
 L'alternative suivante r\'eunit tous les cas qui peuvent se pr\'esenter.
\begin{enumerate}
\item Le radical unipotent  de $\g$ est ou bien  nul, ou bien une alg\`ebre
de  Heisenberg de di\-men\-sion $2n+1$ ( $n \geq 0$ ) avec centre, not\'e $\goth
z$, central dans $\g$ et  la restriction de
$l$ \`a $\goth z$ est non nulle.
\item L'alg\`ebre de Lie  $\g$ contient un id\'eal ab\'elien $\u$,
 $G$-invariant et inclus dans le radical unipotent de $\g$  tel que,
 si on pose $u = l_{|\u}$, on ait
$$\dim (\g(u)/  \ker  u) < \dim \g.$$

\end{enumerate}
\end{lem}
\end{subsection}
\begin{subsection}{}\label{rtl}
Dans cette section, nous allons d\'emontrer le th\'eor\`eme
\ref{yhd} dans le cas o\`u
 ${}^uG = \{1\}$.
 Soit $0<\e< a_G $.
  Soient $\tau$ une
repr\'esentation unitaire irr\'eductible de $ G$ et
$\lambda_\tau \in \g^*$ tel que l'on ait
$$\tau(\exp_{G}(X))=
\varsigma(<\lambda_\tau, X>)\mbox{Id}, \, \mbox{ pour tout }  X \in
\g_{\e|n_F|_p}. $$
On se donne \'egalement  $s \in G$.  On a :   
  $$\mbox{ Tr} \, \tau(s\exp_{G}(X)) =\varsigma(<\lambda_\tau, X>)\mbox{Tr}\, \tau(s),
   \, \mbox{ pour tout }  X \in \g_{\e|n_F|_p} . $$
  D'autre part, $0$ est l'unique forme lin\'eaire de type unipotent et  l'orbite co-adjointe $\O_{0, \lambda_\tau}$  est r\'eduite \`a $\{\lambda_\tau\}$ qui est fix\'e par $s $.
   De plus, $\phi_{0,\tau,s}(\lambda_\tau) = \mbox{Tr}\, \tau(s)$ et la mesure de Liouville $d\mu_{0,\lambda_\tau,s}$  sur $\O_{0, \lambda_\tau, s}$ est la mesure de Dirac en $\lambda_\tau$.
 D'o\`u le th\'eor\`eme dans ce cas.

\end{subsection}
\end{section}
\begin{section}{}\label{5} Dans ce paragraphe, nous allons d\'emontrer le th\'eor\`eme
\ref{yhd} au voisinage de l'\'el\'ement neutre de $G$.
\begin{subsection}{} \label{ma}
Dans cette section,  on suppose que  ${}^uG$ est un groupe de Heisenberg  dont
 le centre est central dans $G$.  
 On se donne
un facteur r\'eductif $T$ de $G$. Alors, le centre $\goth z$ de ${}^u\g$  est fix\'e 
par $T$. Par suite, $\goth z$ admet  un suppl\'ementaire  $V$ dans ${}^u\g$,
 stable par $T$.
 On fixe  un vecteur non nul $E$ de $\goth z$.
On d\'efinit une forme symplectique $B$ sur $V$ par la formule
$$[v, w] = B(v, w)E, \, \, v, w \in V.$$  Le groupe ${}^uG$
s'identifie canoniquement au groupe de Heisenberg $H$ associ\'e
\`a $(V, B)$.

\noindent Soit  $g \in \g^*$ de type unipotent tel que $g_{|\goth z} \neq 0$.
On note $a_0= g(E)$.  On d\'efinit une forme lin\'eaire
$E_{a_0}^*$ sur ${}^u\g = V + kE$ par la formule suivante
$$E_{a_0}^*(E) = a_0, \, E_{a_0}^*( V ) =  0.$$ En la prolongeant
par $0$ sur  $\t$,  l'alg\`ebre de Lie de $T$, on obtient une
forme lin\'eaire sur $\g$ not\'ee aussi $E_{a_0}^*$. Elle est
 de type unipotent  et on a  $G.g =
G.E_{a_0}^*$.  Si bien que
 l'on peut supposer  que $g=E_{a_0}^*$.
 
D\'esignons par $ T^V$ l'extension m\'etaplectique
 de $T$ associ\'e \`a son action dans $(V, B)$~: 
 c'est l'ensemble des couples $(x, \psi)$,
 o\`u $x\in T$ et $\psi$ est
 une fonction sur les lagrangiens de $V$, tels que, notant $\bar{x}$ l'image de $x$
 dans $Sp(V)$, on ait $(\bar{x}, \psi) \in Mp(V)$.
Le morphisme de groupes de $ T^V$ dans $Mp(V)$ qui \`a $(x,
\psi)$ associe $(\bar{x}, \psi)$ est not\'e $\eta_{T}^V$.
 On a une projection naturelle de $T^V$ sur $T$, qui au couple $(x, \psi)$ de $T^V$ fait
  correspondre $x$ de $T$. Le noyau de cette projection
   est constitu\'e de deux \'el\'ements  qui sont centraux dans $T^V$ dont on note 
$(1, -1)$ l'\'el\'ement non trivial.

 Soit  $\tau \in  ( T^{V}) \, \widehat{}$ \,  telle que $\tau(1, -1) =
-\mbox{Id}$
 et  soit  $\pi_{E_{a_0}^*, \tau}$ 
  la classe
de repr\'e\-sen\-tations unitaires irr\'eductibles de $G$ 
 associ\'ee \`a  $(E_{a_0}^*, \tau)$ : 
notons $\pi_{E_{a_0}^*}$ la classe de repr\'esentations unitaires
irr\'eductibles de ${}^uG$ associ\'ee \`a $E_{a_0}^*$ par la
m\'ethode des orbites de Kirillov et  $S_{E_{a_0}^*}$ la
repr\'esentation m\'etaplectique associ\'ee au caract\`ere
$\varsigma_{a_0}$, r\'ealis\'ee dans l'espace 
de $\pi_{E_{a_0}^*}$. On a~: $$\pi_{E_{a_0}^*,  \tau} = \tau
\otimes S_{E_{a_0}^*}\pi_{E_{a_0}^*}, $$ o\`u  $\tau \otimes
S_{E_{a_0}^*}\pi_{E_{a_0}^*}(x.y) = \tau(\hat{x}) \otimes
S_{E_{a_0}^*}(\eta_{T}^V(\hat{x}))\pi_{E_{a_0}^*}(y), \, \, x \in T, \, \, y
\in {}^uG, $ \, $\hat{x}$ \'etant un rel\`evement de $x$ dans
$T^V$.
\begin{subsubsection}{}\label{ala}
 Consid\'erons la d\'ecomposition, $T$-invariante,
  $V= V_1 \oplus V_2$, o\`u $V_1$ est le sous-espace symplectique
$$V_1 =\{v \in V \, / \, adX.v = 0 \, , \, \mbox{ pour tout } X \in \t \} $$
et $V_2$ son orthogonal  par rapport \`a $B$. Notons
 $$ \h_i = V_i + kE, \, H_i = \exp(\h_i),
  \, E_i^* = E_{a_0|\h_i}^*, \, i = 1,2.   $$
Soient $\pi_{E_i^*}$ la classe de repr\'esentations unitaires
irr\'eductibles de $H_i$ associ\'ee \`a $E_i^*$ par la m\'ethode
des orbites de Kirillov et 
 $S_{E_i^*}$ la  repr\'esentation m\'etaplectique correspondante
 de $Mp(V_i)$.

\noindent L'application $(y_1, y_2) \longmapsto y_1 y_2$ de $H_1 \times H_2$
dans $H = {}^uG$ est un homomorphisme surjectif de groupes de Lie
de noyau \'egale \`a $\Delta= \{(z, z^{-1}), \, z \in Z\}$, avec $Z = \exp(\goth z)$. La
repr\'esentation  $\pi_{E_1^*} \otimes \pi_{E_2^*}$  de $H_1
\times H_2$ est unitaire irr\'eductible dont la restriction  \`a
$\Delta$ est triviale. Elle passe au quotient en une
repr\'esentation unitaire irr\'eductible de $H$ \'equivalente \`a
$\pi_{E_{a_0}^*}$. Si $x$ est un \'el\'ement de $ Sp(V)$ qui
laisse invariant $V_i$ et $\phi$ (resp. $\phi_i$) une fonction
sur les lagrangiens de $V$ (resp. $V_i$) telle que $(x, \phi) \in
Mp(V)$ (resp. $(x_{|V_i}, \phi_i) \in Mp(V_i)$), on note $
\phi \phi_1^{-1}\phi_2^{-1}$ le nombre  $\phi(\l_1 + \l_2)
\phi_1(\l_1)^{-1}\phi_2(\l_2)^{-1}$, o\`u $\l_i$ est un lagrangien de $V_i$, $i= 1, 2$.
On a~:
$$S_{E_{a_0}^*}(x, \phi) = \phi
\phi_1^{-1}\phi_2^{-1}S_{E_1^*}(x_{|V_1}, \phi_1) \otimes
S_{E_2^*}(x_{|V_2}, \phi_2). $$
 D\'esignons par ${T}^{V_i}$ le rev\^etement m\'etaplectique au
dessus de $T$ par rapport \`a $(V_i, B)$ et par $\eta_{T}^{V_i}$
l'application canonique de ${T}^{V_i}$ dans $Mp(V_i)$.
 On d\'efinit une repr\'esentation (unitaire irr\'eductible)
$\tilde \tau$ de $({T}^{V_1})^{V_2}$ par la formule $$\tilde
\tau(x, \phi_1, \phi_2)=\phi \phi_1^{-1}\phi_2^{-1}\tau(x,
\phi).$$
    Fixons
$0 <\e<{\bf a}_G$.
 On a,  $$S_{E_1^*}(\eta_{T}^{V_1}(\tilde x))= \mbox{ Id},
\, \mbox{ pour tout  } \tilde x \in T^{V_1}_\e.$$
Rappelons que
$T_\e$ (resp. $T_\e^V $)
   est
 un sous-groupe compact ouvert, central,  de $T$ (resp. $T^V$).
 D'apr\`es le lemme de Schur, il existe un caract\`ere
$\Psi_\tau$ de $T_\e^V $ v\'erifiant
 $$\tau_{|T_\e^V} =  \Psi_\tau \mbox{Id}_{\mathcal H_\tau},$$
$\mathcal H_\tau$ \'etant l'espace de $\tau$.
On d\'efinit, par la suite,  un caract\`ere $\psi_\tau$ de $\t_{\e|n_F|_p} $, en
posant, $$ \psi_\tau(X)=  \Psi_\tau  \circ \exp_{T^V}(X)  \, ,\, X
\in \t_{\e|n_F|_p} .$$ Soit $\lambda_\tau \in \t^*$ telle que
$$\psi_\tau(X) = \varsigma(<\lambda_\tau , X>) \, ,\, \mbox{ pour tout } X
\in \t_{\e|n_F|_p}.  $$ 
 Pour $X \in \t_{\e|n_F|_p}$, on pose $$\exp_{{T}^V}(X)=
(x_X, \var_X), \, \, \exp_{{T}^{V_1}}(X)=(x_X, \var_{1, X}), \, \,
\exp_{{T}^{V_2}}(X)=(x_X, \var_{2, X}).$$
 Comme $p \neq 2$, on a $\var_X \var_{1, X}^{-1} \var_{2, X}^{-1}=
 1$. Il s'en suit que l'on a :
 $$\tilde \tau(\exp_{({T}^{V_1})^{V_2}}(X)) = \varsigma(<\lambda_\tau , X>) \mbox{Id}_{\mathcal H_\tau} \, ,\, 
\mbox{ pour tout } \,  X 
\in \t_{\e|n_F|_p}.$$
 On consid\`ere alors  la forme lin\'eaire $f$ sur $\g= \t + \h$ d\'efinie par~:
 $$f_{|\h} =E_{a_0}^* \, , \quad f_{|\t} = \lambda_\tau.$$
  Si bien que l'on a : 
$$\O_{E_{a_0}^*,\lambda_\tau }= G.f.$$  
 On munit  $V_i$
 de la mesure de Haar $d_{V_i}v$
autoduale, relativement \`a $\varsigma$, et $V$ de la mesure de
Haar produit  $d_V u =d_{V_1}vd_{V_2}w$.
Alors $d_V u$ est  la mesure de
Haar autoduale, relativement \`a $\varsigma$,  sur $V$.
On fixe une mesure de
   Haar $d_Tx$ sur $T$. On munit $G$ de la mesure de Haar
   $d_Gx = d_Tx d_Vv d\mu(t)$. La mesure de Haar sur $\g$ tangente \`a $d_Gx$ est donn\'ee par $d_\g X = d_\t X d_Vv d\mu(t)$, o\`u  $d_\t X$ est la mesure de Haar sur $\t$
   tangente \`a $d_Tx$.
   Il s'agit de d\'emontrer la formule suivante : pour tout $\alpha \in C_c^\infty(G_\e)$,
   \begin{eqnarray}\label{Adem}
 \Theta_{E_{a_0}^*, \tau}(\alpha d_{G}x)
 & =&  \dim \tau\int_{\O_{E_{a_0}^*,\lambda_\tau}} (\alpha\circ
\exp_{G} d_\g X)\, \widehat{}_\g(l)
 d\mu_{\O_{E_{a_0}^*,\lambda_\tau}}(l).
 \end{eqnarray}
Comme l'application  $V_1 \times  T_\e H_2 \longrightarrow G_\e$, $(v, x)\longmapsto x\exp(v)$ est  un diff\'eomorphisme, il suffit de montrer la formule  (\ref{Adem}) dans le cas o\`u 
 \begin{center}$\alpha = \alpha_1 \otimes \var, \quad \alpha_1 \in
C_c^\infty(V_1), \, \, \var \in
C_c^\infty(T_\e H_2).$\end{center}
\end{subsubsection}
\begin{subsubsection}{} On d\'esigne par 
$$G_2 = TH_2, \, \, \g_2 = \t + \h_2, \, \,  d_{G_2}x = d_Txd_{V_2}v d\mu(t),$$
et $d_{\g_2}X$ la mesure de Haar sur $\g_2$ tangente \`a $d_{G_2}x$.
On a : 
$$\widehat{(\alpha\circ
\exp_{G} d_\g X)}_\g  = \widehat{(\alpha_1 d_{V_1} v)}_{V_1} \otimes  \widehat{(\var\circ
\exp_{G_2} d_{\g_2}X)}_{\g_2} . $$
Donc 
$$
\int_{\O_{E_{a_0}^*,\lambda_\tau }}   (\alpha\circ
\exp_{G} d_\g X)\,  \widehat{}_\g(l)
 d\mu_{\O_{E_{a_0}^*,\lambda_\tau}}(l)
 = |a_0|_p^{- \frac{\dim
V_1}{2}}\alpha_1(0) I, $$
o\`u 
$$
I =  |a_0|_p^{ \frac{\dim
V_2}{2}} \int_{V_2}  (\var\circ
\exp_{G_2} d_{\g_2}X)\, \widehat{}_{\g_2} (\exp v.f_{|\g_2})d_{V_2}v.
$$
De plus, on a, pour tout $v \in V_2$,  
$$\exp v.f_{|\g_2} = f_{|\g_2} - a_0\tilde{B}_2(v)+ \frac{1}{2} a_0 l_v ,$$
 o\`u $\tilde{B}_2(v)$ (resp. $l_v$) est la forme lin\'eaire
sur $\g_2$ d\'efinie par~: $\tilde{B}_2(v)(X + w + tE)= B(v,w)$,\, $ X \in \t, \,  w \in V_2, \,  t \in k$ (resp.
$l_v(X + w + tE)= B(adX.v, v)$).
\end{subsubsection}
\begin{subsubsection}{}
 On a~:
\begin{eqnarray*}
   \pi_{E_{a_0}^*, \tau}(\alpha_1 \otimes \var d_Gx)&=&
    \int_{\t_{\e|n_F|_p}} \tilde \tau(\exp_{({T}^{V_1})^{V_2}}(X))
   \otimes \int_{V_1} \alpha_1 ( \exp(v))\pi_{E_1^*}(\exp(v))
   d_{V_1}v \\&&
   \otimes \int_{H_2}\var(\exp_T(X)y)
   S_{E_2^*}(\eta_{T}^{V_2}(\exp_{{T}^{V_2}}(X))) \pi_{E_2^*}(y)d_{H_2}y d_\t X.
   \end{eqnarray*}
   De plus, on a : 
   $$
   \mbox{ Tr}\left(\int_{V_1} \alpha_{1} ( \exp(v))\pi_{E_1^*}(\exp(v))d_{V_1}v\right) = |a_0|_p^{-\frac{\dim V_1}{2}}\alpha_{1} ( 0).
   $$
   Il reste \`a calculer la trace de l'op\'erateur   $\mathcal I_{E_{2}^*, \tilde \tau}(\var d_{G_2}x)$ donn\'e par :
   \begin{eqnarray*}
   \int_{\t_{\e|n_F|_p}} \tilde \tau(\exp_{({T}^{V_1})^{V_2}}(X))
   \otimes \int_{H_2}\var(\exp_T(X)y)
   S_{E_2^*}(\eta_{T}^{V_2}(\exp_{{T}^{V_2}}(X))) \pi_{E_2^*}(y)d_{H_2}y d_\t X.
    \end{eqnarray*}
   Pour cela,  nous allons
introduire le mod\`ele latticiel pour la r\'ealisation de $ \pi_{E_2^*}$.
\end{subsubsection}
\begin{subsubsection}{} \label{lattice}  Soit
$r$ un r\'eseau de $V_2$, autodual relativement \`a
$\varsigma_{a_0}$, c'est-\`a-dire que $$r = r^\perp : = \{v \in
V_2 \mbox{ tel que } \varsigma_{a_0}( B(v, r))=1\} = \{v \in V_2
\mbox{ tel que }  B(v, r) \in a_0^{-1}\O \}.$$ On  choisit $r$ de
la mani\`ere suivante :  on note $2m = \dim V_2$. On fixe une base symplectique $(e_1,
\ldots, e_m, f_1, \ldots , f_m)$   de $V_2$.
 On pose $$r_+ =\O e_1 +
  \cdots +\O e_m + \O f_1 + \cdots +\O f_m \mbox{ et } r_- = \O e_1 +
 \cdots +\O e_m + \O \varpi^{-1} f_1 + \cdots +\O
 \varpi^{-1} f_m .$$
Alors, 
$$r=  \left\{
\begin{array}{ll}
\varpi^{-\frac{ v(a_0)}{2}} r_+ \, ,& \mbox{ si } v(a_0)  \mbox{ est impair }\\
\varpi^{-\frac{v(a_0)-1}{2}} r_- \, , &\mbox{ si } v(a_0)  \mbox{ est pair }
\end{array}
\right..$$
 On d\'esigne par  $ R =\exp(r + kE) $. C'est un sous-groupe ferm\'e de $H_2$ dont
l'image dans $H_2/\exp(a_0^{-1}\O E)$ est un sous-groupe ab\'elien
maximal. On d\'efinit un  caract\`ere  $\chi_{ R}$  de $ R$ par la
formule~:
$$
 \chi_{ R}(\exp(\gamma + tE)) = \varsigma_{a_0}(t) \mbox{, pour tout }
 \gamma \in r,   t \in k.
$$
D'apr\`es \cite{[Moe-Vi-Wa]}, la repr\'esentation induite ${\rm
Ind}_{R}^{H_2} \chi_{R} $ est unitaire irr\'eductible, on la note
$\pi_{{R}}$. 
On d\'esigne par $d_r\gamma$ la mesure de Haar sur
$r$, restriction de la mesure de Haar $d_{V_2}v$ sur ${V_2}$, et
par $d_{{V_2}/r}\dot v$ la mesure de Haar quotient sur ${V_2}/r$.
Alors $\pi_{{R}}$ agit dans l'espace $\mathcal H_r$ des fonctions
$\var$  de  $ V_2$ \`a valeurs dans $\cit$ 
v\'erifiant~: $$ \var(v+m) = \varsigma_{a_0}(\frac{1}{2}B(v, m))
\var(v) \mbox{, pour }
 v \in V_2 \mbox{, } m \in  r  \mbox{, et }
 \int_{V_2/ r}|\var(v)|_p^{2}d_{V_2/ r}\dot{v}< \infty, $$
par la formule
\begin{eqnarray}
\label {ec1} (\pi_{{R}}(w,t)\var)(v)
=\varsigma_{a_0}(t-\frac{1}{2}B(v, w)) \var(v-w)\mbox{ , pour tout
}
 v, w  \in V_2, t \in k .
\end{eqnarray}
Comme le caract\`ere central de  $\pi_{{R}}$ est $\varsigma_{a_0}$ donc,
d'apr\`es le th\'eor\`eme de Stone-Von Neumann, les deux
repr\'esenta\-tions $\pi_{E_2^*}$ et
 $\pi_{{R}}$ de $H_2$ sont \'equivalentes ; si bien que
  la repr\'esentation m\'etap\-lectique $S_{E_2^*}$ peut \^etre r\'ealis\'ee  dans
 l'espace $\mathcal H_r$ de $\pi_{{R}}$.
On note  $ S_{R}$ cette r\'ealisation.
On pose
\begin{eqnarray*}
Sp(V_2)(r_{+,-}) = Sp(V_2)(r_+) \cap Sp(V_2)(r_-) \mbox{ et }
Mp(V_2)(r_{+,-}) =p_{V_2}^{-1}(Sp(V_2)(r_{+,-})),
\end{eqnarray*}
o\`u $p_{V_2}$ est la projection canonique de $Mp(V_2)$ dans $Sp(V_2)$. 
 Alors
$Sp(V_2)(r_{+,-}) $ (resp. $Mp(V_2)(r_{+,-})$)  est un
sous-groupe compact ouvert de $Sp(V_2)$ (resp. $Mp(V_2)$). Pour tout $x
\in Sp(V_2)(r_{+,-}) $, on consid\`ere  l'op\'erateur
$\sigma_{R}(x)$ agissant dans $\mathcal H_r$~:
\begin{eqnarray} \label{metaplectique}
  \sigma_{R}(x)\var(v) &=& \var(x^{-1}.v),  \var \in \mathcal H_r, v \in V_2.
\end{eqnarray}
Alors, on
a~:
\begin{eqnarray}
\label{Per} \sigma_{R}(x)\pi_{{R}}(y)
\sigma_{R}(x)^{-1}&=&\pi_{{R}}(x(y)) \mbox{, pour tout } y \in H_2
.
\end{eqnarray}
 De plus,
$\sigma_{R}(x)$ est un op\'erateur unitaire de $\mathcal H_r$. Il
en r\'esulte que la formule (\ref{metaplectique}) d\'efinit une
repr\'esentation unitaire de $Sp(V_2)(r_{+,-})$, v\'erifiant
l'\'egalit\'e (\ref{Per}). Si bien qu'il existe  un caract\`ere
unitaire $\tau_R$ de $Mp(V_2)(r_{+,-})$ v\'erifiant
  $$S_{R}(\tilde x) = \tau_R(\tilde x) \sigma_{R}(x),
\mbox{ pour tout } \tilde x \in Mp(V_2)(r_{+,-}).$$ D'apr\`es
(\cite{[Moe-Vi-Wa]}, Lemme II.10), on a~:  $\tau_R(\tilde x) = \pm
1, \mbox{ pour tout } \tilde x \in Mp(V_2)(r_{+,-}).$
\end{subsubsection}
\begin{subsubsection}{} \label{x24}
 Soit $\var \in
C_c^\infty(T_\e H_2)$. On a : 
$$\mbox{ Tr}(\mathcal I_{E_{2}^*, \tilde \tau}(\var d_{G_2}x)) = \mbox{ Tr}\left(\int_{\t_{\e|n_F|_p}} \tilde \tau(\exp_{({T}^{V_1})^{V_2}}(X))
   \otimes \mathcal T_\var( \exp_{{T}^{V_2}}(X))d_\t X\right), $$
   o\`u 
   $\mathcal T_\var( \exp_{{T}^{V_2}}(X))$
est
l'op\'erateur de rang fini, agissant dans $\mathcal H_r$,
 $$\mathcal T_\var( \exp_{{T}^{V_2}}(X)) = \int_{H_2}
\var(\exp_T(X)y) S_R(\eta_{T}^{V_2}(\exp_{{T}^{V_2}}(X))) \pi_R(y) d_{H_2}y.$$
Soit $X \in
\t_{\e|n_F|_p}$, on pose $x = \exp_T(X)$,
$\tilde{{x}}=\exp_{{T}^{V_2}}(X)$, et   $\var_x(y)= \var(xy)$,  $y \in H_2$. 

\noindent On note $T_{2}$ l'image de $T$ dans $Sp(V_2)$ par
 la repr\'esentation adjointe.
 Pour  une  base symplectique  convenable
$(e_1, \ldots, e_m, f_1, \ldots , f_m)$ de $V_2$ et en vertu   des
conditions impos\'ees sur ${\bf a}_G$ (Condition (\ref{cond})), on a : 
$$T_{2, \e} \subset Sp(V_2)(r_{+,-}). $$ 
Si bien que l'on a : 
$$S_R(\eta_{T}^{V_2}(\tilde{{x}}))= \sigma_{R}(x_{|V_2}).$$
Maintenant, on a, pour tout $\alpha \in \mathcal H_r$,
$$
\mathcal T_\var(\tilde x)  \alpha(v)
= \int_{{V_2}/r} A_{\var}(v , w) \alpha(w) d_{{V_2}/r}\dot w ,
$$
 o\`u $
 A_{\var}(v , w)
= \int_{  r \times k}  \var_x((x^{-1}.v -(w +\gamma), t))\varsigma_{a_0}( t - \frac{1}{2} B(w +
\gamma, x^{-1}.v - \gamma))
d_r\gamma d\mu(t),
$
pour tout $ v , w \in {V_2}$. Autrement dit, $ \mathcal T_\var(\tilde x)$ est un op\'erateur \`a noyau (\'egal \`a $ A_{\var}$).
 D'apr\`es le th\'eor\`eme de
Mercer, on a~:
 $$\mbox{Tr}(\mathcal T_\var(\tilde x)) =
\int_{{V_2}/r}A_{\var}(v , v) d_{{V_2}/r}\dot v.$$
 Comme $  A_{\var}(v +\gamma  , v +\gamma  ) =  A_{\var}(v , v),$
pour tout $ v  \in {V_2}$ et $\gamma \in r$, on a :
\begin{eqnarray*}
\mbox{Tr}(\mathcal T_\var(\tilde x))& =&|a_0|_p^{\frac{\dim {V_2}}{2}} \int_{V_2}
A_{\var}(v , v) d_{V_2}v\\
 &=& |a_0|_p^{\frac{\dim {V_2}}{2}} \int_{{V_2}}\int_{ r \times k} \var_x((x^{-1} - 1). v -
\gamma,t)\\
 && \times \varsigma_{a_0}(t- \frac{1}{2}B(v +\gamma , x^{-1}.v -
\gamma)) d\mu(t) d_r\gamma  d_{V_2}v.
 \end{eqnarray*}
On suppose, d\'esormais, que $\det(1-x)_{V_2} \neq 0$. On pose
$$q_{\pm x} = \frac{1+\pm x_{V_2}}{2}(1-\pm x_{V_2})^{-1}$$ et on d\'esigne par
$Q_{\pm x}$ la forme quadratique sur ${V_2}$ d\'efinie par $$Q_{\pm
x}(v) = B(q_{\pm x}v,v) \, , \,  v \in {V_2}. $$
Alors, d'apr\`es ( \cite{[M3]}, Lemme 39.2.2), on a : 
 \begin{eqnarray*}
\mbox{Tr}(\mathcal T_\var(\tilde x))
 & = & |a_0|_p^{\frac{\dim V_2}{2}}|\det(1-x)_{V_2}|_p^{-1}
 \left\{\int_{r}\varsigma_{a_0}( \frac{1}{2}Q_x(\gamma))d_{r}\gamma\right \}\times \\
 &&\times \int_{{V_2}  \times k} \var_x(v, t)\varsigma_{a_0}(t - \frac{1}{2}Q_x(v))
 d\mu(t)  d_{V_2}v.
 \end{eqnarray*}
Comme  $\frac{1+x_{V_2}}{2}r \subset r$ et  $|\det
(\frac{1+x_{V_2}}{2})_{V_2}|_p = 1$ donc 
$\frac{1+x_{V_2}}{2}r= r$. Il s'en suit que $$\varsigma_{a_0}\left(
\frac{1}{2}B\left(\left(1-x_{V_2}\right)\left(\frac{1+x_{V_2}}{2}\right)^{-1}v, v\right)\right)= 1\, , \,
\mbox{ pour tout } v \in r. $$ En utilisant les r\'esultats de (\cite{[M3]},
Section 39.2.3), on a~:
$$
|a_0|_p^{\frac{\dim V_2}{2}}\int_{r}\varsigma_{a_0}(
\frac{1}{2}Q_x(\gamma))d_{r}\gamma
=|\det(1-x)_{V_2}|_p^{\frac{1}{2}} \overline{\gamma(a_0Q_{-x})}.
$$
On note $u= \frac{- 1 + x_{V_2} }{\log (x_{V_2}) }
(\frac{1+x_{V_2}}{2})^{-1}$.
 Si $ \lambda \in \bar k$ est une
valeur propre de $x_{V_2}$ alors $$\left|\frac{- 1 + \lambda }{\log
(\lambda ) } \left(\frac{1+\lambda }{2}\right)^{-1}-1\right|_{ \bar k} < \e .$$
Il s'ensuit que $u \in GL(V_2)_\e$.
Ainsi,  $u$ poss\`ede  une  racine carr\'ee   qui commute avec
$\log(x_{V_2})$, \`a savoir  $$u^{\frac{1}{2}}=
\exp_{GL(V_2)}(\frac{1}{2}\log(u)).$$
Maintenant, par un  calcul \'el\'ementaire 
 (\cite{[M3]}, Paragraphes 39.2.5 et  39.2.6), on obtient~: 
$$\mbox{Tr}(\mathcal T_\var(\tilde x))= |a_0|_p^{ \frac{\dim V_2}{2}}
C^{V_2}_\var(X),$$ o\`u
$$
C^{V_2}_\var(X) = \int_{V_2}\int_{{V_2}\times k} (
\var\circ\exp_{G_2})(X +w
 +tE)
\varsigma_{a_0}(t+ B(w , v) +\frac{1}{2}B(X.v , v))d\mu( t)
d_{V_2}wd_{V_2}v.
$$
Comme 
 $\{x \in T_\e \ , \, \det(1-x)_{V_2} \neq
0\}$  est un ouvert  de $T_\e$  et
$T_\e \smallsetminus \{x \in T_\e \ , \, \det(1-x)_{V_2} \neq
0\}$ est de mesure nulle par rapport \`a la mesure de Haar $d_Tx$ sur $T$, 
on a : 
\begin{eqnarray*}
\mbox{ Tr}(\mathcal I_{E_{2}^*, \tilde \tau}(\var d_{G_2}x))
 & =&   |a_0|_p^{ \frac{\dim V_2}{2}} \int_{\t_{\e|n_F|_p}}
 \mbox{Tr}(\tilde \tau(\exp_{({T}^{V_1})^{V_2}}(X)))C^{V_2}_\var(X)d_\t X\\ &=&
\dim \tau|a_0|_p^{ \frac{\dim V_2}{2}}
\int_{\t_{\e|n_F|_p}}\varsigma(<\lambda_\tau, X>)C^{V_2}_\var(X)d_\t
X\\
 &=& \dim \tau\int_{\O_{E_{2}^*,\lambda_\tau}} (\var\circ
\exp_{G_2} d_{\g_2}X)\, \widehat{}_{\g_2}(l)
 d\mu_{\O_{E_{2}^*,\lambda_\tau}}(l)\\
 &=&  (\dim \tau) I.
 \end{eqnarray*}
\end{subsubsection}
\end{subsection}
\begin{subsection}{}\label{mlk}
 Dans cette section, on se donne  un id\'eal ab\'elien $\u$, $G$-invariant,  non nul  et
 inclus dans ${}^u\g$. 
 On reprend les notations des num\'eros \ref{ddv} et \ref{repr}.
 Fixons $0 < \e < {\bf a}_G$. Soit  
$\lambda_\tau \in \g(g)^*$ associ\'e \`a $\tau$. 
D'apr\`es les r\'esultats du paragraphe \ref{kh22},
le  prolongement, not\'e   de m\^eme, de  $\lambda_\tau $ en une forme lin\'eaire sur $\h(h)$ 
par $g_{|\u}$ sur $\u$, 
 est associ\'e  \`a $\tilde \tau$. Remarquons que si  $\O_{g, \lambda_\tau}$ est
ferm\'ee dans $\g^*$ alors $\O_{h, \lambda_\tau}$ (resp. $\O_{g_1, \lambda_\tau}$) est ferm\'ee dans $\h^*$ (resp. $\g_1^*$).
\begin{pro}\label{ymlk} On suppose que   $\O_{g, \lambda_\tau}$ est
ferm\'ee dans $\g^*$.
 Si le th\'eor\`eme \ref{yhd} est vrai pour $(G_1, g_1,
\tilde \tau,\lambda_\tau, \e )$ alors il est vrai pour $(G, g,
\tau, \lambda_\tau, \e)$.
\end{pro}
\begin{dem} 
 Soit $\var \in C_c^\infty(G_\e)$.  Alors,  $\pi_{g,
\tau}(\var d_Gx)$  est un op\'erateur
\`a noyau continu $K_\var$, d\'efini par $$K_{\var}(x, y) =
\Delta_{G}(y)^{-1} \int_H \var(xzy^{-1}) \Delta_{H,
G}(z)^{-\frac{1}{2}} \pi_{h, \tilde \tau}(z) d_Hz,$$ $d_Hz$
\'etant une mesure de Haar sur $H$ (voir \cite{[B]}, Chapitre V).    Pour $x \in G$, on
d\'esigne par $\var_H^x $ la fonction sur $H$ d\'efinie par
$$\var_H^x(y)= \var (xyx^{-1}), \, y \in H.$$ 
Alors, $\var_H^x \in C_c^\infty(H)$ et on a : 
$$\mbox{ support} (\var_H^x ) \subset G_\e \cap H = H_\e.$$
 Puisque
$(\Delta_{H, G})_{|H_\e} = 1$, on a :
 $$K_{\var}(x, x) =  
\Delta_{G}(x)^{-1}\pi_{h, \tilde \tau}(\var_H^x d_Hz).$$ Choisissons une mesure de Haar $d_Qx$ sur $Q$, notons $d_\q X$
la mesure de Haar sur $\q$ tangente \`a $d_Qx$ et posons, 
 $$\var_{G_1}^x(y) = \int_{Q} \var^x(yz)d_Qz = \int_\q
\var^x(y \exp(X))d_\q X, \quad y \in G_1.$$ Alors, 
$\var_{G_1}^x \in C_c^\infty(G_1)$, $\mbox{ support}( 
\var_{G_1}^x) \subset H_\e/Q =: G_{1, \e}$, et on a :
$$K_{\var}(x, x) = \Delta_{G}(x)^{-1}\pi_{g_1, \tilde \tau}
(\var_{G_1}^x d_{G_1}z),$$ $ d_{G_1}z$ \'etant la mesure de Haar
sur $G_1$, quotient de $d_Hx$ par $d_Qx$.

\noindent Maintenant, on suppose que   $\O_{g, \lambda_\tau}$ est
ferm\'ee dans $\g^*$. Alors, $\O_{g}$ est ferm\'ee dans $\g^*$ et d'apr\`es 
(\cite{[M3]}, Th\'eor\`eme 37), $\pi_{g, \tau}$ est admissible.
 Il en r\'esulte que l'on a : 
$$
\Theta_{g, \tau}(\var d_Gx)=\int_{G/H}\mbox{Tr} K_{\var}(x, x)d_{G/H}\dot x.
$$
 Supposons que le th\'eor\`eme \ref{yhd} est vrai pour $(G_1, g_1,
\tilde \tau,\lambda_\tau, \e )$. Alors, pour tout $x \in G$, on a :
\begin{eqnarray} \label{form}
\mbox{Tr} K_{\var}(x, x) &=&\Delta_{G}(x)^{-1} \dim \tau
\int_{\O_{g_1, \lambda_\tau}} (\var_{G_1}^x \circ \exp_{G_1}
d_{\g_1}X)\, \widehat{}_{\g_1} \,(l)  d\mu_{\O_{g_1,
\lambda_\tau}}(l),
\end{eqnarray}
o\`u $d_{\g_1}X = d_\h X / d_\q X$ et $d_\h X$ est  la
mesure de Haar sur $\h$ tangente \`a $d_H z$. Etant donn\'e que  $\O_{h,
\lambda_\tau}$   est contenu dans
$\q^\perp $, elle s'identifie canoniquement \`a
 l'orbite co-adjointe  $\O_{g_1, \lambda_\tau}$.
On a, pour tout $l \in \O_{h, \lambda_\tau}$,
\begin{eqnarray*}
(\var_{G_1}^x \circ \exp_{G_1} d_{\g_1}X)\, \widehat{}_{\g_1} \, (l)
&=&
 \int_{\g_1} (\var_{G_1}^x \circ \exp_{G_1})(X) \varsigma(<l, X>)
 d_{\g_1}X\\
 &=&\int_{\g_1} \{\int_\q \var^x(\exp_H(X)\exp_H(Y))d_\q Y\} \varsigma(<l, X>)
 d_{\g_1}X\\
 &=&\int_{\g_1} \int_\q \var^x(\exp_H(X)\exp_H( \frac{e^{ad X} - 1}{adX}.Y))d_\q Y \\ && \times \varsigma(<l, X>)
 d_{\g_1}X.
\end{eqnarray*}
Cependant, pour tout $X \in \h_{\e|n_F|_p}$, $Y \in \q$, on a, en
r\'earrangeant convenablement les termes dans la formule de
Campbell Hausdorff, 
\begin{eqnarray}\label{duf}
\exp_H(X) \exp_H(\frac{e^{ad X} - 1}{adX}.Y) = \exp_{H} (X + Y).
\end{eqnarray}
Si bien que l'on a : 
\begin{eqnarray*}
(\var_{G_1}^x \circ \exp_{G_1} d_{\g_1}X)\, \widehat{}_{\g_1}  (l)
&=&
 \int_{\h} (\var^x \circ \exp_{H})(X) \varsigma(<l, X>)
 d_{\h}X \\
 &=& \int_{\h^\perp}(\var^x \circ \exp_{G} d_{\g}X)\, \widehat{}_{\g}  (\tilde l +
 t)d_{\h^\perp}t\\
 &=& |\det Adx_\g|_p^{-1} \int_{\h^\perp}(\var \circ \exp_{G} d_{\g}X)\, \widehat{}_{\g} (x.(\tilde l +
 t))d_{\h^\perp}t,
\end{eqnarray*}
o\`u $d_{\h^\perp}t$ est la mesure de Haar sur $\h^\perp$ duale de
la mesure de Haar $d_\g X/ d_\h X$ sur $\g / \h$ et $\tilde l$ est
un \'el\'ement  de $\g^*$ dont la restriction \`a $\h$ est $l$. 

\noindent Il r\'esulte de ce qui pr\'ec\`ede que l'on a : 
\begin{eqnarray*}
&& \int_{G/H}\mbox{Tr} K_{\var}(x, x)d_{G/H}\dot x\\ &&\qquad \qquad = \dim \tau
\int_{G/H}\int_{\O_{h, \lambda_\tau}} \int_{\h^\perp}(\var \circ \exp_{G} d_{\g}X)\, \widehat{}_{\g} (x.(\tilde l +
 t))d_{\h^\perp}t  d\mu_{\O_{h, \lambda_\tau}}(l)d_{G/H}\dot x \\
 && \qquad \qquad = \dim \tau \int_{\O_{g,
\lambda_\tau}} (\var \circ \exp_{G} d_{\g}X)\, \widehat{}_{\g}  (l)
d\mu_{\O_{g, \lambda_\tau}}(l),
\end{eqnarray*}
o\`u la derni\`ere \'egalit\'e d\'ecoule de la proposition \ref{ymeet}.
\end{dem}
\end{subsection}
\end{section}
\begin{section}{}\label{B}
Dans ce paragraphe, nous donnons une d\'emonstration du
th\'eor\`eme \ref{yhd} au voisinage d'un \'el\'ement semi-simple
$s$ de $G$.
\begin{subsection}{La d\'emonstration du th\'eor\`eme  \ref{yhd} dans
 le premier cas de l'alternative du lemme \ref{ym}}\label{A}
   Le cas o\`u ${}^uG $ est trivial a \'et\'e trait\'e
dans la section \ref{rtl}. Il reste \`a montrer le th\'eor\`eme \ref{yhd} dans les conditions  du paragraphe \ref{ma}.
  \begin{subsubsection}{}\label{Med-Adem} On se place dans les conditions du paragraphe \ref{ma}. Soit $s$ un \'el\'ement semi-simple de $G$. 
  On peut supposer sans restriction que $s \in T$.  On pose~:
$$ V_{1,s} =  V(s) \mbox{ , } V_{2,s} =   (1 - s)V \mbox{ , } \h_{i,s} =
V_{i,s} + kE \mbox {, } H_{i,s} = \exp(V_{i,s} + kE) \mbox {, } E^*_{i,s}=
E^*_{a_0|\h_{i,s}}  \mbox{, }$$
$ i = 1\mbox{, } 2$.
 On se donne un lagrangien $\l_i$ de
 $V_{i,s}$, $i= 1,2$. Alors $ \goth l_i =\l_i \oplus kE $
 est une polarisation en $E^*_{i,s}$.
On consid\`ere $(\pi_{E^*_{i,s}}, \mathcal H_i)$ la repr\'esentation
unitaire irr\'eductible de $H_{i,s}$
 associ\'ee \`a $E^*_{i,s}$ par la
m\'ethode des orbites de Kirillov  et r\'ealis\'ee dans l'espace de Hilbert $
\mathcal H_i=  \mathcal H_{\goth l_i}$, et soit
$S_{E^*_{i,s}}$ la repr\'esentation m\'etaplectique  de $Mp(V_{i,s})$
associ\'ee au caract\`ere $\varsigma_{a_0}$ et r\'ealis\'ee dans
le m\^eme espace de
 Hilbert $\mathcal H_i$.

L'application $(y_1, y_2) \longmapsto y_1 y_2$ de $H_{1,s} \times H_{2,s}$
dans $H$ est un morphisme  surjectif de groupes de Lie, de noyau
\'egal $\Delta = \{(z, z^{-1}) , z \in Z\}$. La repr\'esentation
$\pi_{E^*_{1,s}} \otimes \pi_{E^*_{2,s}} $
 de $H_{1,s} \times H_{2,s}$ est unitaire irr\'eductible et est triviale sur $\Delta$.
   Par passage au quotient,
 on obtient  une repr\'esentation
unitaire
 irr\'eductible  de $ H$ \'equivalente \`a $\pi_{E_{a_0}^*}$.
 
On note $T(s)^{V_{i,s}}$ (resp. $T(s)^V$) l'extension m\'etaplectique
de $T(s)$ associ\'ee \`a l'action symplectique de $T(s)$ sur $V_{i,s}$
(resp. $V$) et $\eta_{T(s)}^{V_{i,s}}$ (resp. $\eta_{T(s)}^{V}$)
l'application naturelle de $T(s)^{V_{i,s}}$ (resp.  $T(s)^V$)   dans
$Mp(V_{i,s})$ (resp. $Mp(V)$). On choisit un rel\`evement $\tilde
s_i = (s, \psi_{i,s}) $ (resp. $\tilde s = (s, \psi_{s})$) de $s$
dans  $T(s)^{V_{i,s}}$ (resp.  $T(s)^V$).
 Avec ces notations, on~a~: \begin{center}$
 S_{E_{a_0}^*}(\eta_{T(s)}^{V}(\tilde s))=
 \psi_s( \l_1 +\l_2) \psi_{2,s}(\l_2)^{-1}
 \mbox{Id}_{ \mathcal H_1 }  \otimes S_{E^*_{2,s}}(\eta_{T(s)}^{V_{2,s}}(\tilde s_2)).$\end{center}
 \end{subsubsection}
\begin{subsubsection}{} Soit
$0<\e(s)<\inf\{a_G'(s), a_G''(s), \e\} $ tel que  $ G \times_{G(s)}
G(s)_{\e(s)}\longrightarrow \W_G(s , \e(s))$, \, $[x, y]\longmapsto
xs y x^{-1}$ soit un diff\'eomorphisme. La  restriction de  $\Theta_{E_{a_0}^*, \tau}$ \`a $\W_G(s
, \e(s) )$ est une fonction
g\'en\'eralis\'ee, $G$-invariante.  
D\'esignons par $\theta_{E_{a_0}^*, \tau,s}$ la fonction
g\'en\'eralis\'ee sur $G(s)_{\e(s)}$  qui lui est associ\'ee  par
le th\'eor\`eme \ref{yPTo}. 
Soit   $\beta \in
C_c^\infty(G(s)_{\e(s)})$ et $d_{G(s)}x = d_Txd_{V_{1,s}}v d\mu(t)$ une mesure de Haar sur $G(s)$.
Alors, il existe  un sous-groupe
compact ouvert $K$ de $G(s)$ tel que $\beta$ soit $K$-bi-invariant.
On se donne  $ \gamma \in C_c^\infty(G(s))$ tel que $\int_{G(s)}
\gamma(y) d_{G(s)}y=1$ dont le support  est contenu dans $K$.
 On se donne \'egalement $\alpha_s \in C_c^\infty(V_{2,s})$ tel que $\int_{V_{2,s}} \alpha_s(v)
d_{V_{2,s}}v = 1$. 
Comme l'application $\xi: G(s) \times V_{2,s} \longrightarrow T(s)H$, $(
 y, v)\longmapsto y \exp(v)
 $, est un diff\'eomorphisme et comme $T(s)H$ est ouvert dans  $G$, donc la fonction $\alpha = (\gamma\otimes \alpha_s ) \circ\xi^{-1}$ appartient \`a $ C_c^\infty(G)$. 
  On  munit $G$  de la  mesure de Haar $d_{G}x = d_Txd_{V_{1,s}}v d_{V_{2,s}}v d\mu(t)$ et on consid\`ere  l'application $  \Psi : G
\times G(s)_{\e(s)} \longrightarrow  G$, \, $ (x , y) \longmapsto
xs y x^{-1}$. On  obtient par la suite l'\'el\'ement $\Psi_{\alpha\otimes \beta}$ de $C_c^\infty(\W(s , \e(s)))$ donn\'e par la proposition \ref{yHarish}. On a alors $$\theta_{E_{a_0}^*, \tau,s}(\beta d_{G(s)}y) = \Theta_{E_{a_0}^*, \tau}(\Psi_{\alpha\otimes \beta}d_Gx).$$
En utilisant la r\'ealisation de $\pi_{E_{a_0}^*, \tau}$ donn\'ee dans le paragraphe \ref{Med-Adem}, 
on a~:
\begin{eqnarray*}
\pi_{E_{a_0}^*, \tau}(\Psi_{\alpha\otimes \beta}d_Gx)&=& \int_{\t_{\e(s)|n_F|_p}}\tau(\tilde s \exp_{T(s)^{V}}(X))
\otimes J(\tilde s\exp_{T(s)^{V}}(X)) d_{\t}X,
\end{eqnarray*}
o\`u
 $$\exp_{T(s)^{V}}(X)= (\exp_{T(s)}(X), \phi),\,
\exp_{T(s)^{V_{i,s}}}(X)= (\exp_{T(s)}(X), \phi_{i,s}), \, i= 1, 2, $$
$$J(\tilde s\exp_{T(s)^{V}}(X))= (\psi_s \psi_{ 1,s}^{-1}
\psi_{2,s}^{-1})( \phi \phi_{1,s}^{-1} \phi_{2,s}^{-1})\psi_{ 1,s}( \l_1)J_1(\exp_{T(s)^{V_{1,s}}}(X)) \otimes J_2(\tilde
s_2\exp_{T(s)^{V_{2,s}}}(X))$$ $J_1(\exp_{T(s)^{V_{1,s}}}(X))$ \'etant
l'op\'erateur
\begin{center}$ 
\int_{V_{1,s}\times
k}\beta(\exp_{T(s)}(X)\exp(Y+tE))S_{E_{1,s}^*}(\eta_{T(s)^{V_{1,s}}}(\exp_{T(s)^{V_{1,s}}}(X)))
\pi_{E_{1,s}^*}(\exp(Y+tE))d_{V_{1,s}}Y d\mu(t),
$\end{center}
et $J_2(\tilde s_2\exp_{T(s)^{V_{2,s}}}(X))=
J_{\alpha_s}(\eta^{V_{2,s}}_{T(s)}(\tilde s_2\exp_{T(s)^{V_{2,s}}}(X)))$, o\`u ce dernier op\'erateur est donn\'e par la formule (\ref{op}).

\noindent Il s'en suit  que l'on a :
\begin{eqnarray*}
\theta_{E_{a_0}^*, \tau,s}(\beta d_{G(s)}y) &=&
\psi_s( \l_1 + \l_2)
 \psi_{ 2,s}(\l_2)^{-1}
\int_{\t_{\e(s)|n_F|_p}}\mbox{Tr}(\tau(\tilde s
\exp_{T(s)^{V}}(X))) \\ && \times \mbox{Tr}(J_1(\exp_{T(s)^{V_{1,s}}}(X)) )
 \mbox{Tr}(J_2(\tilde s_2\exp_{T(s)^{V_{2,s}}}(X)))    d_\t X.
\end{eqnarray*}
D'apr\`es la proposition \ref{yexc} et compte tenu des conditions impos\'ees sur $a_G''(s)$, on a, pour tout $ X \in  \t_{\e(s)|n_F|_p}$,  $$ \psi_s( \l_1 + \l_2)
\psi_{2,s}(\l_2)^{-1}\mbox{Tr}(J_2(\tilde
s_2\exp_{T(s)^{V_{2,s}}}(X)))=|\det(s-1)_{V_{2,s}}|_p^{-\frac{1}{2}}\Phi_{a_0}(\eta_{T(s)}^{V}(\tilde
s)).$$ Il en r\'esulte que  l'on a~:
 \begin{eqnarray*}
\theta_{E_{a_0}^*, \tau,s}(\beta d_{G(s)}y)
&=&|\det(s-1)_{V_{2,s}}|_p^{-\frac{1}{2}}\Phi_{a_0}(\eta_{T(s)}^{V}(\tilde
s)) \int_{\t_{\e(s)|n_F|_p}}\mbox{Tr}(\tau(\tilde s
\exp_{T(s)^{V}}(X)))\\ && \times
\mbox{Tr}(J_1(\exp_{T(s)^{V_{1,s}}}(X)) )    d_\t X.
\end{eqnarray*}
Cependant, on a :
 $$\mbox{Tr}(\tau(\tilde s \exp_{T(s)^{V}}(X)))= \mbox{Tr}\tau(\tilde s)
\varsigma(<\lambda_\tau, X>), \mbox{ pour tout } X \in
\t_{\e(s)|n_F|_p}. $$  Si bien  que  l'on obtient :
\begin{eqnarray*}
\theta_{E_{a_0}^*, \tau,s}(\beta d_{G(s)}y)
&=&|\det(s-1)_{V_{2,s}}|_p^{-\frac{1}{2}}
\Phi_{a_0}(\eta_{T(s)}^{V}(\tilde s)) \mbox{Tr}\tau(\tilde s)
\int_{\t_{\e(s)|n_F|_p}}\varsigma(<\lambda_\tau, X>) \\ && \times
\mbox{Tr}(J_1(\exp_{T(s)^{V_{1,s}}}(X)) )d_\t X\\
&=& |\det(s-1)_{V_{2,s}}|_p^{-\frac{1}{2}}
\Phi_{a_0}(\eta_{T(s)}^{V}(\tilde s)) \mbox{Tr}\tau(\tilde s) \times \\ && \times \int_{\O_{E_{1,s}^*,\lambda_\tau}}(\beta \circ \exp
d_{\g(s)}Y)\,\widehat{}_{\g(s)} \, (l)
d\mu_{\O_{E_{1,s}^*,\lambda_\tau}}(l),
\end{eqnarray*}
o\`u $d_{\g(s)}Y$ est la mesure de Haar sur $\g(s)$ tangente \`a  la mesure de Haar $d_{G(s)}y$ sur $G(s)$.
\end{subsubsection}

\end{subsection}
\begin{subsection}{La d\'emonstration du th\'eor\`eme  \ref{yhd} dans
 le deuxi\`eme cas de l'alternative du lemme \ref{ym}}
 On se place dans la situation indiqu\'ee dans
la section \ref{mlk}.
 Soit $s \in G$, semi-simple et  $0<\e_G(s)<\inf\{a_G'(s), a_G''(s), \e\} $ tel que l'application $
G \times_{G(s)} G(s)_{\e_G(s)}\longrightarrow \W_G(s , \e(s))$, \,
$[x, y]\longmapsto  xs y x^{-1}$ soit un diff\'eomorphisme.
\begin{subsubsection}{}\label{reduct} Dans ce paragraphe, on suppose que
l'id\'eal $\u$ est contenu dans $\ker g$. Dans ce cas,  $\q = \u$,
$G_1 = G/Q$ et $\g_1= \g/\q$ ; via l'identification $\q^\perp =
\g_1^*$, on a : $g= g_1$, $\O_{g, \lambda_\tau} = \O_{g_1,
\lambda_\tau}$ et les mesures de Liouville $d\mu_{\O_{g,
\lambda_\tau}}$ et $d\mu_{\O_{g_1, \lambda_\tau}}$ sont \'egales.
D\'esignons par  $p_1$  la projection canonique de $G$ sur $G_1$
et par la m\^eme lettre (i.e. par  $p_1$) la projection canonique
de $\g$ sur $\g_1$.

{\it Remarque  : On a : $G_1(g_1)= p_1(G(g))$ et l'application
 $G(g)^\g \longrightarrow  G_1(g_1)^{\g_1}$, $(x, \psi) \longmapsto  (p_1(x), \psi)$
  est un morphisme
surjectif, que nous noterons encore $p_1$,  de noyau $Q$.
 Si bien que $\tau = \tilde \tau \circ p_1$ et
\begin{eqnarray}\label{njuu}
\pi_{g, \tau}& =& \pi_{g_1,\tilde \tau  } \circ p_1.
\end{eqnarray}
}
 On pose $s_1 = p_1(s)$.
Alors, $s_1$ est un \'el\'ement semi-simple de $G_1$. De plus, on
a : $\O_{g, \lambda_\tau, s} =\O_{g_1, \lambda_\tau, s_1}$ et
 $p_1(\W_G(s,\e_G(s) ))$ est un ouvert $G_1$-invariant, contenu dans
 $\W_{G_1}(s_1, \e_G(s))$. D'autre part,
on choisit une mesure de Haar $d_G x$ (resp. $d_Qx$) sur $G$
(resp. $Q$)
 et $d_\g X$ (resp. $d_\q X$) la mesure de Haar sur $\g$
  (resp. $\q$)  tangente \`a $d_G x$ (resp. $d_Qx$).  On munit $G_1$ (resp. $\g_1$)  de
la mesure de Haar quotient $d_{G_1}x =d_Gx/d_Qx $ (resp.
$d_{\g_1}X = d_\g X /d_\q X$). Soit $\var \in
C_c^\infty(\W_G(s,\e_G(s)))$. On pose,  $$
\var_{G_1}(x) = \int_Q \var(xy) d_Qy = \int_\q \var(x \exp(Y))d_\q Y, \, \,  x \in G.
$$ Alors $\var_{G_1} \in C_c^\infty(p_1(\W_G(s,\e_G(s) )))$ et on a~:
\begin{eqnarray}\label{israala}
\Theta_{g, \tau}(\var d_Gx ) &=& \Theta_{g_1,\tilde \tau}(\var_{G_1}
d_{G_1}x ).
\end{eqnarray}
On suppose d\'esormais  que le th\'eor\`eme \ref{yhd} est vrai pour
$G_1$ et $\pi_{g_1, \tilde \tau}$. 

\noindent $-$ Ou bien   $\O_{g, \lambda_\tau,
s} = \emptyset $, dans ce cas $\O_{g_1, \lambda_\tau, s_1} =
\emptyset $ et la formule (\ref{israala})
 montre que le th\'eor\`eme~\ref{yhd} est vrai   pour $G$ et $\pi_{g,  \tau}$.

\noindent $-$ Ou bien $\O_{g,
\lambda_\tau, s} \neq \emptyset $. On note  $\theta_{g_1, \tilde \tau, s_1}$ la
 fonction g\'en\'eralis\'ee $G_1(s_1)$-invariante sur $G_1(s_1)_{\e_G(s)}$ associ\'ee \`a
   $\Theta_{g_1,\tilde \tau}$ par le th\'eor\`eme \ref{yPTo}.
On se donne une mesure de Haar $d_{\g(s)}X$ (resp.
$d_{\g_1(s_1)}X$) sur $\g(s)$ (resp. $\g_1(s_1)$) . Le groupe $
G(s)$ (resp. $G_1(s_1)$) est muni de la mesure de Haar $d_{G(s)}x$
(resp. $d_{G_1(s_1)}x$) tangente \`a $d_{\g(s)}X$ (resp.
$d_{\g_1(s_1)}X$). En utilisant  le th\'eor\`eme \ref{yPTo}, on a~:
\begin{eqnarray*}
 \Theta_{g_1,\tilde \tau} (\var_{G_1} d_{G_1} x)&=&
\int_{G_1/G_1(s_1)}\Delta_{G_1}(x)^{-1}\int_{
G_1(s_1)_{\e_{G}(s)}}
 \theta_{g_1, \tilde \tau, s_1}(y)
\var_{G_1}^{x}(s_1y)\\&& \times
|\det(1-(s_1y)^{-1})_{(\g_1)_{s_1}}|_pd_{G_1(s_1)}yd_{G_1/G_1(s_1)}\dot
x,
\end{eqnarray*}
o\`u $d_{G_1/G_1(s_1)}\dot x$ est la mesure $G_1$-invariante sur
$G_1/G_1(s_1)$ tangente \`a la mesure de Haar quotient $d_{\g_1}
X/ d_{\g_1(s_1)} X$ sur $\g_1/\g_1(s_1)$ et $\var_{G_1}^{x}$ est la
fonction sur $G_1$ d\'efinie par $\var_{G_1}^{x}(y)=
\var_{G_1}(xyx^{-1})$, \, $y \in G_1$. On a~:
\begin{eqnarray*}
 &&\Theta_{g_1,\tilde \tau} (\var_{G_1} d_{G_1} x)= |\det(1-s_1^{-1})_{(\g_1)_{s_1}}|_p
\int_{G_1/G_1(s_1)}\Delta_{G_1}(x)^{-1}\times \\ &&\quad 
  \int_{ \O_{g_1,\lambda_\tau,s_1}} (1_{\g_1(s_1)_{\e_{G}(s)|n_F|_p}}
\var_{G_1}^{x}(s_1 \exp( . )) d_{\g_1(s_1)}X)\,\widehat{}_{\g_1(s_1)}
\,
 \phi_{g_1,\tilde\tau,s_1}
d\mu_{\O_{g_1,\lambda_\tau, s_1}} d_{G_1/G_1(s_1)}\dot x.
\end{eqnarray*}
Maintenant, consid\'erons la d\'ecomposition $\q = \q(s) \oplus
\q_s$ et choisissons des mesures de Haar $d_{\q(s)}X$ et
$d_{\q_s}X$ sur $\q(s)$ et $\q_s$ respectivement de sorte que
$d_\q X =d_{\q(s)}Xd_{\q_s}X$.  En utilisant la formule (\ref{duf})
appliqu\'ee \`a $\g(s)$ et $\q(s)$,
 nous pouvons
\'ecrire, pour $x \in G$ et  $X \in \g$,
\begin{eqnarray*}
&& 1_{(\g_1(s_1)_{\e_{G}(s)|n_F|_p})}(p_1(X)) \var_{G_1}^{p_1(x)}(s_1
\exp(p_1 (X) ))= \\&& \quad |\det Adx_\q|_p
1_{(\g(s)_{\e_{G}(s)|n_F|_p} + \q_s)}(X)
\int_{\q(s)}\{\int_{\q_s}\var^{x}(s \exp( X + Y )\exp(Z))
d_{\q_s}Z \}d_{\q(s)}Y.
\end{eqnarray*}
Or, l'application $\g(s) \longrightarrow \g_1(s_1)$, \,
$X\longmapsto p_1(X)$, est lin\'eaire surjective de noyau \'egal
$\q(s)$. Ceci permet d'identifier $\g(s)/ \q(s)$ \`a $\g_1(s_1)$.
Pour $l \in \g_1(s_1)^*$, on a~:
\begin{eqnarray*}
&&(1_{(\g_1(s_1)_{\e_{G}(s)|n_F|_p})} \var_{G_1}^{p_1(x)}(s_1 \exp( . ))
d_{\g_1(s_1)}X)\,\widehat{}_{\g_1(s_1)}\, (l)
 =     |\det Adx_\q|_p \times \\&&\qquad \qquad\int_{\g(s)\times \q_s  } \var^{x}(s\exp( X )
 \exp(Z))  1_{(\g(s)_{\e_{G}(s)|n_F|_p})}(X) \varsigma(<l, X>)d_{\q_s}Zd_{\g(s)}X.
\end{eqnarray*}
Comme, pour tout $Z \in \q_s$ et $X \in \g(s)_{\e_{G}(s)|n_F|_p}$,
$$\exp(Z) s\exp(X) \exp(-Z) = s \exp(X) \exp(((s\exp(X))^{-1}-
1)Z),$$ on a~:
\begin{eqnarray*}
&&(1_{(\g_1(s_1)_{\e_{G}(s)|n_F|_p})} \var_{G_1}^{x}(s_1 \exp( .
)) d_{\g_1(s_1)}X)\,\widehat{}_{\g_1(s_1)}\, (l)
 =|\det Adx_\q|_p|\det(1- s^{-1})_{\q_s}|_p\\ && \qquad \qquad  \quad  \times
  \int_{\g(s)\times \q_s}1_{(\g(s)_{\e_{G}(s)|n_F|_p})}(X) \var^{x\exp(Z)}
 (s\exp( X )) \varsigma(<l, X>)d_{\q_s}Zd_{\g(s)}X\\ &&
\qquad \qquad  \quad =|\det Adx_\q|_p|\det(1- s^{-1})_{\q_s}|_p
\\ && \qquad \qquad  \qquad \qquad    \times
  \int_{\q_s}(1_{(\g(s)_{\e_{G}(s)|n_F|_p})} \var^{x\exp(Z)}(s\exp( . )d_{\g(s)}X))_{\g(s)}\widehat{}\, (l) d_{\q_s}Z.
\end{eqnarray*}
Cependant, $$|\det(1-s_1^{-1})_{(\g_1)_{s_1}}|_p|\det(1-
s^{-1})_{\q_s}|_p = |\det(1- s^{-1})_{\g_s}|_p,$$  
$$|\det Ad (p_1(x))_{\g_1}|_p|\det Adx_\q|_p = |\det Adx_\g|_p =
\Delta_G(x)^{-1}, \mbox{ pour tout } x \in G,$$ et  $$\phi_{g_1,\tilde \tau,s_1} =
\phi_{g,\tau,s}.$$
 Maintenant, en utilisant le fait que  $G_1 /G_1(s_1)
= G/G(s)\exp(\q_s)$, nous obtenons  de ce qui pr\'ec\`ede,
\begin{eqnarray*}
 \Theta_{g, \tau} (\var d_{G}x)&=& |\det(1-s^{-1})_{\g_{s}}|_p
\int_{G/G(s)}\Delta_{G}(x)^{-1} \\&& \times
  \int_{ \O_{g,\lambda_\tau,s}} (1_{(\g(s)_{\e_{G}(s)|n_F|_p})}
\var^{x}(s \exp( . )) d_{\g(s)}X)\,\widehat{}_{\g(s)} \,
 \phi_{g,\tau,s}
d\mu_{\O_{g,\lambda_\tau, s}} d_{G/G(s)}\dot x.
\end{eqnarray*}
C'est la formule cherch\'ee.
\end{subsubsection}
\begin{subsubsection}{}
Nous retournons d\'esormais \`a la situation g\'en\'erale du cas
consid\'er\'e en \ref{mlk}.  Fixons une mesure de Haar \`a gauche
$d_Gx$ (resp. $d_Hx$) sur $G$ (resp. $H$) et  notons  $d_\g X$ (resp.
$d_\h X$) la mesure de
 Haar tangente sur $\g$ (resp. $\h$) et  r\'ealisons la repr\'esentation
 induite $\pi_{g , \tau }
= Ind_H^G \pi_{h,  \tilde \tau} =  Ind_H^G (\pi_{g_1,  \tilde
\tau} \circ p_1) $ (voir les  formules (\ref{represe}) et
(\ref{repres}))
 en utilisant la mesure invariante $d_{G/H}\dot x$ sur $G/H$,
quotient de $d_Gx$ par $d_Hx$.
Soit $\var \in C_c^\infty(G)$. On sait que l'op\'erateur $ \pi_{g
,\tau }(\var d_Gx)$ est un op\'erateur \`a noyau continu donn\'e
par $$K_\var(x , y) =\Delta_G(y)^{-1} \int_{H}\Delta{_{H,
G}}(z)^{-\frac{1}{2}}\var(xzy^{-1})\pi_{h,  \tilde \tau}(z)
d_Hz,$$ et gr\^ace au th\'eor\`eme de Mercer, on a~:
\begin{eqnarray}\label{mercer}
\Theta_{g ,\tau }(\var d_Gx) &=& \int_{G/H} \Delta_G(x)^{-1}
\Theta_{h, \tilde \tau} (\Delta{_{H, G}}^{-\frac{1}{2}}\var_H^x
d_{H}z) d_{G/H}\dot x.
\end{eqnarray}
\end{subsubsection}
\begin{subsubsection}{} Dans cette section, on suppose que
$G.s \cap H = \emptyset$. Vu le choix de $\e_G(s)$, on a~: $\W_G(s
, \e_G(s)) \cap H = \emptyset$. Ainsi,   si le support de $ \var$
est inclus dans $ \W_G(s , \e_G(s))$ alors, pour tout $x \in G$,
$\var_H^x= 0$. Si tel est le cas,
 la formule (\ref{mercer})  donne
 $\Theta_{g ,\tau }(\var d_Gx)  =
0$. D'autre part, notre hypoth\`ese implique  aussi que $\O_{g,
\lambda_\tau, s} = \emptyset$ ; le th\'eor\`eme est donc d\'emontr\'e dans
ce cas.
\end{subsubsection}
\begin{subsubsection}{} Maintenant, on se place dans le cas o\`u  $G.s \cap H \neq
\emptyset$.  On peut supposer que $s \in H$.
On pose
$$G_{s , H} = \{x \in G \mbox { tel que } x^{-1}s x \in H \}.$$
Alors $G_{s , H}$ est une sous-vari\'et\'e ferm\'ee de $G$.
Le groupe produit $G(s) \times H$ op\`ere \`a gauche dans $ G_{s , H}$
 par la formule~: $(x ,h).y = xyh^{-1}$. De plus,  l'action conjugaison  de $H$ dans $G.s \cap
H$  se prolonge en une action de $G(s) \times H$ qui soit triviale sur $G(s) \times
\{1\}$.
Avec ces notations, l'application
$\eta : G_{s , H} \longrightarrow G.s \cap H$, $x \longmapsto x^{-1}sx$,  est $G(s) \times
H$-\'equivariante. Elle induit  une  bijection~:
 $$ \begin{array}[t]{clclc}
            {\overline{ \eta}}  : & G(s)\setminus G_{s , H}/H  &
\longrightarrow & H\setminus (G.s \cap H)  \\
                   &  G(s) xH & \longmapsto    & H.(x^{-1}s x) = \{
hx^{-1}s xh^{-1} , h \in H \}
                     \end{array}.$$
 En utilisant (\cite {[Ri]} et  \cite {[Bo-Se]}),
 on  montre que $G.s \cap H $ est une r\'eunion disjointe d'un  nombre  fini
de $H$-orbites : soit $x_1 = 1 , x_2, \ldots , x_m  \in G_{s , H}$
tels que $$ Gs \cap H = \bigsqcup_{i = 1}^m H.(x_i^{-1}sx_i) \mbox{
et }
 G_{s, H} = \bigsqcup_{i = 1}^m
G(s)x_iH.$$
 Pour  $ i \in \{1 , 2 , \ldots , m \} $, on pose $ s_i =
x_i^{-1}s x_i$. D'apr\`es (\cite{[M3]}, Proposition 43.3.1), on a
\begin{eqnarray}\label{l}
 \mathcal W_G(s , \e_G(s)) \cap H = \bigsqcup_{i = 1}^m {\mathcal
W_H}(s_i , \e_G(s)).
\end{eqnarray}
Soit $\var \in C_c^\infty(\mathcal W_G(s , \e_G(s) ))$. On pose,
pour
 $1 \leq i \leq m$,
$\var_{s_i} = 1_{{\mathcal W_H}(s_i , \e_G(s))}\var_{|H}$. En
utilisant ce qui pr\'ec\`ede, la formule (\ref{mercer}) s'\'ecrit
aussi~:
\begin{eqnarray}\label{mercers}
\Theta_{g ,\tau }(\var d_Gx) &=&\sum_{1 \leq i \leq m
}I_{s_i}(\var d_Gx),
\end{eqnarray}
o\`u
\begin{eqnarray}\label{Imercers}
 I_{s_i}(\var d_Gx)&=& \int_{G/H} \Delta_G(x)^{-1} \Theta_{h,
\tilde \tau} (\Delta{_{H, G}}^{-\frac{1}{2}}(\var^x)_{s_i} d_{H}z)
d_{G/H}\dot x.
\end{eqnarray}
\end{subsubsection}
\begin{subsubsection}{Description de $\O_{g, \lambda_\tau, s} $ en terme
des $\O_{h, \lambda_\tau, s_i}$}\label{Richa}  On
consid\`ere  les d\'ecompo\-si\-tions suivantes. $$\g^* =
\g({s_i})^* \oplus \g_{s_i}^* \mbox{ et } \h^* = \h({s_i})^*
\oplus \h_{s_i}^*.$$ On d\'esigne par $r$ l'application  de
restriction  de $\g^*$  \`a $\h^*$ et par $r_{s_i}$ la restriction
de $r$ \`a $\g({s_i})^*$.
 Si $E$ est une partie de
$\O_{h,\lambda_\tau, {s_i}}$, on pose $[E]_{s_i} =
G({s_i})r_{s_i}^{-1}(E)$. Alors, on le r\'esultat  suivant~:  soit $f \in \g^*$  dont la restriction \`a
 ${}^u\g$  co\"{\i}ncide avec  celle de $g$, et dont la restriction \`a
 $\g(g)$ est \'egale \`a $\lambda_\tau$.
\begin{lem}\label{yRicha}
On a~:
\begin{itemize}
    \item $ \quad  r^{-1}(\O_{h,\lambda_\tau})= H.f$
    \item $ \quad  r^{-1}(\O_{h,\lambda_\tau})\cap \g({s_i})^* = r_{s_i}^{-1}
    (\O_{h,\lambda_\tau, {s_i}})$
    \item $ \displaystyle \quad \O_{g,\lambda_\tau,  s} =
     \bigsqcup_{1\leq i \leq m}x_i.[\O_{h,\lambda_\tau, {s_i}} ]_{s_i}$
\end{itemize}
De plus, chaque $x_i.[\O_{h,\lambda_\tau, {s_i}} ]_{s_i}$ est
un ouvert de $\O_{g,\lambda_\tau,  s}$, r\'eunion finie de
$G(s)$-orbites. L'application $\omega \longmapsto [\omega]_{s_i}$
induit une bijection de $H({s_i})\backslash \O_{h,\lambda_\tau,
{s_i}}$ sur $G({s_i})\backslash [\O_{h,\lambda_\tau,
{s_i}}]_{s_i}$.
\end{lem}
\medskip
 L'analogue de ce r\'esultat dans le cas r\'eel  a
 \'et\'e d\'emontr\'e dans (\cite{[KT]}, Lemme 8.2.1). Leur d\'emonstration s'\'etend  \`a notre cas.
\end{subsubsection}

\begin{subsubsection}{} \label{lm} Dans ce num\'ero,  on suppose
 $\O_{g,\lambda_\tau, s} = \emptyset $. D'apr\`es le lemme
 pr\'ec\'edent, on a ${\O_{h,\lambda_\tau, s_i}}= \emptyset$,
 pour tout $i \in \{ 1, \ldots, m\}$. Supposons que le th\'eor\`eme \ref{yhd}
  est vrai pour $(H, h, \tilde \tau)$. Alors,
  $ I_{s_i}(\var d_Gx)=0 $, pour tout $i \in \{ 1, \ldots, m\}$.
  Si bien que
$ \Theta_{g , \tau}(\var d_Gx)=0$.
\end{subsubsection}
\begin{subsubsection}{} Dans ce paragraphe,
on suppose d\'esormais que  $\O_{g,\lambda_\tau, s} \neq \emptyset
$. Comme la fonction g\'en\'eralis\'ee $\Theta_{g , \tau}$ est
$G$-invariante, on peut supposer que $Ad^*s.f = f$.
On choisit une mesure de Haar $d_{\g(s)}X$ (resp. $d_{\h(s)}X$)
sur $\g(s)$ (resp. $\h(s)$). On munit chaque $\g(s_i)$ (resp.
$\h(s_i)$) de la mesure de Haar $d_{\g(s_i)}X$ (resp.
$d_{\h(s_i)}X$) image de $d_{\g(s)}X$ (resp. $d_{\h(s)}X$)  par
l'isomorphisme : $ x_i^{-1} : \g(s)\longrightarrow \g(s_i)$ (resp. $
x_i^{-1} : \h(s)\longrightarrow \h(s_i)$). Le groupe $ G(s_i)$ (resp.
$H(s_i)$) est muni de la mesure de Haar $d_{G(s_i)}x$ (resp.
$d_{H(s_i)}x$) tangente \`a $d_{\g(s_i)}X$ (resp. $d_{\h(s_i)}X$).
 Pour chaque $i \in \{ 1, \ldots,
m\}$, on d\'esigne par $\theta_{h,\tilde \tau, s_i}$ la fonction
g\'en\'eralis\'ee $H(s_i)$-invariante sur $H(s_i)_{\e_G(s)}$
associ\'ee la fonction g\'en\'eralis\'ee $\Theta_{h,\tilde \tau}$
par le th\'eor\`eme \ref{yPTo}. Alors, on a :
\begin{eqnarray*}
\Theta_{h,\tilde \tau}(\Delta{_{H, G}}^{-\frac{1}{2}}
 (\var^x)_{s_i} d_Hz)
 &=&
\int_{H/H(s_i)}\Delta_H(y)^{-1}\int_{H(s_i)_{\e_G(s)}}\theta_{h,\tilde
\tau, s_i}(z) \Delta{_{H, G}}^{-\frac{1}{2}}(s_iz)\\ && \times
\var^{xy}(s_iz)
|\det(1-(s_iz)^{-1})_{\h_{s_i}}|_pd_{H(s_i)}zd_{H/H(s_i)}\dot y,
\end{eqnarray*}
o\`u $d_{H/H(s_i)}\dot y$ est la mesure $H$-invariante sur $H/H(s_i)$
tangente \`a la mesure de Haar quotient $d_\h X/ d_{\h(s_i)} X$
sur $\h/\h(s_i)$.

Mais,  pour tout $z \in H(s_i)_{\e_G(s)}$, on a : $$\Delta{_{H,
G}}(z) = 1 \mbox{ et } |\det(1-(s_iz)^{-1})_{\h_{s_i}}|_p =
|\det(1-s_i^{-1})_{\h_{s_i}}|_p.$$ On pose $c(s_i) =\Delta{_{H,
G}}^{-\frac{1}{2}}(s_i)|\det(1-s_i^{-1})_{\h_{s_i}}|_p$.
 Alors, on peut
\'ecrire,
\begin{eqnarray*}
\Theta_{h,\tilde \tau}(\Delta{_{H, G}}^{-\frac{1}{2}}
 (\var^x)_{s_i} d_Hz)
 &=&c(s_i)
\int_{H/H(s_i)}\Delta_H(y)^{-1}\\ &&
\times\int_{H(s_i)_{\e_G(s)}}\theta_{h,\tilde \tau, s_i}(z)
\var^{xy}(s_iz) d_{H(s_i)}zd_{H/H(s_i)}\dot y.
\end{eqnarray*}
Par hypoth\`ese le th\'eor\`eme \ref{yhd} est vrai pour $(G_1, g_1,
\tilde \tau, \lambda_{\tau}, \e_G(s), p_1(s_i)$). En appliquant
les r\'esultats de la section  \ref{reduct}, on a~:
\begin{eqnarray*}
 I_{s_i}(\var d_Gx) & =& c(s_i)\int_{G/H}\Delta_G(x)^{-1}
\int_{H/H(s_i)}\Delta_H(y)^{-1} \, \, \times \\&&   \times  \int_{
\O_{h,\lambda_\tau,s_i}} (1_{\h(s_i)_{\e_G(s)|n_F|_p}} \var^{xy}(s_i
\exp( . )) d_{\h(s_i)}X)\,\widehat{}_{\h(s_i)} \, \, \times \\ &&
\times \phi_{h,\tilde \tau,s_i} d\mu_{\O_{h,\lambda_\tau,
s_i}}d_{H/H(s_i)}\dot y d_{G/H}\dot x.
\end{eqnarray*}
La fonction $ \phi_{h,\tilde \tau, s_i}$ (resp. $\phi_{g,\tau,
s_i}$) est constante sur  les $H(s_i)$-orbites (resp.
$G(s_i)$-orbites) donc les nombres $ \phi_{h,\tilde \tau,
s_i}(\omega)$, $\omega \in H(s_i)\backslash \O_{h,\lambda_\tau,
s_i}$ (resp. $ \phi_{g, \tau, s_i}(\omega)$, $\omega \in
G(s_i)\backslash \O_{g,\lambda_\tau, s_i}$) sont bien d\'efinis.
D'apr\`es la proposition \ref{ymlhgf}, on a : $$\phi_{h,\tilde
\tau, s_i}(\omega)
 = \phi_{g, \tau,
s_i}([\omega]_{s_i})\Delta_{H,
G}^{\frac{1}{2}}(s_i)|\det(1-s_i^{-1})_{(1-s_i)\g/\h}|_p, \, \omega
\in H(s_i)\backslash \O_{h,\lambda_\tau, s_i}.$$ Remarquons que
$$c(s_i)\Delta_{H,
G}^{\frac{1}{2}}(s_i)|\det(1-s_i^{-1})_{(1-s_i)\g/\h}|_p
=|\det(1-s_i^{-1})_{\g_{s_i}}|_p. $$
 Ainsi,

\begin{eqnarray} \label{bgg}
  I_{s_i}(\var d_Gx) &=& |\det(1-s_i^{-1})_{\g_{s_i}}|_p\sum_{\omega \in H(s_i)
  \backslash \O_{h,\lambda_\tau, s_i}}
\phi_{g,\tau , s_i}([\omega]_{s_i}) I^{s_i}_{\omega}(\var
d_G x)
\end{eqnarray}
o\`u on a pos\'e, pour $\omega \in H(s_i)\backslash
\O_{h,\lambda_\tau, s_i}$,
\begin{eqnarray*}
&&I^{s_i}_{\omega}(\var d_G x) =
\int_{G/H}\Delta_G(x)^{-1}\int_{H/H(s_i)}\Delta_H(y)^{-1} \times\\&& \quad  \qquad 
 \int_{ \omega}
((1_{\h(s_i)_{\e_G(s)|n_F|_p}}\var^{xy}(s_i\exp_{G(s_i)}(\, . \,
)))_{\h(s_i)} d_{\h(s_i)}X)\,\widehat{}_{\h(s_i)} \,
d\mu_{\omega}d_{H/H(s_i)}\dot y d_{G/H}\dot x.
\end{eqnarray*}
 D\'esignons par $\h(s_i)^\perp$ l'espace des formes lin\'eaires
 sur $\g(s_i)$ nulle sur $\h(s_i)$ (que l'on identifie canoniquement
 avec $(\g(s_i)/\h(s_i))^*$)
   et par $d_{\h(s_i)^\perp}t$ la mesure de Haar sur
$\h(s_i)^\perp$ duale de la mesure de Haar
$d_{\g(s_i)}X/d_{\h(s_i)}Y$ sur $\g(s_i)/\h(s_i)$. Soit $\omega
\in H(s_i)\backslash \O_{h,\lambda_\tau, s_i}$. On peut \'ecrire,
\begin{eqnarray*}
\quad I^{s_i}_{\omega}(\var d_G x) &=&
\int_{G/G(s_i)}
\Delta_G(x)^{-1} \{\int_{G(s_i)/H(s_i)}
\\&& \times \{\int_{\omega}\int_{\h(s_i)^\perp}
(1_{\g(s_i)_{\e_G(s)|n_F|_p}}\var^{x} (s_i\exp_{G(s_i)}(\, . \,
))d_{\g(s_i)}X)\,\widehat{}_{\g(s_i)}(y.(l_1 + t))  \\&&
d_{\h(s_i)^\perp}t d\mu_{\omega}(l)\}d_{G(s_i)/H(s_i)} \dot y \}
d_{G/G(s_i)}\dot x.
\end{eqnarray*}
En appliquant la proposition \ref{ymeet} aux orbites $\omega$ et
$[\omega]_{s_i}$ dans $\h(s_i)^*$ et $\g(s_i)^*$, on obtient,
\begin{eqnarray} \label{ghj}
&& I^{s_i}_{\omega}(\var d_G x) =\int_{G/G(s_i)}
\Delta_G(x)^{-1}\times \\ && \qquad \times \left\{\int_{[\omega]_{s_i}}
(1_{\g(s_i)_{\e_G(s)|n_F|_p}}\var^{x} (s_i\exp_{G(s_i)}(\, . \,
))d_{\g(s_i)}X)\, \widehat{}_{\g(s_i)} (l)
d\mu_{[\omega]_{s_i}}(l)\right\}d_{G/G(s_i)}\dot x. \nonumber
\end{eqnarray}
 Reportant la  formule (\ref{ghj}) dans
(\ref{bgg}), on obtient,
\begin{eqnarray*} 
 \quad I_{s_i}(\var d_G x) &=& |\det(1-s_i^{-1})_{\g_{s_i}}|_p
\int_{G/G(s_i)} \Delta_G(x)^{-1}\times \\ &&\times
\left\{\int_{[\O_{h,\lambda_\tau, s_i}]_{s_i}}
(1_{\g(s_i)_{\e_G(s)|n_F|_p}}\var^{x} (s_i\exp_{G(s_i)}(\, . \,
))d_{\g(s_i)}X)\, \widehat{}_{\g(s_i)} (l) \, \times\right. \\&&
\times \left. \phi_{g, \tau, s_i}(l)  d\mu_{[\O_{h,\lambda_\tau,
s_i}]_{s_i}}(l)\right\}d_{G/G(s_i)}\dot x 
\end{eqnarray*}
Maintenant, pour tout $l \in \g(s)^*$, on a : $x_i^{-1}.l \in
\g(s_i)^*$ et
\begin{eqnarray*}
&&(1_{\g(s_i)_{\e_G(s)|n_F|_p}}\var^{x} (s_i\exp_{G(s_i)}(\, . \,
))d_{\g(s_i)}X)\, \widehat{}_{\g(s_i)} (x_i^{-1}.l) = \\&& \qquad \qquad
\qquad \qquad \qquad \qquad c_i(1_{\g(s)_{\e_G(s)|n_F|_p}
}\var^{xx_i^{-1}} (s\exp_{G(s)}(\, . \,
))d_{\g(s)}X)\, \widehat{}_{\g(s)} (l),
 \end{eqnarray*}
  o\`u $c_i =
|\det(x_i^{-1} : \g(s) \longrightarrow \g(s_i))|_p$, le d\'eterminant
\'etant calcul\'e relativement \`a une base, univolumique pour
$d_{\g(s)}Y$ (resp. $d_{\g(s_i)}Y$), de $\g(s)$ (resp. $\g(s_i)$).

\noindent Comme $x_i^{-1}$ induit un isomorphisme de vari\'et\'es symplectique de
$x_i.[\O_{h,\lambda_\tau, s_i}]_{s_i}$ sur
$[\O_{h,\lambda_\tau, {s_i}}]_{s_i}$, il s'en suit que
\begin{eqnarray*}
 && \int_{[\O_{h,\lambda_\tau, s_i}]_{s_i}}
(1_{\g(s_i)_{\e_G(s)|n_F|_p}}\var^{x} (s_i\exp_{G(s_i)}(\, . \,
))d_{\g(s_i)}X)\, \widehat{}_{\g(s_i)} (l)\phi_{g, \tau, s_i}(l)
d\mu_{[\O_{h,\lambda_\tau, s_i}]_{s_i}}(l)= \\&& \quad  
c_i\int_{x_i.[\O_{h,\lambda_\tau,s_i}]_{s_i}}
(1_{\g(s)_{\e_G(s)|n_F|_p}}\var^{xx_i^{-1}} (s\exp_{G(s)}(\, . \,
))d_{\g(s)}X)\, \widehat{}_{\g(s)}(l) \phi_{g, \tau, s}(l)
d\mu_{\O_{g,\lambda_\tau, s}}(l)
\end{eqnarray*}
Cependant, pour tout $h \in L^1(G ; G(s))$, on a
$$\int_{G/G(s_i)}h(xx_i^{-1})d_{G/G(s_i)}\dot x =
c_i'\int_{G/G(s)}h(x)d_{G/G(s)}\dot x, $$ avec $c_i'= |\det(x_i^{-1} :
\g_s \longrightarrow \g_{s_i})|_p$, le d\'eterminant \'etant
calcul\'e relativement \`a une base de $\g_s$ (resp. $\g_{s_i}$),
univolumique pour  la mesure de Haar tangente \`a la mesure
quotient sur $G/G(s)$ (resp. $G/G(s_i)$). Il est facile de montrer
que $c_ic_i' =\Delta_G(x_i)$.
 Il en r\'esulte que
\begin{eqnarray*}
  I_{s_i}^x(\var d_G x) &=& |\det(1-s^{-1})_{\g_s}|_p \times \\ && \times
\int_{G/G(s)} \int_{x_i.[\O_{h,\lambda_\tau,s_i}]_{s_i}}
(1_{\g(s)_{\e_G(s)|n_F|_p}}\var^{xx_i^{-1}} (s\exp_{G(s)}(\, . \,
))d_{\g(s)}X)\, \widehat{}_{\g(s)} (l) \, \,  \times \\ &&
\times \, \phi_{g, \tau, s}(l) d\mu_{\O_{g,\lambda_\tau, s}}(l)
d_{G/G(s)}\dot x
\end{eqnarray*}
Compte tenu  de ce qui pr\'ec\`ede et du lemme \ref{yRicha},
on obtient
\begin{eqnarray*}
&&\Theta_{g, \tau}(\var d_Gx) =\int_{G/G(s)}
\Delta_G(x)^{-1}\left\{\int_{\O_{g,\lambda_\tau,s}}
(1_{\g(s)_{\e_G(s)|n_F|_p} }\var^{x} (s\exp_{G(s)}(\, . \,
))d_{\g(s)}X)\, \widehat{}_{\g(s)} (l) \right.
\\&& \qquad \qquad \qquad  \times\left. \phi_{g, \tau, s}(l)
d\mu_{\O_{g,\lambda_\tau, s}}(l)\right\}|\det(1-s^{-1})_{\g_s}|_p\
d_{G/G(s)}\dot x.
\end{eqnarray*}
D'o\`u la formule d\'esir\'ee.
\end{subsubsection}
\end{subsection}
\end{section}
\begin{center}{\bf  Mesure de Plancherel.}
\end{center}

\begin{section}{}\label{a} On
sait, depuis les travaux de (\cite{[G]}, \cite{[Di2]}), que si $G$ est unimodulaire alors il
existe  une mesure $d_{\widehat{G}}\mu$ sur $\widehat{G}$, appel\'ee mesure de
Plancherel de $G$, qui est enti\`erement d\'etermin\'ee une fois
fix\'ee une mesure de Haar $d_Gx$ sur $G$ par~:
 $$ \var(1) =
\int_{\widehat{G}} \mbox{Tr}(\pi(\var d_Gx))d_{\widehat{G}}\mu(\pi), \, \mbox{ pour tout }
\, \var \in C_c^\infty(G) .$$
 L'objet  de cette  partie du
pr\'esent travail est de  d\'ecrire explicitement la mesure $d_{\widehat{G}}\mu$.
\end{section}
\begin{section}{Formes lin\'eaires fortement r\'eguli\`eres.} Dans ce paragraphe nous allons \'etudier les formes lin\'eaires fortement r\'eguli\`eres de $\g^*$.
\begin{subsection}{}\label{regular}
On d\'esigne par $k'$ soit le  corps   $k$ soit le corps    $\bar k$.
 On pose $\g_{k'} =  {k'}\otimes_k\g$.
Un \'el\'ement $g$ de $\g_{k'}^*$ est dit r\'egulier si 
$\dim\g_{k'}(g) \leq \dim\g_{k'}(g')$ pour tout $g' \in  \g_{k'}^*$. Si tel est le cas, $\g_{k'}(g)$ est une
sous-alg\`ebre de Lie alg\'ebrique commutative (\cite{[Di]},
Proposition~1.11.7). On note alors $\j_{ g}$ l'ensemble des
\'el\'ements semi-simple de $\g_{k'}(g)$.
 Ainsi,  $\j_{ g}$  est l'unique facteur
r\'eductif de $\g_{k'}(g)$. La forme lin\'eaire $g$ est dite
fortement r\'eguli\`ere si $\j_{ g}$ est de dimension maximale. Dans
ce cas, la sous-alg\`ebre $\j_{ g}$ est appel\'ee sous-alg\`ebre de
Cartan-Duflo. L'ensemble des \'el\'ements fortement
r\'eguliers de $\g_{k'}^*$ est not\'e $\g^*_{k', t.r}$. Lorsque $k'= k$, l'ensemble $\g^*_{k', t.r}$ est not\'e simplement $\g^*_{ t.r}$.
\begin{subsubsection}{}  \label{reg1} Soit $\j_{ \bar k}$ une sous-alg\`ebre de
Cartan-Duflo de $\g_{ \bar k}$. On note $ \h_{ \bar k}$ le
centralisateur de $\j_{ \bar k}$ dans $\g_{ \bar k}$. 
 Puisque $\j_{ \bar k}$ est une alg\`ebre de Lie r\'eductive,
 on a   $\g_{ \bar k}= \h_{ \bar k} \oplus [\j_{ \bar k}, \g_{ \bar k}] $
  ce qui permet d'identifier $\g_{ \bar k}^*$ \`a  $ \h_{ \bar k}^* \oplus[\j_{ \bar k}, \g_{ \bar k}]^* $.
\begin{lem} \label{yreg1} Soit $ l \in \h_{  \bar k}^*$ r\'egulier comme \'el\'ement
de  $\g^*_{{  \bar k}}$. On a $\g_{{  \bar k}}(l) \subset \h_{  \bar k}$ et
$\j_{l} = \j_{ \bar k}$, ce qui donne que $ l \in\g^*_{{  \bar k}, t.r}$.
\end{lem}
\begin{dem} Puisque $[\j_{ \bar k},\h_{  \bar k} ]= 0$, $[\j_{ \bar k},[\j_{  \bar k}, \g_{  \bar k}] ]
 \subset [\j_{  \bar k}, \g_{ \bar k}]$,  et le fait que $l$ est nulle
 sur $[\j_{  \bar k}, \g_{  \bar k}]$, on a $\j_{ \bar k} \subset \g_{{
  \bar k}}(l)$. Comme $l$ est un  \'el\'ement r\'egulier de  $\g^*_{{  \bar k}}$, on a
 $\g_{{  \bar k}}(l) \subset \h_{  \bar k} $ et $  \j_{ \bar k} \subset  \j_{l}$. A
 cause de la maximalit\'e de $\j_{ \bar k}$, on a $  \j_{ \bar k} =
 \j_{l}$.
\end{dem}
\end{subsubsection}
\begin{subsubsection}{} \label{reg2}
 Soit  $B= (e_1, \ldots, e_{2d})$ une base  de $ [\j_{ \bar k}, \g_{ \bar k}]$ et  $(e_1^*, \ldots, e_{2d}^*)$ la base duale.  Si $l \in
 \h_{ \bar k}^*$, on note $\pi_B(l)$ le pfaffien de la restriction de la
 forme $\beta_l$ \`a $[\j_{ \bar k}, \g_{ \bar k}]$ calcul\'e dans la base $B$~: $\frac{1}{d \, !}\, {\beta_l}_{|[\j_{ \bar k}, \g_{ \bar k}]}^{\, d} = \pi_B(l) e_1^*{\wedge} \ldots \wedge  e_{2d}^*$.
 \begin{lem}\label{yreg2} Soit $f \in \h_{ \bar k}^*$.
Les assertions suivantes sont \'equivalantes :
\begin{enumerate}
\item $f \in \g^*_{{  \bar k},t.r}$,
\item $f$ est r\'egulier dans $\h_{ \bar k}^*$  et $\pi_B(f) \neq 0$.
\end{enumerate}
   \end{lem}
\begin{dem}
Supposons que $f \in \g^*_{{  \bar k},t.r}$. D'apr\`es le lemme
\ref{yreg1}, on a $\g_{ \bar k}(f) \subset  \h_{ \bar k}$, $\g_{ \bar k}(f) =
\h_{ \bar k}(f)$,  et $  \j_{ \bar k} = \j_{f}$, ce qui prouve que la
restriction de $\beta_f$ \`a
 $[\j_{ \bar k}, \g_{ \bar k}] $ est non
 d\'eg\'en\'er\'ee et par suite  $\pi_B(f) \neq 0$.
 Soit maintenant
 $l$ un \'el\'ement r\'egulier de $\h_{ \bar k}^*$ tel que  $\pi_B(l) \neq 0$.
 Comme $[\h_{ \bar k},[\j_{ \bar k}, \g_{ \bar k}] ] \subset [\j_{ \bar k}, \g_{ \bar k}]$, on a
 $\h_{ \bar k}(l) \subset \g_{ \bar k}(l)$. Soit $X \in \g_{ \bar k}(l)$,
 \'ecrivons $X= X_1 + X_2$, $X_1 \in \h_{ \bar k}$ et $X_2 \in  [\j_{ \bar k},
 \g_{ \bar k}]$. Pour tout $Y \in  [\j_{ \bar k}, \g_{ \bar k}]$,  $<l, [X,Y ]> =
 0$ et ceci implique que $<l, [X_2,Y ]> = 0$. Comme la restriction
  de $\beta_l$ \`a  $[\j_{ \bar k}, \g_{ \bar k}] $
  est non d\'eg\'en\'er\'ee, $X_2=0$. Ainsi $\g_{ \bar k}(l) \subset
  \h_{ \bar k}$ et par cons\'equent  $  \g_{ \bar k}(l) \subset \h_{ \bar k}(l)$. Si bien que l'on a  $  \g_{ \bar k}(l) = \h_{ \bar k}(l)$.
D'autre part, on a $$\dim \h_{ \bar k}(l)\leq \dim \h_{ \bar k}(f) = \dim \g_{ \bar k}(f) \leq
\dim \g_{ \bar k}(l) .$$ 
Il en r\'esulte que $$\dim \h_{ \bar k}(l) =  \dim \h_{ \bar k}(f) = \dim \g_{ \bar k}(f) =
\dim \g_{ \bar k}(l) .$$
Ce qui prouve le r\'esultat.
\end{dem}

\medskip
On pose $ \goth H^*_{{\bar k}, t.r} = \h_{ \bar k}^* \cap \g^*_{{\bar k},
 t.r}$.  Comme l'ensemble des \'el\'ements r\'eguliers de $\h_{ \bar k}^*$ est un ouvert de Zariski non vide de  $\h_{ \bar k}^*$, on d\'eduit  du lemme ci-dessus que  $ \goth H^*_{{\bar k}, t.r}$
est aussi  un ouvert de Zariski  non vide  de $ \h_{\bar k}^*$.
\end{subsubsection}
\begin{subsubsection}{}  \label{reg3} Le r\'esultat suivant a un r\^ole
important pour la suite de ce travail.
\begin{pro} \label{yreg3}
Avec les notations ci-dessus, toutes les sous-alg\`ebres de
Cartan-Duflo de $\g$ sont $G$-conjugu\'ees.
\end{pro}
\begin{dem}
L'application $\sigma : [\j_{\bar k}, \g_{\bar k}] \times \h_{\bar
k}^* \longrightarrow \g_{\bar k}^*$, $(X, f)\longmapsto\exp(X).f$, est \'evidemment polynomiale. Sa diff\'erentielle en
$(0, l_0), \, l_0 \in  \goth H^*_{{\bar k}, t.r}$ est
l'application~: $[\j_{\bar k}, \g_{\bar k}] \times \h_{\bar k}^*
\longrightarrow \g_{\bar k}^*$, $(X, l)\longmapsto X.l_0 + l$.
Supposons que l'on ait, pour tout $Y \in \g_{\bar k}$, $<X.l_0 +
l\, , \,  Y> = 0$. En prenant $Y \in \h_{\bar k}$, on trouve que
$l= 0$. Il en r\'esulte que $X \in \g_{\bar k}(l_0)$. Par le lemme
\ref{yreg1}, $X \in \h_{\bar k} \cap [\j_{\bar k}, \g_{\bar k}]$,
soit $X= 0$. La diff\'erentielle est donc injective, et surjective
par argument de dimension. Ainsi, $\sigma$ est \'etale en $(0,
l_0)$. Si bien que l'image de $[\j_{\bar k}, \g_{\bar k}] \times
\goth H^*_{{\bar k},t.r}$ par $\sigma$ contient un ouvert de
Zariski  non vide $\mathcal P$ de $\g_{\bar k}^*$ form\'e d'\'el\'ements fortement
r\'eguliers. On en d\'eduit qu'il existe $g_0 \in \g^*$ tel que, notant $\tilde g_0$ l'\'el\'ement de $\g_{\bar k}^*$ qui s'en d\'eduit par   extension  des scalaires \`a $\bar k$, on ait  $g_0 \in \mathcal P$. On a alors $\j_{\tilde g_0} = \bar k \otimes_k \j_{g_0}$.
 
On se donne maintenant $f \in \g^*_{\bar k,t.r}$. Soit $\h_{ \bar
k}'$ le centralisateur de $\j_f$ dans $\g_{\bar k}$ et ${{\goth
H}'}^*_{{\bar k},t.r} =  {\h'_{ \bar k}}^* \cap \g^*_{{ \bar k},
 t.r}$.
 Le raisonnement
ci-dessus affirme l'existence des \'el\'ements $X_1 \in [\j_{\bar
k}, \g_{\bar k}]$, $l_1 \in \goth H^*_{{\bar k},t.r}$,
 $X_2 \in [\j_{f}, \g_{\bar k}]$, $l_2 \in {{\goth
H}'}^*_{{\bar k},t.r}$  tels que $$\exp(X_1).l_1 =
\exp(X_2).l_2.$$ D'apr\`es le lemme \ref{yreg1}, on a $\j_{l_1} =
\j_{\bar k}$ et $\j_{l_2} =\j_f$. On en d\'eduit que $\j_{\bar k}$
et $\j_f$ sont conjugu\'es par un \'el\'ement de ${}^u{\bf G}$. On
vient de montrer que l'ensemble $\bf V$ des sous-alg\`ebres de
Cartan-Duflo de $\g_{\bar k}$ est un espace homog\`ene, d\'efini sur $k$,  de
${}^u{\bf G}$. On note ${\bf H}$ le centralisateur de $\j_{\tilde g_0}$ dans $\bf G$. 
 L'application canonique ${}^u{\bf G}/ {}^u{\bf H}
\longrightarrow {\bf V}$ induit une bijection ${}^u{
G}$-\'equivariante de $({}^u{\bf G}/ {}^u{\bf H})_k$ sur $ {\bf
V}_k$. En utilisant (\cite{[Bo-Se]}, Proposition 1.12), les
orbites de ${}^u{ G}$ dans $({}^u{\bf G}/ {}^u{\bf H})_k$
correspondent bijectivement aux \'el\'ements du noyau de
l'application canonique $H^1(k,{}^u{\bf H })\longrightarrow
H^1(k,{}^u{\bf G})$. Or,   $H^1(k,{ \bf U})=0$, pour tout
 groupe unipotent $ \bf U$ d\'efini sur $k$. Il s'ensuit que  les \'el\'ements de $ {\bf
V}_k$ constituent une seule orbite sous l'action de ${}^u G$.
 \end{dem}
Arriv\'e \`a ce stade, nous allons d\'emontrer le r\'esultat suivant. 
\begin{lem}  \label{yyreg3}  
   L'ensemble $\g^*_{t.r}$ est un ouvert de Zariski, non vide,  $G$-invariant de $\g^*$.
\end{lem}   
\begin{dem} On suppose que $\j_{\bar k}$ est d\'efini sur $k$. On note $\bf H$ le centralisateur de $\j_{\bar k}$ dans $\G$. C'est un sous-groupe ferm\'e, d\'efini sur $k$, de $\G$. 
 Consid\'erons le morphisme $\kappa : \G\times_{\bf H}  \goth H^*_{{\bar k},t.r} \longrightarrow \g_{\bar k}^*$, qui \`a $[x, f]$ associe
$x.f$. Puisque $\kappa$ est injectif et dominant, il est de 
type fini. Comme $\kappa$ est \'etale, d'apr\`es (\cite{[Mu]}, III.10. Th\'eor\`eme  3) le morphisme  $\kappa$ est plat.   D'apr\`es (\cite{[Mi]}, Th\'eor\`eme 2.12) le morphisme $\kappa$ est ouvert. Ainsi, $\g^*_{\bar k,t.r}$, qui est l'image de  $\kappa$, est un ouvert de $\g_{\bar k}^*$.  Il est d\'efini sur $k$ puisque $\kappa$ et $\goth H^*_{{\bar k},t.r}$ sont d\'efinis sur $k$. 
Si bien que $\g^*_{t.r}$, qui est l'ensemble des points rationnels de  $\g^*_{\bar k,t.r}$, est un ouvert de Zariski, non vide, de $\g^*$.  Il est clair qu'il est $G$-invariant.
   \end{dem}
\end{subsubsection}
\end{subsection}
 \begin{subsection}{}\label{V}
 On suppose d\'esormais que $G$ est unimodulaire. 
\begin{pro}\label{yV}
  Il existe un ouvert de Zariski  $\Omega_G$, non vide, $G$-invariant
  de $\g^*$ contenu dans
   $\g^*_{t.r}$
  tel que,
 pour tout $f \in \Omega_G$, l'orbite co-adjointe  $G.f$ soit ferm\'ee dans~$\g^*$.
\end{pro}
\begin{dem}
 D'apr\`es
\cite{[CD]}, il existe une
 partie $\G$-invariante $\bf U$ de $ \g_{\bar k}^*$ ouverte et dense dans
  $\g_{\bar k}^*$, telle que toute $\G$-orbite contenue dans $\bf U$
  soit ferm\'ee dans $\g_{\bar k}^*$. On sait qu'il existe $f \in \bar k[
\g_{\bar k}^*]$, non nul, tel que  $D_{f} = \{X \in \g_{\bar k}^*
\, / \, f(X) \neq 0 \}$ soit contenu dans $\bf U$. Soit $(v_1, \ldots, v_m)$ une base de $\g$. Alors $(v_1, \ldots, v_m)$ est une base de $\g_{\bar k}$. D\'esignons par $(v_1^*, \ldots, v_m^*)$ la base duale dans $\g_{\bar k}^*$. Relativement \`a cette base, on a $f \in \bar k[x_1, \ldots, x_m]$.
 Soit $k_1$ une extension galoisienne  finie  de $k$ contenant  les coefficients de $f$.
 Notons ${\goth \Gamma}$ le groupe de Galois de
$ k_1$ sur $k$. Posons $\tilde f =
\prod_{\sigma \in \goth \Gamma } f^\sigma$.  Il est clair que $ \tilde f \neq 0$ et que  $\tilde f \in  k[x_1, \ldots, x_m]$.  Si bien que  $D_{\tilde f}$ est un ouvert de Zariski
(non vide) de $\g_{\bar k}^*$,  d\'efini sur $k$, et contenu dans
$\bf U$.  Il en r\'esulte que  $ \cup_{x \in {\bf G}} \, x.D_{\tilde f}$ est un
ouvert de Zariski,  d\'efini sur $k$, ${\G}$-invariant, et il est
contenu dans $\bf U$.
   La proposition r\'esulte alors  du lemme suivant.   
\end{dem}

\begin{lem} \label{yyV} Soit $g \in \g^*$. On note $\tilde g$ l'\'el\'ement de $ \g_{\bar k}^*$ obtenu \`a partir de $g$ par extension des scalaires \`a ${\bar k}$.
 Si ${\G}.\tilde g$
est ferm\'ee (de Zariski) dans $ \g_{\bar k}^*$ alors $G.g$ est 
ferm\'ee dans ${\bf \g}^*$ pour la topologie $p$-adique.
\end{lem}

\begin{dem} Puisque la sous-vari\'et\'e ${\G}.\tilde g$ est lisse, on d\'emontre par une m\'ethode analogue \`a celle utilis\'ee dans (\cite{[Wh]}, Th\'eor\`eme 1) que 
le ferm\'e $({\G}. \tilde g)_k$ de $\g^* $ est une sous-vari\'et\'e  
$k$-analytique de $\g^*$. Pour tout $f \in ({\G}. \tilde g)_k$,
l'application de ${\G}_k$ dans $({\G}. \tilde g)_k$ qui \`a $x$ associe
$Ad^*x.f$ est une submersion. Il s'en suit que l'orbite de chaque
\'el\'ement de $({\G}.\tilde  g)_k$ sous l'action de $ {\G}_k$ est ouverte
dans $({\G}.\tilde  g)_k$. Comme l'orbite $ {\G}_k.g$ est le
compl\'ementaire d'une r\'eunion de $ {\G}_k$-orbites dans  $({\G}.\tilde  g)_k$, elle est
ferm\'ee dans $\g^*$. D'autre part, pour
tout $f \in {\G}_k.g$,  l'application de $G$ dans $
{\G}_k.g$ qui \`a $x$ associe $Ad^*x.f$ est une submersion. Il en r\'esulte que l'orbite $G.g$ de $g$ sous l'action co-adjointe  de
$G$ est ferm\'ee dans $\g^*$.
\end{dem}
\end{subsection}
\end{section}
\begin{section}{D\'esint\'egration d'une mesure.}\label{adem Med}
\begin{subsection}{} Posons $U = {}^uG$ et $\u = {}^u\g$. Fixons $g_0 \in \Omega_G$.  D\'esignons par  $u_0$ sa restriction  \`a $\u$. 
 Notons $\j$ la sous-alg\`ebre de Cartan-Duflo de $\g(g_0)$.
Fixons un facteur  r\'eductif $\t$ de $\g$ contenant $\j$. Consid\'e\-rons la d\'ecomposition $\g = \t \oplus \u$ et la d\'ecomposition duale $\g^* = \t^* \oplus \u^*$. Alors, pour tout $g \in \Omega_G$ et pour tout $t \in \t^*$, on a : 
\begin{itemize}
	\item l'orbite co-adjointe de 
$t+ g $ sous l'action de $G$ est ferm\'ee dans $\g^*$ ;
\item $\g(t+ g) = \g(g)$. 
\end{itemize}
  On peut supposer d\'esormais que 
$$\Omega_G = \t^* +\Omega_G.$$
On note $\Omega_U$ l'image de $\Omega_G$ par l'application restriction de $\g^*$ \`a $\u^*$. Alors $\Omega_U$ 
est un ouvert de Zariski  non vide, $G$-invariant, de $\u^*$ et on a : 
$$\Omega_G = \t^* +\Omega_U.$$
On d\'esigne par $\Omega_{G, u}$ les formes lin\'eaires  de type unipotent  de $\Omega_G$. On a alors 
$$\Omega_G = \t^* +\Omega_{G, u}.$$

 \begin{lem}  La sous-alg\`ebre 
$\j$ est un facteur r\'eductif de $\g(u_0)$ et  tout facteur r\'eductif de $G(g_0)$  est un facteur r\'eductif de $G(u_0)$. 
\end{lem}
\begin{dem}
Soit $\r_{u_0}$ un facteur r\'eductif de $\g(u_0)$ contenant $\j$ et $\r$ un facteur r\'eductif de $\g$ contenant $\r_{u_0}$. Soit $X \in \r_{u_0}$. Si $Y \in \r$ et $Z \in \u$,  on a : $[X, Y]= 0$ et $<u_0,[X,  Z] >= 0$. Donc 
$$
  <g_0,[X, Y + Z] > = <g_0,[X, Z] > = <u_0,[X,  Z] >= 0.
  $$
Si bien que $X \in \j$. 

\noindent  Soit $R_{g_0}$ un facteur r\'eductif de $G(g_0)$ et  $R_{u_0}$ un facteur r\'eductif de $G(u_0)$ contenant $R_{g_0}$. 
Soit $x \in R_{u_0}$. On a : 
\begin{eqnarray*}
  <x. g_0, Y + Z > = <g_0,x^{-1}Y + x^{-1}Z > &= & <g_0,Y> +  <u_0, x^{-1}Z >\\
  &=&<g_0,Y> +  <u_0, Z >\\
  &=& <g_0, Y + Z >,
  \end{eqnarray*}
pour tout $Y \in \r$ et $Z \in \u$. Donc $x \in  R_{g_0}$. 
\end{dem}
\begin{pro}
L'orbite $G.u_0$ admet une mesure positive $G$-invariante. 
\end{pro}
\begin{dem}
Il suffit de d\'emontrer que le sous-groupe $G(u_0)$ est unimodulaire. Posons $\b = \j + \u$ et $b$ la restriction de $g_0$ \`a $\b$. On voit $\b$ est un id\'eal de $\g$ et que $\b(b) = \g(u_0)$. 

D'autre part, l'application $\g/\b \times \b(b)/\g(g_0)\longrightarrow k, (\dot{X}, \dot{Y}) \longmapsto \beta_{g_0}(X, Y)$ est non d\'eg\'en\'er\'ee, $G(g_0)$-invariante.  Comme l'action de $R_{g_0}$ dans $\g/\b$ est unimodulaire donc 
l'action de $R_{g_0}$ dans $\b(b)/\g(g_0)$ est unimodulaire. Etant donn\'e que $G$ est unimodulaire et l'action de $R_{g_0}$ dans $\g/\g(g_0)$ est unimodulaire, l'action de $R_{g_0}$ dans $\g(g_0)$ est unimodulaire. Il en r\'esulte que l'action de $R_{g_0}$ dans $\g(u_0)$ est unimodulaire.
\end{dem}
\end{subsection}
\begin{subsection}{}\label{integrti}
Fixons des mesures de Haar $d_{\g}X$, $d_{\u}X$, et $d_{\j}X$  sur $\g$, $\u$,  et $\j$ respectivement. D\'esignons par 
$d_Gx$ la mesure de Haar  sur $G$ tangente \`a $d_{\g}X$ et par 
 $d_{B_{u_0}}x$ la mesure de Haar sur $B_{u_0}= G(u_0)U$ tangente \`a la mesure de Haar $d_{\b}X = d_{\j}Xd_{\u}Y$ sur $\b$. 
 La mesure de Haar quotient sur $G/B_{u_0}$ est not\'ee $d_{G/B_{u_0}}\dot x$.   D\'esignons par   $d\mu_{U.u_0}$ la mesure de Liouville sur l'orbite co-adjointe $U.u_0$.
\begin{lem} Si $\var $ est une fonction bor\'elienne positive sur $G.u_0$, on pose  
 $I_{u_0}(\var)=  \int_{G/B_{u_0}}\int_{U.u_0}\var(xl)d\mu_{U.u_0}(l)d_{G/B_{u_0}}\dot x$. Alors, $I_{u_0}$ 
d\'efinit une mesure positive $G$-invariante sur $G.u_0$. On la notera $\beta_{G.u_0}$.

\end{lem}
Maintenant, nous allons d\'emontrer le r\'esultat suivant. 
\begin{pro}\label{yyintegrti} On suppose que $g_0$ est de type unipotent et que $\beta_{G.u_0}$ est une mesure de Radon sur $\u^*$.
Soit $\var \in C_c^\infty(\g)$. On a : 
$$
 \int_{G.u_0} \widehat{(\var_{|\u} d_\u X)}_{\u}(l)d\beta_{G.u_0}(l)=  \int_{\j^* } \int_{\O_{g_0, \lambda}}\widehat{(\var d_\g X)}_{\g}(l) d\mu_{\O_{g_0, \lambda}} (l)d_{\j^*}\lambda.
$$
\end{pro}
\begin{dem} 
On a : 
$$
 \widehat{(\var_{|\u} d_\u X)}_{\u}(l)=  \int_{\b^\perp \times \j^*}\widehat{(\var d_\g X)}_{\g}(\tilde l+ \lambda_1 + \lambda_2) d_{\b^\perp}\lambda_2d_{\j^*}\lambda_1,
 $$
 o\`u $\tilde l$ est la forme lin\'eaire sur $\g$ nulle sur $\t$ et dont la restriction \`a $\u$ est $l$. 
 Donc 
 \begin{eqnarray*}
 && \int_{G.u_0} \widehat{(\var_{|\u} d_\u X)}_{\u}(l)d\beta_{G.u_0}(l) \\
  & &\quad =  \int_{G/B_{u_0}}\int_{U/U(u_0)} \int_{\b^\perp \times \j^*}    \widehat{(\var d_\g X)}_{\g}(xy.\tilde u_0+ \lambda_1 + \lambda_2)d_{\b^\perp}\lambda_2d_{\j^*}\lambda_1 d_{U/U(u_0)}\dot yd_{G/B_{u_0}}\dot x\\
  && \quad = \int_{ \j^*}\int_{G/B_{u_0}}\int_{U/U(u_0)} \int_{\b^\perp}    \widehat{(\var d_\g X)}_{\g}(xy.\tilde u_0 + \lambda_1+ \lambda_2)d_{\b^\perp}\lambda_2d_{U/U(u_0)}\dot yd_{G/B_{u_0}}\dot xd_{\j^*}\lambda_1.
 \end{eqnarray*}
 Cependant, $B_{u_0}(b) = R_{g_0}U(u_0)$,  il s'ensuit que 
  \begin{eqnarray*}
 &&\int_{G/B_{u_0}}\int_{U/U(u_0)} \int_{\b^\perp}    \widehat{(\var d_\g X)}_{\g}(xy.\tilde u_0 + \lambda_1+ \lambda_2)d_{\b^\perp}\lambda_2d_{U/U(u_0)}\dot yd_{G/B_{u_0}}\dot x\\
 &&\qquad \qquad\qquad\qquad\qquad=\int_{\O_{g_0, \lambda_1}}\widehat{(\var d_\g X)}_{\g}(l) d\mu_{\O_{g_0, \lambda_1}} (l).
  \end{eqnarray*}
 Si bien que l'on a :  
  $$
  \int_{G.u_0} \widehat{(\var_{|\u} d_\u X)}_{\u}(l)d\beta_{G.u_0}(l) = \int_{\j^* } \int_{\O_{g_0, \lambda}}\widehat{(\var d_\g X)}_{\g}(l) d\mu_{\O_{g_0, \lambda}} (l)d_{\j^*}\lambda.
 $$
\end{dem}

\end{subsection}
\begin{subsection}{}
Le r\'esultat suivant a son propre int\'er\^et.
\begin{pro}
Soit $\var $ une fonction bor\'elienne positive sur $\u^*$. 
L'application $$u \in \Omega_U \longrightarrow I_{u}(\var)$$ est bor\'elienne.
\end{pro}
\begin{dem} D\'esignons par $\bf J$ le sous-groupe connexe de $\bf G$ d'alg\`ebre de Lie $\j$ et par $J$ l'image r\'eciproque de ${\bf J}_k$ dans $G$. On note $B = JU$. Alors $B$ est un    sous-groupe, alg\'ebrique,  normal de $G$, d'alg\`ebre de Lie $\j + \u$. On voit que $B$ est un sous-groupe de $B_{u_0}$ d'indice fini. On munit $B$ de la mesure de Haar tangente \`a $d_{\j}Xd_{\u}Y$.  Par un calcul \'el\'ementaire, on a  : 
$$I_{u}(\var) =\frac{1}{[B_u: B]}\int_{G/B}\int_{U.u}\var(xl)d\mu_{U.u}(l)d_{G/B}\dot x, \mbox{ pour tout }u \in \Omega_U. $$
Remarquons que $$\int_{U.u}\var(xl)d\mu_{U.u}(l)= \int_{U.xu}\var(l)d\mu_{U.xu}(l).$$
Il suit que l'application $(x, u) \in G \times \Omega_U \longmapsto \int_{U.u}\var(xl)d\mu_{U.u}(l)$ est bor\'elienne.
Si bien que l'application $u \in \Omega_U \longrightarrow \int_{G/B}\int_{U.u}\var(xl)d\mu_{U.u}(l)d_{G/B}\dot x$ est bor\'elienne.

Dans la suite, nous allons montrer que $u \in \Omega_U \longmapsto\frac{1}{[B_u: B]}$ est bor\'elienne. On d\'esigne par $\mathcal F(G)$ l'ensemble des sous-groupes ferm\'es de $G$ que l'on munit de la topologie de Fell \cite{[F]}.  L'application $u \in  \Omega_U \longmapsto G(u)$ est bor\'elienne. Donc $A: = \{(u, G(u)), \, u \in  \Omega_U \}$ est un bor\'elien de $ \Omega_U \times \mathcal F(G)$, o\`u $\Omega_U \times \mathcal F(G)$ est muni de la structure bor\'elienne produit.  Il suffit maintenant de d\'emontrer que $\alpha :(u, G(u)) \in  A \longmapsto [B_u: B] \in \N$ est  bor\'elienne. Nous allons prouver, en fait, que $\alpha$ est s\'equentiellement semi-continue inf\'erieurement. Soit $(u_n, G(u_n))$ une suite convergente dans $A$ vers $(u, G(u))$. On a : 
$$\alpha(u, G(u)) \leq \lim \inf_{n \rightarrow + \infty} \alpha(u_n, G(u_n)).$$
En effet, soit $m$ une valeur d'adh\'erence de la suite $(\alpha(u_n, G(u_n)))$. Il existe une sous-suite $(u_{\beta(n)}, G(u_{\beta(n)}))$ telle que $\lim_{n \rightarrow + \infty}\alpha(u_{\beta(n)}, G(u_{\beta(n)}))= m$. 
Soit $x^1 = e, \ldots, x^d$ ($d= \alpha(u, G(u))$) des \'el\'ements de $B_u= G(u)U$ tels que 
$$B_u = \bigsqcup_{i=1}^d x^i B.$$
Comme $U$ est un sous-groupe de $B$, on peut supposer que $x^i \in G(u)$. Il existe une suite $x^i_{\beta(n)} \in  G(u_{\beta(n)})$ convergente vers $x^i$, $1 \leq i \leq d$.
Soit $1\leq i \neq j \leq d$. On a $ x^i_{\beta(n)}(x^j_{\beta(n)})^{-1}\longrightarrow x^i(x^j)^{-1} \not \in B$. Comme $B$ est ferm\'e dans $G$, 
il suit que, pour $n$ assez grand, $  x^i_{\beta(n)}B \cap x^j_{\beta(n)}B = \emptyset$. Si bien que,  pour $n$ assez grand, $\alpha(u_{\beta(n)}, G(u_{\beta(n)})) \geq d$. D'o\`u $m\geq d$.  
\end{dem}
\end{subsection}
\begin{subsection}{}\label{MAISALCHKH}
On d\'esigne par $d_{\u^*}l$ la mesure de Haar  sur $\u^*$ duale de la mesure de Haar $d_{\u}X$ sur $\u$. 
On a le r\'esultat suivant. 
\begin{pro}\label{yMAISALCHKH}
Il existe une unique mesure bor\'elienne positive $\mu_u$ sur $G\backslash \u^*$ telle que l'on ait 
$$
 \int_{\u^*}\var(l) d_{\u^*}l = \int_{G\backslash \u^*}d\mu_u(\omega)\int_{\omega}\var(l)d\beta_\omega(l)
  $$
  pour toute fonction $\var$,  bor\'elienne positive, ou int\'egrable pour $d_{\u^*}l$ sur $\u^*$.
\end{pro} 

\begin{dem}
La d\'emonstration de cette proposition est de m\^eme style que celle donn\'ee par (\cite{[D-R]}, Lemme 5.1.7).\end{dem}

Soit $\a^*$ un r\'eseau dans $\u^*$. La proposition ci-dessus montre que pour presque tout $\omega \in  G\backslash \u^*$, on a : 
$$
 \int_{\omega}1_{\varpi^n\a^*}(l)d\beta_\omega(l) < + \infty, \mbox{ pour tout } n \in \Z.
  $$
Si bien que, pour presque tout $\omega \in  G\backslash \u^*$, la mesure $d\beta_\omega$ est une mesure de Radon.
\end{subsection}
\end{section}

 \begin{section}{La mesure de Plancherel de $G$.} \label{amin}
 \begin{subsection}{}\label{plancar}   Rappelons que $\Omega_{G, u}$ est l'ensemble des formes de type unipotent contenu dans $\Omega_G$.
Soit $g \in \Omega_{G, u}$.  D'apr\`es la proposition \ref{yreg3}, il existe $x \in G$ 
tel que   $\j_g = x.\j$ soit l'unique facteur r\'eductif de $\g(g)$.  
 On munit  $\j_g$ de la mesure de Haar $d_{\j_g}X$, image de la mesure de Haar
 $d_\j X$ sur $\j$ par  l'application $ \j \longrightarrow \j_g, X \longmapsto x.X$. 
 Soit $R_g$ un facteur r\'eductif de $G(g)$ (d'alg\`ebre de
Lie $\j_g$). 
On note $\widehat{(R_g^\g)}_{-}$ l'ensemble des classes  des
repr\'esentations unitaires irr\'eductibles $\tau$ de $ R_g^\g$ 
telle que $\tau (1,-1)= -\mbox{Id}$.  Alors $\widehat{(R_g^\g)}_{-}$ est un ouvert de $\widehat{R_g^\g}$ (voir \cite{[F2]}, Corollaire du Th\'eor\`eme 1.3). 
On d\'esigne par  $d_{ R_g^\g}x$  la
mesure de Haar sur  $ R_g^\g$  tangente \`a  $d_{\j_g}X$ 
   et par $d_{ \widehat{R_g^\g}}\tau$  la mesure de Plancherel de
$\widehat{R_g^\g}$  correspondante. 
Si $\tau \in \widehat{(R_g^\g)}_{-}$, on d\'esigne par $ \tau \otimes \chi_g$ l'\'el\'ement de $Y_G(g)$ d\'efini par 
$$\tau \otimes \chi_g(xy)= \chi_g(y)\tau(x), \, \, x \in R_g^\g, \, y \in {}^uG(g)^\goth g.$$
L'application   $ \widehat{(R_g^\g)}_{-}\longrightarrow  Y_G(g), \, \,  \tau \longmapsto \tau \otimes \chi_g$ est un hom\'eomorphisme : elle est continue bijective, l'application r\'eciproque \'etant  $ Y_G(g) \longrightarrow \widehat{(R_g^\g)}_{-} , \, \,  \tau \longmapsto \tau_{|R_g^\g}$ qui est continue d'apr\`es le corollaire du th\'eor\`eme 1.3 de \cite{[F2]}.  On note  alors $d_g\tau$ la mesure image de  $d_{ \widehat{R_g^\g}}\tau$ sur $Y_G(g)$. Remarquons que cette mesure ne d\'epend pas du choix de  $R_g$. 
\begin{subsubsection}{}\label{desmes}
 On d\'esigne par $Z_g$ le centre de $R_g^\g $ et par  $N= \{ (1,-1), \, (1,1)\}$.
 Alors $Z_g$ est  un sous-groupe ferm\'e,   d'indice fini, de $R_g^\g$ contenant $N$.
 Soit  $\chi_0$ le caract\`ere 
  non trivial de $N$ et 
  $\widehat{Z_g}_{-} = \{\chi \in \widehat{Z_g}, \chi_{|N} = \chi_0 \}$. Alors,  $\widehat{Z_g}_{-}$  est
ouvert et ferm\'e dans $\widehat{Z_g}$ (voir \cite{[We3]}, Chapitre 6). 
  On munit  $Z_g$
de  la mesure de Haar $d_{Z_g}z$  restriction de la mesure de
Haar $d_{ R_g^\g}x$ sur $ R_g^\g$
 et  $\widehat{Z_g}$ de la
mesure de Haar $d_{Z_g}^*\chi$ duale de $d_{Z_g}z$. 
Pour chaque
  $\chi \in \widehat{Z_g}_{-}$, on pose
 $$ (R_g^\g)\, \widehat{}_\chi  =\{\tau \in   \widehat{R_g^\g} \mbox{ tel que }
 \tau_{| Z_g} = (\dim \tau) \chi\}. $$
Si $I$ est un
ensemble fini, on note $|I|$ le nombre de ses \'el\'ements. D'apr\`es (\cite{[G-M]}, Th\'eor\`eme 4.1),
on a, pour tout $\beta \in C_c^\infty(R_g^\g)$,
\begin{eqnarray}\label{aloula}
 \int_{\widehat{(R_g^\g)}_{-}}  \mbox{ Tr}(\tau(\beta d_{ R_g^\g}x))d_{\widehat{R_g^\g}}\tau &=&
  \int_{\widehat{Z_g}_{-}} \sum_{\tau \in (R_g^\g)\, \widehat{}_\chi }\dim \tau
 \mbox{ Tr}(\tau(\beta d_{
  R_g^\g}x))
    \frac{1}{|R_g^\g/Z_g|} d_{Z_g}^*\chi.
   \end{eqnarray}
\end{subsubsection}
\begin{subsubsection}{} Fixons $0< \e<{\bf a}_G$.
 Alors,  $ { R_{g,\e}^\g}N$ est un sous-groupe ouvert compact de $Z_g$.
 On  pose $( {R_{g,\e}^\g}N)^\perp  = \{\chi \in  \widehat{Z_g}\, , \,
   \chi_{|  {R_{g,\e}^\g}N} = 1\}$. C'est un  sous-groupe compact de $
 \widehat{Z_g}$.
  On   munit  $( {R_{g,\e}^\g}N)^\perp$ de la mesure de Haar normalis\'ee $d_\e \chi$  et $\widehat{Z_g}/( {R_{g,\e}^\g}N)^\perp
  $ de la mesure de Haar quotient $d_{Z_g}^* \dot \chi =  d_{Z_g}^*  \chi/d_\e \chi$. On a : 
  $$d_{Z_g}^* \dot \chi = \frac{1}{mes({R_{g,\e}^\g}N)}\sum_{\chi  \in \widehat{Z_g}/( {R_{g,\e}^\g}N)^\perp}\delta_{\chi}. $$
   On munit \'egalement le sous-groupe  compact ${(\j_{g,\e|n_F|_p})}^\perp$ de
$\j_g^*$  de la mesure de Haar normalis\'ee et
$\j_g^*/{(\j_{g,\e|n_F|_p})}^\perp$  de la mesure de Haar quotient
$d_{\j_g}^*\bar{\lambda}$. On a : 
$$d_{\j_g}^*\bar{\lambda} = \frac{1}{mes(\j_{g,\e|n_F|_p})}\sum_{\bar \lambda \in \j_g^*/{(\j_{g,\e|n_F|_p})}^\perp}\delta_{\bar \lambda}.$$
Remarquons que 
$$mes({R_{g,\e}^\g}N) = 2\, mes(\j_{g,\e|n_F|_p}). $$ D'autre part,  chaque \'el\'ement $\bar {\lambda} \in
\j_g^*/{(\j_{g,\e|n_F|_p})}^\perp$ d\'efinit
   un caract\`ere $\psi_{\bar {\lambda}}$ 
   de ${ R_{g,\e}^\g}N$  en posant,
    $$\psi_{\bar {\lambda}} (1,-1) = -1 \mbox{ et }
    \psi_{\bar {\lambda}}(\exp_{ R_g^\g}(X)) = \varsigma(<\lambda,X>),
    \mbox{ pour tout } X \in \j_{g,\e|n_F|_p},$$
   $\lambda$ \'etant un rel\`evement de $\bar \lambda$ dans
   $\j_g^*$. Alors, l'application $\bar {\lambda}\longmapsto \psi_{\bar {\lambda}}$ est une bijection   de $\j_g^*/{(\j_{g,\e|n_F|_p})}^\perp$ sur $\widehat{({ R_{g,\e}^\g}N)}_{-}: = \{\chi \in \widehat{({ R_{g,\e}^\g}N)}, \, \chi (1,-1)= -1\}$.

\noindent
Il en r\'esulte que, pour toute fonction $\alpha$
mesurable positive ou int\'egrable sur $\widehat{Z_g}_{-}$, l'on~a~: 
\begin{eqnarray}\label{plan}
\int_{\widehat{Z_g}_{-}} \alpha(\chi)  d_{Z_g}^*\chi &=&\frac{1}{2}\int_{\j_g^*/{(\j_{g,\e|n_F|_p})}^\perp}
  \int_{({R_{g,\e}^\g}N)^\perp}\alpha(\tilde{\psi}_{\bar {\lambda}}.\chi) d_\e\chi d_{\j_g}^* \bar
  \lambda, 
\end{eqnarray}
$\tilde{\psi}_{\bar {\lambda}}$ \'etant un caract\`ere de ${Z_g}$
prolongeant $ \psi_{\bar {\lambda}}$.

\noindent Si bien que, pour toute fonction  $\alpha$ 
mesurable positive ou int\'egrable sur $\widehat{(R_g^\g)}_{-}$, l'on a~: 
\begin{eqnarray}\label{fields}
  \int_{\widehat{(R_g^\g)}_{-}} \alpha(\tau) d_{\widehat{R_g^\g}}\tau 
&=&
 \frac{1}{2}\int_{\j_g^*/{(\j_{g,\e|n_F|_p})}^\perp}
 \int_{({R_{g,\e}^\g}N)^\perp}
  \sum_{\tau \in (R_g^\g)\, \widehat{}_{\tilde \psi_{\bar {\lambda}}.\chi }}
  \frac{\dim \tau}{|R_g^\g/Z_g|} \alpha(\tau)\,
  d_\e\chi d_{\j_g}^* \bar \lambda .
 \end{eqnarray}
\end{subsubsection}
\end{subsection}
\begin{subsection}{} 
Pour $g \in \g^*$ de type unipotent et $\tau \in Y_G(g)$,   on rappelle que   $\pi_{g, \tau}$ est la classe de repr\'esentations unitaires irr\'eductibes de $G$, associ\'ee \`a $(g, \tau)$ par la   m\'ethode des orbites de
Kirillov-Duflo. Si $\pi_{g, \tau}$  est admissible 
alors $ \Theta_{g, \tau}$
   d\'esigne son caract\`ere, c'est la fonction g\'en\'eralis\'ee sur $G$ d\'efinie par la formule (\ref{anniver}).
\begin{pro}  Soit $g_0 \in \Omega_{G, u}$. On note $u_0$ sa restriction \`a $\u$. On suppose que $d\beta_{G.u_0}$ est une mesure de Radon.
Soit $\var \in C_c^\infty(G)$.
Alors, pour tout $g \in G.g_0$ , la fonction $\tau \longmapsto \Theta_{g, \tau}(\var d_G x)$, d\'efinie sur $Y_G(g) $,  est mesurable et  on a, 
$$\int_{Y_G(g) }  |\Theta_{g, \tau}(\var d_G x)|_{\cit}d_g\tau
 < + \infty.
$$
\end{pro}

\medskip

Dans ce qui suit nous allons   d\'emontrer  
cette proposition.
\begin{subsubsection}{}    Pour $m \in \N^\times$ et  $\chi \in \widehat{( {(R_g^\g)_\e}N)}_{-}$, on note
$${}_m\widehat{(R_g^\g)}_{\chi} = \{\tau \in
\widehat{(R_g^\g)}_{-}\, , \, \dim \tau = m \mbox{ et } \tau_{|R_{g,\e}^\g N } = \chi \}.$$ Alors
${}_m\widehat{(R_g^\g)}_{\chi}$  est une partie localement ferm\'ee de $\widehat{(R_g^\g)}_{-}$. Elle est
   localement compacte.

\noindent $\bullet$ On suppose que le support de $\var $ est contenu dans $G_\e$.
On a, pour tout $\tau \in {}_m\widehat{(R_g^\g)}_{\chi}$,
$$
 \Theta_{g, \tau}(\var d_Gx)= \dim \tau \int_{\O_{g, \lambda}}
    (\var \circ \exp d_\g X)\, \widehat{}_{\g}  (l  )
    d\mu_{\O_{g, \lambda}}(l),
 $$
o\`u $\lambda \in \j_g^*$ v\'erifiant $\chi(\exp_G(X)) = \varsigma(<\lambda, X>),\mbox{ pour tout } X \in \j_{g, \e|n_F|_p}$. Ainsi, la fonction $\tau \longmapsto \Theta_{g, \tau}(\var d_G x)$ est constante sur ${}_m\widehat{(R_g^\g)}_{\chi}$. Comme   $\widehat{( {(R_g^\g)_\e}N)}_{-}$ est au plus d\'enombrable, on en d\'eduit que la fonction $\tau \longmapsto \Theta_{g, \tau}(\var d_G x)$
  est mesurable sur $\widehat{ (R_g^\g)}_{-} $. 

\noindent $\bullet$ Soit $s\neq 1$ un \'el\'ement semi-simple de $G$  et
 $0<\e(s)< \min \{\e, a_G'(s),  a_G''(s)\}$ tel que l'application~:  $G
\times_{G(s)} G(s)_{\e(s)}\longrightarrow  \W_G(s , \e(s) )$, $[x, y]  \longmapsto 
  xs y x^{-1}$,  soit un diff\'eo\-morphisme. 
On suppose que  $ \var \in C_c^\infty(\W_G(s, \e(s)))$. On a, pour tout $\tau \in {}_m\widehat{(R_g^\g)}_{\chi}$, 
\begin{eqnarray*}
\Theta_{g, \tau }(\var d_Gx) 
 &=& |\det(1-s^{-1})_{\g/\g(s)}|_p
\sum_{\omega_s \in \O_{g,\lambda}\cap \g^*(s) }\phi_{g, \tau,
s}(\omega_s)
\\&& \int_{G/G(s)} \int_{\omega_s}(1_{\g(s)_{\e(s)|n_F|_p}}\var_s^x \circ \exp
d_{\g(s)}X)\, \widehat{}_{\g(s)} (l)
d\mu_{\omega_s}(l)d_{G/G(s)}\dot x.
 \end{eqnarray*}
 Cependant,  pour chaque $\tilde x \in R_g^\g$, la fonction $\tau \longmapsto \mbox{Tr}\,   \tau(\tilde x)$ est continue sur ${}_m\widehat{(R_g^\g)}_{\chi}$. Ainsi, 
 pour chaque $\omega_s \in \O_{g,\lambda}\cap \g^*(s)$, la fonction $\tau \longmapsto \phi_{g, \tau,
s}(\omega_s)$ est continue sur ${}_m\widehat{(R_g^\g)}_{\chi}$. 
Si bien que la fonction $\tau \longmapsto \Theta_{g, \tau}(\var d_G x)$
  est continue sur ${}_m\widehat{(R_g^\g)}_{\chi}$. On conclut que la fonction $\tau \longmapsto \Theta_{g, \tau}(\var d_G x)$
  est mesurable sur $\widehat{ (R_g^\g)}_{-} $.
 
 \noindent $\bullet$ Soit $S$ un syst\`eme de repr\'esentants des $G$-orbites semi-simples de $G$. 
 En utilisant  une partition de l'unit\'e subordonn\'ee \`a 
 $(\W_G(s , \e(s) ))_{s \in S}$, on voit que la fonction $\tau \longmapsto \Theta_{g, \tau}(\var d_G x)$ est mesurable sur $\widehat{ (R_g^\g)}_{-} $.
\end{subsubsection}
\begin{subsubsection}{}
 Soit $\a$ un r\'eseau de $\g$ contenu dans $\g_\e$ tel que 
 $K = \exp_G(\a) $ soit  un sous-groupe (compact ouvert) de $G$  et 
$\var$ soit $K$-bi-invariante.
On a~: 
$$\frac{1}{(mes(K))^2} 1_{K}*
\var * 1_{K}=\var, $$
o\`u on a not\'e, pour $\var_1$, $\var_2 \in C_c^\infty(G)$, \, $\var_1 * \var_2 (x)= \int_G \var_1 (y) \var_2 (y^{-1}x)d_Gy$, \, $x \in G$.
Soit $\mathcal H_{g,  \tau}$ une r\'ealisation  de $\pi_{g, 
\tau}$ et soit $v \in \mathcal H_{g,  \tau}$.
 On a~:
\begin{eqnarray*}
 < \pi_{g,  \tau}(\var d_Gx)v,v> &=&\frac{1}{(mes(K))^2} < \pi_{g,  \tau}(\var d_Gx) \pi_{g,  \tau}(1_K d_Gx)v, \pi_{g,  \tau}(1_K d_Gx)v>.
\end{eqnarray*}
En utilisant l'in\'egalit\'e de Cauchy-Schwarz ($|<\xi, \zeta>|_\cit \leq ||\xi||.||\zeta||, \, \, \xi, \zeta \in 
\mathcal H_{g,  \tau}$) et le fait que
$\pi_{g,  \tau}$ est unitaire, on a~:
\begin{eqnarray*}
|< \pi_{g,  \tau}(\var d_Gx)v,v> |_\cit &\leq& \frac{M_\var}{mes(K)}  <  \pi_{g,  \tau}(1_K d_Gx)v, \pi_{g,  \tau}(1_K d_Gx)v>
\\
&\leq& M_\var<  \pi_{g,  \tau}(1_K d_Gx)v, v>,
\end{eqnarray*}
o\`u  $M_\var =\frac{1}{mes(K)} \int_{G}|\var(x)|_\cit  d_Gx  $. 
En utilisant une base hilbertienne de $\mathcal H_{g,  \tau}$, on obtient,
 $$ |\Theta_{g, \tau}(\var d_Gx)|_\cit \leq
  M_\var \Theta_{g, \tau}(1_K d_Gx). $$ 
Soit $\lambda \in \j_g^*$ tel que $ \tau_{| {R_{g,\e}^\g}N} = (\dim \tau)\psi_{_{\bar {\lambda}}}$.
On a~: 
$$
 |\Theta_{g, \tau}(\var d_Gx)|_\cit  \leq M_\var  \dim \tau  
 \int_{\O_{g,\lambda}} (1_K \circ \exp_{G} d_\g
X)\, \widehat{}_{\g} (m) d\mu_{\O_{g,\lambda}}(m). 
$$
Par cons\'equent, en utilisant  la formule (\ref{fields}) puis  la proposition \ref{yyintegrti}, on a~:  
\begin{eqnarray*} 
\int_{ \widehat{ (R_g^\g)}_{-} }  |\Theta_{g, \tau}(\var d_G x)|_{\cit}d_{ \widehat{ R_g^\g}}\tau
& \leq &  \frac{1}{2} M_\var
\int_{\j_g^*} \int_{\O_{g,\lambda}} (1_K \circ \exp_{G} d_\g
X)\, \widehat{}_{\g} (m)d\mu_{\O_{g,\lambda}}(m)
d_{\j_g^*}\lambda \\
&\leq& \frac{1}{2} M_\var    \int_{G.u_0} (1_{K\cap U} \circ \exp d_\u
X)\, \widehat{}_{\u}(l) d\beta_{G.u_0} (l) <+ \infty.
\end{eqnarray*}
\end{subsubsection}
\begin{pro} \label{yAdam-Khemais} Soit $g_0 \in \Omega_{G, u}$ et $u_0$ sa restriction \`a $\u$.  On suppose que $d\beta_{G.u_0}$ est une mesure de Radon.
Alors, pour tout $g \in G.g_0$, on a :  
\begin{eqnarray} \label{MC1}
 \int_{Y_G(g)}  \Theta_{g, \tau}(\var d_Gx)\, d_g\tau&=& \frac{1}{2}\int_{G.u_0}\widehat{(\var_{|U}\circ \exp d_{\u}X)}_{\u}(l) d\beta_{G.u_0}(l) 
 \end{eqnarray}
 pour tout $\var \in C_c^\infty(G)$. 
 \end{pro}

\medskip

\noindent Gr\^ace \`a la proposition ci-dessus, on pose,  pour  $\omega \in  G \backslash  \Omega_U$, 
$$
\Theta_\omega(\var d_Gx) = 2 \int_{Y_G(g)}  \Theta_{g, \tau}(\var d_Gx)\, d_g\tau, 
 $$ 
 o\`u $g \in \Omega_{G, u}$ dont la restriction \`a $\u$ appartient \`a $\omega$.

\begin{cor} Soit $\var \in C_c^\infty(G)$.  La fonction 
$\omega \longmapsto \Theta_\omega(\var d_Gx) $
est $d\mu_{ u}$-int\'egrable.
\end{cor}
\end{subsection}

\begin{subsection}{}\label{principal}
Le r\'esultat principal de cette  partie est le suivant.
\begin{tho}\label{yprincipal}  Pour tout $\var
\in C_c^\infty(G)$, on a~:
\begin{eqnarray} \label{isra-ala}
 \var(1) &=& \int_{ G \backslash\u^*}\Theta_\omega(\var d_Gx)
d\mu_u(\omega ).
 \end{eqnarray}
\end{tho}
\begin{dem} 
 Le th\'eor\`eme \ref{yprincipal} r\'esulte  de la  formule
(\ref{MC1}) et de la proposition \ref{yMAISALCHKH}.
\end{dem}
\end{subsection}
 
 \end{section}
 \begin{section}{ D\'emonstration de la proposition \ref{yAdam-Khemais}}{}
 Soit $g \in \Omega_{G,u}$ telle que $d\beta_{G.(g_{|u})}$ soit une mesure de Radon. Soit $R_g$ un facteur r\'eductif de $G(g)$.
\begin{subsection}{}\label{finish} Dans ce num\'ero, nous allons d\'emontrer la proposition \ref{yAdam-Khemais} lorsque  $\var \in C_c^\infty(G_\e)$. 
 Soit  $\bar \lambda \in \j_g^*/{(\j_{g,\e|n_F|_p})}^\perp$, $\lambda \in
\j_g^*$  un rel\`evement de $\bar \lambda$,  et $\psi_{\bar
{\lambda}}$  le caract\`ere  de ${ R_{g,\e}^\g}N$~:
   $$\psi_{\bar {\lambda}} (1,-1) = -1 \mbox{ et }
    \psi_{\bar {\lambda}}(\exp_{ R_g^\g}(X)) = \varsigma(<\lambda,X>),
    \mbox{ pour tout } X \in \j_{g,\e|n_F|_p}.$$
    On d\'esigne aussi par $\psi_{\bar {\lambda}}$ un prolongement en un  caract\`ere
    de $Z_g$. Soit $\psi \in ({R_{g,\e}^\g}N)^\perp$.
Pour tout $\tau \in (R_g^\g)\,   \widehat{}_{\psi_{\bar {\lambda}}.\psi }$,
 on a, 
$$
 \Theta_{g, \tau}(\var d_Gx)= \dim \tau \int_{\O_{g, \lambda}}
    (\var \circ \exp d_\g X)\, \widehat{}_{\g}  (l )
    d\mu_{\O_{g, \lambda}}(l).
 $$
Comme
   $\sum_{\tau \in ({R_g}^\g)\,  \widehat{}_{\psi_{\bar {\lambda}}.\psi } } (\dim \tau)^2 = |{R_g}^\g/{Z_g}|,$
 on a~:
$$
  \sum_{\tau \in ({R_g}^\g)\,  \widehat{}_{\psi_{\bar {\lambda}}.\psi }}
  \frac{\dim \tau}{|R_g^\g/Z_g|} \Theta_{g, \tau}(\var d_Gx)= \int_{\O_{g, \lambda}}
    (\var \circ \exp d_\g X)\, \widehat{}_{\g}  (l  )
    d\mu_{\O_{g, \lambda}}(l).
   $$
Si bien que l'on a : 
\begin{eqnarray*}
  &&\int_{\j_g^*/{(\j_{g,\e|n_F|_p})}^\perp}\int_{({R_{g,\e}^\g}N)^\perp}
 \sum_{\tau \in (R_g^\g)\,   \widehat{}_{\psi_{\bar {\lambda}}.\psi }}
  \frac{\dim \tau}{|R_g^\g/Z_g|} \Theta_{g, \tau}(\var d_Gx)\,
  d_\e\psi d_{\j_g}^* \bar \lambda\\&&
 \qquad \qquad \qquad\qquad \qquad \qquad =
\int_{\j_g^*}\int_{\O_{g, \lambda}}
    (\var \circ \exp d_\g X)\, \widehat{}_{\g} (l )
    d\mu_{\O_{g, \lambda}}(l) d_{\j_g}^*\lambda.
 \end{eqnarray*}
La formule (\ref{MC1}) r\'esulte alors de ce qui pr\'ec\`ede, de la formule (\ref{fields}), et de la proposition~\ref{yyintegrti}.

\end{subsection}
\begin{subsection}{} \label{lem}
 Soit $s$ un \'el\'ement semi-simple de $G$,  $s \neq 1$, et
 $0<\e(s)< \min \{\e, a_G(s),  a_G''(s)\}$. 
\begin{subsubsection}{}\label{finiss}  Soit $ \var \in C_c^\infty(\W_G(s, \e(s)))$  et  $\bar \lambda \in \j_g^*/{(\j_{g,\e|n_F|_p})}^\perp$.
\begin{pro} \label{yyfiniss} On a~: 
 \begin{eqnarray}
 \sum_{\bar{\lambda'} \in
{(\j_{g,\e(s)|n_F|_p})}^\perp /{(\j_{g,\e|n_F|_p})}^\perp }
  \int_{({R_{g,\e}^\g}N)^\perp}
  \sum_{\tau \in (R_g^\g)\, \widehat{}_{\psi_{\bar {\lambda} + \bar{\lambda'}}.\psi }}
  \frac{\dim \tau}{|R_g^\g/Z_g|} \Theta_{g, \tau}(\var d_Gx)\,
  d_{\e}\psi  &=& 0.
\end{eqnarray}
\end{pro}
Dans ce qui suit, nous allons prouver ce r\'esultat.  Soit 
 $\bar{\lambda'} \in
{(\j_{g,\e(s)|n_F|_p})}^\perp /{(\j_{g,\e|n_F|_p})}^\perp$.  On note
 $ \bar \lambda  + \bar{\lambda'}= \bar\lambda_0  $ et $ \lambda_0 \in \j_g^*$ un repr\'esentant de $\bar\lambda_0$.

\noindent $-$ On suppose que $\O_{g, \lambda_0} \cap \g^*(s)= \emptyset$.  D'apr\`es
 le th\'eor\`eme \ref{yhd}, on a, pour tout  $\psi \in ({R_{g,\e}^\g}N)^\perp$ et pour tout  $\tau \in ({R_g}^\g)\, \widehat{}_{
\psi_{\bar \lambda_0}.\psi} $, 
 $$ \Theta_{g, \tau }(\var d_Gx) = 0.$$
$-$ On suppose  que
  $\O_{g, \lambda_0} \cap \g^*(s)\neq \emptyset$.  Comme $ \Theta_{g, \tau }$ est $G$-invariante,
   on peut supposer d\'esormais $s \in R_g$. On choisit   $\hat{s} \in R_g^\g$  un relev\'e de
  $s$. Soit  $\psi \in ({R_{g,\e}^\g}N)^\perp$. 
  En appliquant successivement  les th\'eor\`emes \ref{yPTo} et \ref{yhd} , on
  a, pour  tout $\tau \in ({R_g}^\g)\, \widehat{}_{(\psi_{\bar \lambda_0}.\psi)} $,
\begin{eqnarray*}
\Theta_{g, \tau }(\var d_Gx) &=&
\int_{G/G(s)}\int_{G(s)_{\e(s)}}\theta_{g,\tau,s}(y)\var(xysx^{-1})
|\det(1-(sy)^{-1})_{\g/\g(s)}|_pd_{G(s)}yd_{G/G(s)}\dot x \\
 &=& |\det(1-s^{-1})_{\g/\g(s)}|_p
\sum_{\omega \in \O_{g,\lambda_0}\cap \g^*(s) }\phi_{g, \tau,
s}(\omega)
\\&& \int_{G/G(s)} \int_{\omega}(1_{\g(s)_{\e(s)|n_F|_p}}\var_s^x \circ \exp
d_{\g(s)}X)\, \widehat{}_{\g(s)}  (l)
d\mu_{\omega}(l)d_{G/G(s)}\dot x.
 \end{eqnarray*}
$\bullet $ En premier lieu,  on suppose que 
 $\hat{s} \not\in {Z_g}$.

\begin{lem}\label{yfiniss}
 On a, pour toute $G(s)$-orbite $\omega$
 dans $\O_{g,\lambda_0}\cap \g^*(s)$,
 $$
 \sum_{\tau \in ({R_g}^\g)\, \widehat{}_{\psi_{\bar \lambda_0}.\psi}}
 \dim\tau\,
 \phi_{g, \tau, s}(\omega) = 0.
 $$
 \end{lem}
 Il r\'esulte du lemme ci-dessus que l'on a~:
    $$\sum_{\tau \in
   (R_g^\g)\, \widehat{}_{\psi_{\bar\lambda_0}.\psi}} \frac{\dim \tau}{|R_g^\g/Z_g|}
  \Theta_{g, \tau }(\var d_Gx)=0.$$
$\bullet$ En deuxi\`eme lieu, on suppose que $\hat{s} \in {Z_g}$.
Soit $ \omega \in \O_{g,\lambda_0}\cap \g^*(s)$ et $x \in G$ tel que
$x.(g+ \lambda_0) \in \omega$.  On peut supposer que $x = 1$. On a~:
 $$
 \phi_{g, \tau, s}(\omega) = c \, {\rm Tr}(\tau (\hat{s}))
 = c \, (\dim \tau) (\psi_{\bar\lambda_0}.\psi)(\hat{s}),
 $$
o\`u $c$ est un nombre complexe bien d\'etermin\'e (voir paragraphe \ref{hcgg} pour la d\'efinition de
$\phi_{g, \tau, s}$). Il en r\'esulte que
  $$
   \int_{({R_{g,\e}^\g}N)^\perp}
 \sum_{\tau \in (R_g^\g)\, \widehat{}_{\psi_{\bar \lambda_0}.\psi }}
 \dim \tau   \mbox{Tr}(\tau (\hat{s}) )   d_{\e}\psi =
     {|R_g^\g/Z_g|}
\psi_{\bar\lambda_0}(\hat{s})
    \int_{({R_{g,\e}^\g}N)^\perp} \psi(\hat{s}) d_\e\psi.
     $$
     \begin{itemize}
\item Premier cas~: on suppose que  $s \in R_{g,\e}$. Alors, on a : 
$$\int_{({R_{g,\e}^\g}N)^\perp} \psi(\hat{s}) d_\e\psi=1.$$ Il
s'en suit que
\begin{eqnarray*}
 && \sum_{\bar{\lambda'} \in
{(\j_{g,\e(s)|n_F|_p})}^\perp /{(\j_{g,\e|n_F|_p})}^\perp }
\int_{({R_{g,\e}^\g}N)^\perp}
 \sum_{\tau \in (R_g^\g)\, \widehat{}_{(\psi_{\bar {\lambda} + \bar{\lambda'}}.\psi )}}
 \dim \tau   \mbox{Tr}(\tau (\hat{s}) )   d_{\e}\psi \\
 &&\qquad \qquad \qquad \qquad \qquad \qquad \qquad \qquad=
     {|R_g^\g/Z_g|}
\psi_{\bar{\lambda}}(\hat{s})\sum_{\psi \in
({R_{g,\e}^\g}N/{R_{g,\e(s)}^\g}N)\, \widehat{}}\psi(\hat{s}).
\end{eqnarray*}
Compte tenu du choix de  $\e(s)$,  on a  $\hat{s} \notin {R_{g,\e(s)}^\g}N$. Ainsi,  le
deuxi\`eme membre de l'\'egalit\'e ci-dessus est \'egale \`a $0$.

\item Deuxi\`eme cas~:  on suppose que  $s \not \in R_{g,\e}$. Alors, il
existe
 $\psi_0 \in ({R_{g,\e}^\g}N)^\perp$ tel que $\psi_0(\hat{s}) \neq 1$.
 Par suite, on a :
  $$
    \int_{({R_{g,\e}^\g}N)^\perp} \psi(\hat{s})
    d_\e\psi=\psi_0(\hat{s}) \int_{({R_{g,\e}^\g}N)^\perp} \psi(\hat{s})
    d_\e\psi \,  = \,  0.
    $$
    \end{itemize}
Ainsi, la proposition   est d\'emontr\'ee.
\end{subsubsection}
\begin{subsubsection}{} Soit $ \var \in C_c^\infty(\W_G(s, \e(s)))$.  Il r\'esulte de la proposition \ref{yyfiniss} que l'on a : 
$$ 
 \int_{Y_G(g)}  \Theta_{g, \tau}(\var d_Gx)\, d_g\tau= 0.
 $$
 D'autre part, on a : $\W_G(s , \e(s) ) \cap U = \emptyset$. Si bien que l'on a : 
 $$
\int_{G.(g_{|\u})}\widehat{(\var_{|U}\circ \exp d_{\u}X)}_{\u}(l) d\mu_{G.(g_{|\u})}(l)=0. 
 $$
 \end{subsubsection}
 \end{subsection}
  \begin{subsection}{}
 Soit $S$ un syst\`eme de repr\'esentants des $G$-orbites semi-simples de $G$. 
 La proposition \ref{yAdam-Khemais} r\'esulte des r\'esultats \'etablis dans les paragraphes \ref{finish} et \ref{lem} et d'une partition de l'unit\'e subordonn\'ee \`a 
 $(\W_G(s , \e(s) ))_{s \in S}$.
 \end{subsection}
\begin{subsection}{D\'emonstration du lemme \ref{yfiniss}}
On note $\chi = \psi.\psi_{\bar\lambda_0}$.
D'apr\`es le th\'eor\`eme de Frobenius, on a~:
 $$
\oplus_{\tau \in (R_g^\g)\,  \widehat{}_\chi}  \, (\dim \tau) \,  \tau =\mbox{ Ind}_{Z_g}^{R_g^\g}\chi
.
 $$
Il en r\'esulte que l'on a, 
$$
\sum_{\tau \in (R_g^\g)\, \widehat{}_{\chi} }\dim \tau \, 
\mbox{Tr}(\tau (\hat{s})) =
 \mbox{ Tr} (({\mbox{ Ind}}_{Z_g}^{R_g^\g}\chi)(\hat{s})).
 $$
  La repr\'esentation  induite ${ \mbox{Ind}}_{Z_g}^{R_g^\g}\chi$ agit dans l'espace
 $\mathcal H_\chi$ des fonctions  continues $\var$ de $R_g^\g$
 \`a valeurs dans $\cit$ v\'erifiant : $$\var(xy) = \chi(y^{-1})\var(x) \mbox{, pour tout }
 x \in R_g^\g ,  y \in {Z_g}.$$
 On d\'esigne par $(y_1, y_2, \ldots, y_m)$
 un syst\`eme repr\'esentatif de $R_g^\g/{Z_g}$.
 On d\'efinit, pour chaque $i \in \{1, 2, \ldots, m\}$,
 une fonction $f_i$ par~:
 $$f_i(x) = 0 \mbox{ si } x \notin g_i{Z} \mbox{ et } f_i(g_i y) = \chi(y^{-1})
 \mbox{, pour tout } y \in {Z_g}.$$
 On v\'erifie que  $\{f_i,\,  i= 1, 2, \ldots, m\}$ est une base de
 $\mathcal H_\chi$. Si $f \in \mathcal H_\chi$,
 on note $f^{\hat{s}}$ la fonction d\'efinie par
 $$f^{\hat{s}}(x)= f({\hat{s}}^{-1}x) \mbox{, pour tout } x \in R_g^\g.$$
  Pour chaque $i \in \{1, 2, \ldots, m\}$, on note $j_i $ l'entier de  $\{1, 2, \ldots, m\}$
  tel que $\hat{s}y_i \in y_{j_i}{Z_g}$. On  pose $\hat{s}y_i  = y_{j_i}z_i$, $z_i \in Z_g$.
   Alors, on a~:
   $$f_i^{\hat{s}} = \chi(z_i) f_{j_i}.$$
   Compte tenu de l'hypoth\`ese  $\hat{s} \notin {Z_g}$,
   on a  $i \neq j_i$, pour tout $i \in \{1, 2, \ldots, m\}$. Si bien que l'on a~:
    $$\mbox{ Tr} (({\mbox{Ind}}_{Z_g}^{R_g^\g}\chi)(\hat{ s})) = 0.$$ D'o\`u le lemme.
 \end{subsection}
 \end{section}
\begin{section}{Exemple. }\label{exempl} On pose
$$G=\left\{(a, x, y)_G = \left( \begin{array}{ccc}
                       a  & \, \, \,0 & \, \, \,x  \\
                       0 & \, \, \,\, \, \, a^{-1} &\, \, \,  y \\
                       0 & \, \, \,0 &\, \, \, 1

                         \end{array}
                         \right)
, \, a  \in k^\times, \, x, y \in k\right\}.$$ 
Alors $G$ est l'ensemble des points rationnels d'un groupe alg\'ebrique d\'efini sur $k$.
On note $\g$ 
l'alg\`ebre de Lie de $G$ :  $$\g=\left\{(x, y, z)_\g =
\left( \begin{array}{ccc}
                       x& \, \, \,0 &\, \, \, y  \\
                       0 & \, \, \,-x &\, \, \, z \\
                       0 & \, \, \,0 & \, \, \,0

                         \end{array}
                         \right)
, \, x, y, z \in k\right\}.$$ On pose $E_1 = (1, 0, 0)_\g\, , E_2 = (0, 1, 0)_\g$, et $  E_3 = (0, 0,
1)_\g.$ Alors, $(E_1, E_2, E_3)$ est une base de $\g$ et on a~:
$$[E_1\, , E_2]= E_2 \, ,  [E_1\, , E_3]= -E_3\, , [E_2\, , E_3]=
0.$$
Si bien que $\g$ est  r\'esoluble. 
 La mesure
de Haar $d_G(a, x,y)_G = |\, a \, |_p^{-1}d\mu(a)d\mu(x)d\mu(y)$ de $G$ est invariante
par translations \`a gauche et \`a droite. Le groupe de Lie $G$
est unimodulaire.
\begin{subsection}{Orbites co-adjointes. }
On d\'esigne par $(E_1^*\, , E_2^* \, , E_3^*)$ la base duale
de $(E_1, E_2, E_3)$. Pour $(\alpha, \beta, \gamma) \in k^3$, on
note $f_{\alpha, \beta, \gamma}=\alpha E_1^* + \beta E_2^* +
\gamma E_3^*$ et $\O_{\alpha, \beta, \gamma}= G.f_{\alpha, \beta,
\gamma}$. On a~: $$\O_{\alpha, \beta, \gamma} = \{f_{\alpha- \beta
x + \gamma y , \beta a , \gamma a^{-1}} \, , a \in k^\times , x, y
\in k \}. $$ Les orbites co-adjointes ferm\'ees  dans $\g^*$ sont~:
\begin{enumerate}
\item les points $f_{\alpha, 0, 0}$\, , $\alpha \in k$ .
\item $\O_{0, 1, \gamma}$\, , $\gamma \in k^\times$.
\end{enumerate}
Les formes lin\'eaires  fortement r\'eguli\`eres sur $\g$  sont
$f_{\alpha, \beta, \gamma}$, avec $(\beta, \gamma) \neq (0, 0)$.
On a, pour $(\beta, \gamma) \neq (0, 0)$, $$G(f_{\alpha, \beta,
\gamma}) = \{(1, x, y)_G \, / \, \beta x -\gamma y = 0\, , \, x, y
\in k \}. $$ Si bien que l'on a $\j = 0$.
 On pose ainsi $$\Omega_G = \{f_{\alpha, \beta, \gamma}\, ,
\alpha \in k , \beta, \gamma \in k^\times \}.$$ Alors,  $\Omega_G$
est un ouvert de Zariski de $\g^*$, $G$-invariant.

\noindent D'autre part, le radical unipotent de $G$ est $U = \{(1,x,y)_G, \, \, x, y \in k \}$, d'alg\`ebre de Lie $$\u = \{(0, y, z)_\g, \, \, y, z \in k\}.$$  
 Pour $( \beta, \gamma) \in k^2$, on
note $u_{ (\beta, \gamma)}= \beta {E_{2|\u}^*} +
\gamma {E_{3|\u}^*}$.
Avec ces notations, on a : 
$$\Omega_U = \{u_{ (\beta, \gamma)}, \, \beta, \gamma \in k^\times\}.$$
Soit
$\Gamma_U$ le graphe de la relation d\'equivalence d\'efinie par
$G$ sur $\Omega_U$. On a~: $$\Gamma_U = \{(u_{(\beta,
\gamma)}, u_{( \beta a, \gamma a^{-1})}) \, , \beta, \gamma, a \in k^\times \}.$$ On voit que $\Gamma_U$
  co\"{\i}ncide avec l'image du plongement~:  $$
{k^\times}^3\longrightarrow \Omega_U \times \Omega_U  , ( \beta, \gamma, a) \longmapsto (u_{( \beta, \gamma)},
u_{(\beta a, \gamma a^{-1})}). $$  Il en r\'esulte que $G
\backslash  \Omega_U$ admet une structure de vari\'et\'e
$k$-analytique telle que la projection canonique de  $\Omega_U$
dans  $G \backslash  \Omega_U$ soit une submersion. 

\noindent Pour $\beta, \gamma \in k^\times$, on note  $\O_{ (\beta, \gamma)}= G.u_{(\beta,
\gamma)}$. On a~: $$\O_{(\beta, \gamma)} = \{u_{( \beta a , \gamma a^{-1})} \, , a \in k^\times  \}=\O_{(1, \beta \gamma)}.$$
 On en d\'eduit que l'application :
$$\psi : k^\times \longrightarrow G \backslash \Omega_U \, ,
\gamma \longmapsto \O_{( 1, \gamma)}, $$  est un diff\'eomophisme. 

\noindent Pour $\gamma \in k^\times$, on a  $\O_{(1,  \gamma)}$ est ferm\'ee dans $\u^*$ et   
$G(u_{ (1, \gamma)})= U. $
On munit $U$ de la mesure de Haar $d_U(1,x,y)_G= d\mu(x)d\mu(y)$, $\u$ de la mesure de Haar $d_\u (0, y, z)_\g = d\mu(y)d\mu(z)$, et ${\u^*}$ de la mesure de Haar $d_{\u^*}l$ duale de $d_\u (0, y, z)_\g$. Alors, on a : 
$$\int_{\O_{ (1, \gamma)}} \var(l) d\beta_{\O_{ (1, \gamma)}}(l) = \int_{k^\times}\var( u_{ (a, a^{-1}\gamma)})\frac{d\mu(a)}{|\, a \, |_p}.$$
 Ainsi, pour
toute fonction $\var $ continue \`a support compact dans $\u^*$, on a~:
$$
\int_{\u*}\var(l)d_{\u^*}l =
\int_{k^\times}\int_{\O_{(1,\gamma)}}\var(l)d\mu_{\O_{(1,\gamma)}}(l)
d\mu(\gamma).
$$
\end{subsection}
\begin{subsection}{Mesure de Plancherel de $G$.}
Soit $\gamma \in k^\times$.
On d\'esigne par $\chi_\gamma$ le caract\`ere
de $U $ d\'efini par :
$$\chi_\gamma((1, x, y)_G)= \varsigma(<f_{0,1, \gamma}, (0, x, y)_\g >)
=\varsigma(x + \gamma y), \, \, x, y \in k. $$ On note $\pi_\gamma =
\mbox{Ind}_{U}^G\chi_\gamma$.
Alors, $\pi_\gamma$ est une repr\'esentation  
  unitaire irr\'eductible de $G$.
 Elle se r\'ealise dans $
L^2(k^\times, |x|_p^{-1}d\mu(x))$, l'action de $G$ \'etant d\'efinie
par~: $$\pi_\gamma((a, x, y)_G) \var(b)=  \varsigma(b^{-1}x +b
\gamma y)\var(a^{-1}b).$$  
On a $\pi_\gamma$ est admissible. En effet, soit $K$
un sous-groupe compact ouvert de $G$. Il s'agit de d\'emontrer que $(L^2(k^\times, |x|_p^{-1}d\mu(x)))^K$ est de dimension finie.
Pour cela,  on note
$K_1 = \{a \in k^\times \, / \, (a, 0,0)_G \in K\}, \,  K_2 =
\{(x, y) \in k^2 \, / \, (1, x,y)_G \in K\}$. Alors, $K_1$ (resp.
$K_2$) est un sous-groupe compact ouvert de $k^\times$ (resp.
$k^2$).
 Soit $\var \in \mathcal
(L^2(k^\times, |x|_p^{-1}d\mu(x)))^K$ que l'on suppose non nulle. Si  $b \in k^\times$ tel que $\var(b) \neq 0$ alors
 on a $b^{-1}x +b \gamma y \in \ker \varsigma$, pour tout $(x, y) \in
K_2$. Ainsi, il existe un compact $M$ de $k^\times$ d\'ependant  de
$K_2$ tel que  le support de $\var$ soit
contenu dans $M$. On note $\bar M$ l'image de $M$ dans
$k^\times/K_1$ par la projection canonique de $k^\times$ dans
$k^\times/K_1$.  Alors  $\bar M$ est un ensemble fini et on a $\dim (L^2(k^\times,
|x|_p^{-1}d\mu(x)))^K \leq |\bar M| $.

\noindent On note $\Theta_\gamma$ le
caract\`ere de $\pi_\gamma$. Pour tout $\var \in C_c^\infty(G)$,
on a~:
$$
\var((1,0,0)_G)=\int_{k^\times}\Theta_\gamma(\var d_Gx) d\mu(\gamma).
$$

\end{subsection}
\end{section}
\section*{Remerciements :}
{\it Ma reconnaissance va au  Professeur  P. Torasso, dont
l'\'etendue des connaissances et la disponibilit\'e m'ont permis
de mener \`a bien ce pr\'esent travail. Je remercie  le
 Professeur A. Bouaziz pour les  discussions enrichissantes que nous avons eues.
 Je remercie  le Professeur M. Duflo pour l'int\'er\^et qu'il a montr\'e pour mon travail.
  Cet article a \'et\'e en partie r\'ealis\'e lors d'un s\'ejour au sein du Laboratoire de Math\'ematiques et Applications, UMR 6086 du CNRS et de l'Universit\'e de
 Poitiers, financ\'e par une bourse de Chercheur invit\'e de la 
R\'egion Poitou-Charentes et aussi  d'un  s\'ejour senior au sein du  m\^eme Laboratoire dans le cadre du PHC Curien $G1504$. 
 }

\addcontentsline{toc}{chapter}{Bibliographie}

 \footnotesize
  \renewcommand{\refname}{R\'{e}f\'{e}rences}

 \end{document}